\documentclass[
    Draft=false,
    Language=english,
    Complement=false,
    Contrast=0,
    BGr=1,BGg=1,BGb=1,
    HAr=0.4,HAg=0.5,HAb=0.6
]{ArticleB}

\usepackage{multicol}
\newcommand{\ForestG}{\mathfrak{g}}
\newcommand{\EasterlyWindLeq}{\Leq}
\newcommand{\EasterlyWindCovering}{\rightharpoonup}
\newcommand{\EasterlyWindCoveringRT}{\overset{*}{\EasterlyWindCovering}}
\newcommand{\EasterlyWindEdgeLabeling}{\lambda}
\newcommand{\Decoration}{\mathrm{dc}}
\newcommand{\Connection}{\mathrm{cnc}}
\newcommand{\Parent}{\mathrm{pa}}
\newcommand{\LocalPosition}{\mathrm{lp}}
\newcommand{\LeftBrothers}{\mathrm{lb}}
\newcommand{\SignatureExample}{\Signature_{\mathrm{e}}}
\newcommand{\NaturalHopfAlgebra}{\mathbf{N}}
\newcommand{\Reduced}{\mathrm{rd}}
\newcommand{\InternalNodes}{\mathrm{N}}
\newcommand{\Edges}{\mathrm{E}}
\newcommand{\Contraction}{\partial}
\newcommand{\SpecialNodes}{\mathcal{X}}
\newcommand{\Tilt}{\mathrm{tlt}}
\newcommand{\TiltReverse}{\Tilt^{\mathrm{r}}}
\newcommand{\InitialInterval}{\kern 0.125em \uparrow \kern -0.25em}
\newcommand{\TerminalInterval}{\kern 0.125em \downarrow \kern -0.25em}
\newcommand{\GeometricRealization}{\mathfrak{G}}
\newcommand{\CompositionFirst}{\bullet}
\newcommand{\SetLeaningForests}{\mathfrak{L}}
\newcommand{\Covered}{\lessdot}
\newcommand{\DecreasingTree}{\mathrm{dt}}
\newcommand{\DecreasingWord}{\mathrm{dw}}
\newcommand{\AdmissibleSet}{\mathrel{\vdash}}
\newcommand{\Top}{\Uparrow \kern -0.25em}
\newcommand{\Bottom}{\Downarrow \kern -0.25em}
\newcommand{\ToRootedTree}{\mathrm{rt}}
\newcommand{\ToFullyTiltedTerm}{\mathrm{ft}}
\newcommand{\ToBinaryTree}{\mathrm{bt}}
\newcommand{\ToForest}{\mathrm{if}}
\newcommand{\Scope}{\mathrm{sc}}
\newcommand{\TamariLeq}{\Leq_{\mathrm{T}}}
\newcommand{\BinaryTreesImage}{\mathfrak{B}}
\newcommand{\ForestMin}{\mathbf{f}^{\Uparrow}}
\newcommand{\ForestMax}{\mathbf{f}^{\Downarrow}}

\newcommand{\LeafSimple}{%
    \kern 0.0625em%
    \scalebox{.75}{%
        \begin{tikzpicture}[Centering,yscale=0.3]
            \node[Leaf]()at(0,0){};
        \end{tikzpicture}%
    }%
    \kern 0.0625em%
}

\newcommand{\Over}{%
    \kern 0.125em%
    \begin{tikzpicture}[Centering,scale=.18]
        \draw[thick](0,0)--(1,1.25);
    \end{tikzpicture}%
    \kern 0.125em%
}

\newcommand{\Under}{%
    \kern 0.125em%
    \begin{tikzpicture}[Centering,scale=.18]
        \draw[thick](0,0)--(1,-1.25);
    \end{tikzpicture}%
    \kern 0.125em%
}

\Title{Fundamental and homogeneous bases of Hopf algebras built from nonsymmetric operads}
\TitleShort{Bases of natural Hopf algebras of operads}
\Author{Samuele Giraudo}
\AuthorShort{S. Giraudo}
\Subjects{%
    05C05, 
    06A07, 
    16T30, 
    18M80. 
}
\Date{\today}
\Address{%
    Université du Québec à Montréal, LACIM, \\
    Pavillon Président-Kennedy, 201 Avenue du Président-Kennedy, Montréal, H2X~3Y7, Canada.
}
\Email{giraudo.samuele@uqam.ca}
\Funding{%
    This research has been partially supported by the Natural Sciences and Engineering
    Research Council of Canada (RGPIN-2024-04465).
}
\Abstract{%
    We introduce new partial order structures on the underlying sets of free nonsymmetric
    operads. These posets involve decorated ordered rooted trees, and their terminal
    intervals are lattices. These lattices are not graded, not self-dual, and not
    semi-distributive, but they are EL-shellable, and their Möbius functions take values in
    $\{-1, 0, 1\}$. They admit sublattices on the families of $m$-Fuss–Catalan objects and
    of forests of trees. This latter order structure is used to construct two new bases for
    the natural Hopf algebras of free nonsymmetric operads: a fundamental basis and a
    homogeneous basis. Along with the already known elementary basis of these Hopf algebras,
    this yields a triple of bases. The situation is similar to what is observed in the Hopf
    algebras of Malvenuto–Reutenauer, Loday–Ronco, and noncommutative symmetric functions,
    each of which presents such triples of bases and basis changes involving, respectively,
    the right weak partial order, the Tamari partial order, and the Boolean lattice partial
    order.
}

\begin{document}

\MakeFirstPage

\pagebreak

\section{Introduction} \label{sec:introduction}
Since the advent of modern Hopf algebra theory around the 1990s, most Hopf algebras defined
on the linear spans of sets $X$ of combinatorial objects are equipped with multiple bases.
This is the case, as main examples, for the Malvenuto-Reutenauer Hopf algebra
$\FQSym$~\cite{MR95,DHT02} of permutations, the Loday-Ronco Hopf algebra
$\PBT$~\cite{LR98,LR02,HNT05} of binary trees, and the Hopf algebra of noncommutative
symmetric functions $\NCSF$~\cite{GKLLRT95} of integer compositions. Each of these
structures is equipped with bases $\Bra{\BasisF_x}_{x \in X}$, $\Bra{\BasisE_x}_{x \in X}$,
and $\Bra{\BasisH_x}_{x \in X}$. A notable fact is that, in each case, there is a partial
order relation $\Leq$ on $X$ such that
\begin{equation} \label{equ:bases_f_e_h}
    \BasisE_x = \sum_{\substack{x' \in X \\ x \, \Leq \, x'}} \BasisF_{x'}
    \quad \text{and} \quad
    \BasisH_x = \sum_{\substack{x' \in X \\ x' \, \Leq \, x}} \BasisF_{x'},
\end{equation}
and two binary operations $\Over$ and $\Under$ on $X$ such that
\begin{equation} \label{equ:product_over_under_bases_e_h}
    \BasisE_{x_1} \Product \BasisE_{x_2} = \BasisE_{x_1 \Over x_2}
    \quad \text{and} \quad
    \BasisH_{x_1} \Product \BasisH_{x_2} = \BasisH_{x_1 \Under x_2},
\end{equation}
and
\begin{equation} \label{equ:product_basis_f}
    \BasisF_{x_1} \Product \BasisF_{x_2}
    = \sum_{\substack{x \in X \\ x_1 \Over x_2 \, \Leq \, x \, \Leq x_1 \, \Under x_2}}
    \BasisF_x,
\end{equation}
where $\Product$ is the product of the Hopf algebra. For $\FQSym$, $\Leq$ is the right weak
partial order, for $\PBT$, $\Leq$ is the Tamari partial order~\cite{Tam62}, and for $\NCSF$,
$\Leq$ is the Boolean lattice partial order. The case of $\NCSF$ is prototypical and, by
analogy with the theory of symmetric functions, the $\BasisF$-basis is termed the
\textit{fundamental} basis, while the $\BasisE$-basis (resp.\ $\BasisH$-basis) is termed the
\textit{elementary} (resp.\ \textit{homogeneous}) basis. There are a lot of other known
examples of Hopf algebras or associative algebras sharing these
properties~\cite{NT06,Gir12,CGM15,CP17,CG22}.

Besides all this, the \textit{natural Hopf algebra} of an operad $\Operad$ is a Hopf algebra
$\NaturalHopfAlgebra \App \Operad$ whose bases are indexed by some words on $\Operad$, where
the coproduct is inherited from the composition map of~$\Operad$. This construction is
considered for instance in~\cite{vdl04,Gir11,BG16,Gir24}, and a noncommutative variation for
nonsymmetric operads is introduced in~\cite{ML14} and already employed in~\cite{Gir11}.
However, surprisingly, no alternative bases are known in general for $\NaturalHopfAlgebra
\App \Operad$. In this work, we focus on natural Hopf algebras $\NaturalHopfAlgebra \App
\ULine{\SetTerms \App \Signature}$ of free nonsymmetric operads $\SetTerms \App \Signature$
generated by a signature $\Signature$. These Hopf algebras are defined on the linear span of
certain forests and come, through the construction $\NaturalHopfAlgebra$, with an
$\BasisE$-basis.

Our main contribution consists in the introduction of a general partial order, the
\textit{$\Signature$-easterly wind partial order}, defined on some treelike structures
decorated on a signature $\Signature$ (called \textit{$\Signature$-terms}). This poset
exhibits notable properties and specializes as a poset on forests as well as on various
other families of combinatorial objects. This structure leads to the definition of
fundamental and homogeneous bases satisfying~\eqref{equ:bases_f_e_h} for
$\NaturalHopfAlgebra \App \ULine{\SetTerms \App \Signature}$. We also introduce two binary
operations $\Over$ and $\Under$ on forests, such that the triple formed by the elementary,
fundamental, and homogeneous bases of $\NaturalHopfAlgebra \App \ULine{\SetTerms \App
\Signature}$ satisfies~\eqref{equ:product_over_under_bases_e_h}
and~\eqref{equ:product_basis_f}.

The contents and the results of this work are presented as follows.

Section~\ref{sec:easterly_wind_partial_orders}, after presenting preliminary notions about
signatures and $\Signature$-terms, introduces the easterly wind partial order relation
$\EasterlyWindLeq$ on the set of $\Signature$-terms. This partial order relation is defined
via \textit{connection words}, which are some sequences of rational numbers associated with
$\Signature$-terms. In parallel, we define a rewrite rule $\EasterlyWindCovering$ on
$\Signature$-terms, and establish, through Theorem~\ref{thm:easterly_wind_posets}, that
$\EasterlyWindCovering$ is the covering relation of $\EasterlyWindLeq$. We then introduce a
map $\Tilt \App \SpecialNodes$, called the \textit{$\SpecialNodes$-tilting map}, on the
easterly wind poset which turns out to be a closure operator~\cite{DP02} of this poset, as
stated by Theorem~\ref{thm:tilt_closure}. This will allow us to consider subposets of
easterly wind posets on closed terms w.r.t.\ $\Tilt \App \SpecialNodes$, called
\textit{tilted terms}.

In Section~\ref{sec:geometric_and_lattice_properties}, we prove via
Theorem~\ref{thm:el_shellability} that the easterly wind posets are
EL-shellable~\cite{Bjo80,BW96} and that their Möbius functions take values in the set
$\Bra{-1, 0, 1}$. Using connections words, we also propose geometric realizations of these
posets. We end this section with Theorem~\ref{thm:easterly_wind_lattices}, which shows that
any terminal interval $\TerminalInterval \App \TermT$ from a term $\TermT$ in the easterly
wind poset is a lattice. This property also holds for its subposets on tilted terms. As a
side remark, these lattices $\TerminalInterval \App \TermT$ are not always join
semi-distributive.

Section~\ref{sec:easterly_wind_lattices_forests} focuses on special cases of easterly wind
posets. We begin by defining a notion of \textit{forest} as a particular kind of term. To
each word $w$ on a signature $\Signature$, we associate an interval $\Han{\ForestMin \App w,
\ForestMax \App w}$ of the $\Signature$-easterly wind poset.
Theorem~\ref{thm:maximal_intervals} shows that such intervals are maximal. We then study
these maximal intervals to construct posets on forests that are in bijection with objects
from the Fuss-Catalan family. This is established in Theorem~\ref{thm:fuss_catalan_posets}.
These resulting posets, which are also lattices, are distinct from the already known posets
on this combinatorial family~\cite{BPR12,CG22}, and to our knowledge, have not appeared
before in the literature. Through a separate construction involving tilted terms, we realize
the Tamari poset~\cite{Tam62} as a maximal interval of a particular easterly wind poset of
tilted terms. Proposition~\ref{prop:rooted_tree_lattices} provides an explicit poset
isomorphism with the Tamari poset using the Knuth realization~\cite{Knu04} involving ordered
rooted trees and scope sequences. We conclude this section by introducing \textit{leaning
forests}, a specific subclass of forests (and thus, of terms), which form the bases of the
natural Hopf algebras of free nonsymmetric operads. We endow this set with two concatenation
operations $\Over$ and $\Under$, and define a shuffle operation $\cshuffle$ on leaning
forests. The properties and concepts established in this section are crucial for the final
one.

In the final part, Section~\ref{sec:natural_hopf_algebras}, we use the easterly wind partial
order on leaning forests to build a fundamental and a homogeneous basis of the natural Hopf
algebra of a free nonsymmetric operad. Theorem~\ref{thm:product_fundamental_basis} shows
that the product of two elements of the fundamental basis can be expressed as an interval of
the easterly wind poset, or equivalently as a shuffle of leaning forests.
Theorem~\ref{thm:product_homogeneous_basis} shows in a similar way that the product of two
elements of the homogenous basis expresses through the $\Under$ operation on leaning
forests.

We conclude in Section~\ref{sec:conclusion} with some open questions raised by this
work.

\paragraph{General notations and conventions}
All functions are written in curried form: given a function $f : A_1 \times \dots \times A_n
\to A$, we denote its application by $f \App a_1 \App \dots \App a_n$ rather than
$f\Par{a_1, \dots, a_n}$. Accordingly, the function type of $f$ is $A_1 \to \dots \to A_n
\to A$; the arrow $\to$ is taken to be right-associative. Rather than enclosing
sub-expressions in parentheses, we use underlining to distinguish these parts within
expressions\footnote{These two notational conventions are particularly useful when working
with treelike structures and operads, as they simplify the handling of compositions and
nested applications.}. For a statement $P$, the Iverson bracket $\Iverson{P}$ takes the
value $1$ if $P$ is true and $0$ otherwise. For two integers $i$ and $j$, $\Han{i, j}$
denotes the interval $\Bra{i, \dots, j}$, $[i]$ denotes the set $\Han{1, i}$ and $\HanL{n}$
denotes the set $\Han{0, i}$. For a set $A$, $A^*$ is the set of words on $A$. For $w \in
A^*$, the length of $w$ is $\Length \App w$. The only word of length $0$ is the empty word
$\epsilon$. For any $i \in [\Length \App w]$, $w \App i$ is the $i$-th letter of $w$. For $a
\in A$ and $n \in \N$, $a^n$ is the word of length $n$ such that $a^n \App i= a$ for all $i
\in [n]$. Given two words $w$ and $w'$, the concatenation of $w$ and $w'$ is denoted by $w
\Conc w'$. In graphical representations of Hasse diagrams of posets, the order relation
progresses from top to bottom.

\section{Easterly wind partial orders} \label{sec:easterly_wind_partial_orders}
This first section is devoted to introducing a partial order on the set of
$\Signature$-terms, which are treelike structures that realize the free nonsymmetric operad
on the generating set $\Signature$. We establish several properties of this partial order.

\subsection{Signatures and terms}
We begin with some preliminary definitions concerning particular treelike structures, called
terms, which are defined from a set of allowed decorations, called signatures.

\subsubsection{Signatures}
A \Def{signature} is a set $\Signature$ endowed with a map $\Arity : \Signature \to \N$. For
any $\GenS \in \Signature$, $\Arity \App \GenS$ is the \Def{arity} of $\GenS$. For any $n
\in \N$, let $\Signature \App n := \Bra{\GenS \in \Signature : \Arity \App \GenS = n}$. For
the examples that will follow, we shall consider the signature $\SignatureExample :=
\Bra{\GenA_{i, j} : i, j \in \N}$ where for any $\GenA_{i, j} \in \SignatureExample$,
$\Arity \App \GenA_{i, j} = i$. To lighten the notations, we shall write simply $\GenA_i$
for $\GenA_{i, 0}$.

In what follows, we define several concepts ``C'' parameterized by a signature $\Signature$,
denoted by ``$\Signature$-C''. To streamline the phrasing, whenever there is no ambiguity,
we simply write ``C''.

\subsubsection{Terms}
Given a signature $\Signature$, an \Def{$\Signature$-term} is either the \Def{leaf}
$\LeafSimple$ or a pair $\Par{\GenS, \Par{\TermT_1, \dots, \TermT_{\Arity \App \GenS}}}$
where $\GenS \in \Signature$, and $\TermT_1$, \dots, $\TermT_{\Arity \App \GenS}$ are
$\Signature$-terms. For brevity, we write $\GenS \, \TermT_1 \dots \TermT_{\Arity \App
\GenS}$ for $\Par{\GenS, \Par{\TermT_1, \dots, \TermT_{\Arity \App \GenS}}}$. By
definition, an $\Signature$-term is therefore a decorated ordered rooted tree where each
internal node having $n$ children is decorated on $\Signature \App n$. The set of
$\Signature$-terms is denoted by $\SetTerms \App \Signature$. For instance,
\begin{math}
    \GenA_3
    \ULine{\GenA_1 \LeafSimple}
    \ULine{\GenA_2 \GenA_0 \LeafSimple}
    \ULine{\GenA_3 \GenA_0 \ULine{\GenA_1 \LeafSimple} \LeafSimple}
\end{math}
is an $\SignatureExample$-term and it writes as the decorated ordered rooted tree
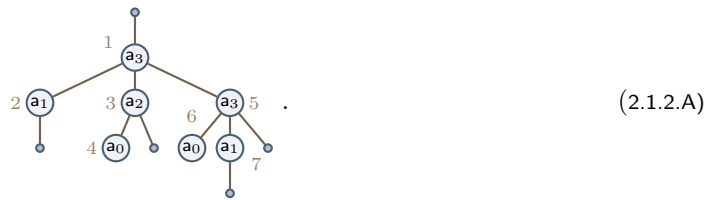
\begin{equation} \label{equ:example_term}
    \scalebox{1}{
        \begin{tikzpicture}[Centering,xscale=0.5,yscale=0.60]
            \node[Leaf](1)at(2.5,0){};
            \node[Node](2)at(2.5,-1){$\GenA_3$};
            \node[Node](3)at(0,-2){$\GenA_1$};
            \node[Leaf](4)at(0,-3){};
            \node[Node](5)at(2.5,-2){$\GenA_2$};
            \node[Node](6)at(2,-3){$\GenA_0$};
            \node[Leaf](7)at(3,-3){};
            \node[Node](8)at(5,-2){$\GenA_3$};
            \node[Node](9)at(4,-3){$\GenA_0$};
            \node[Node](10)at(5,-3){$\GenA_1$};
            \node[Leaf](11)at(5,-4){}; \node[Leaf](12)at(6,-3){};
            \draw[Edge](1)--(2);
            \draw[Edge](2)--(3);
            \draw[Edge](2)--(5);
            \draw[Edge](2)--(8);
            \draw[Edge](3)--(4);
            \draw[Edge](5)--(6);
            \draw[Edge](5)--(7);
            \draw[Edge](8)--(9);
            \draw[Edge](8)--(10);
            \draw[Edge](8)--(12);
            \draw[Edge](10)--(11);
            \node[above of=2,node distance=6pt,font=\scriptsize,text=ColB]{$1\kern 2.5em$};
            \node[left of=3,node distance=9pt,font=\scriptsize,text=ColB]{$2$};
            \node[left of=5,node distance=9pt,font=\scriptsize,text=ColB]{$3$};
            \node[left of=6,node distance=9pt,font=\scriptsize,text=ColB]{$4$};
            \node[right of=8,node distance=9pt,font=\scriptsize,text=ColB]{$5$};
            \node[above of=9,node distance=12pt,font=\scriptsize,text=ColB]{$6$};
            \node[below of=10,node distance=6pt,font=\scriptsize,text=ColB]{$\kern2.5em 7$};
        \end{tikzpicture}
    }.
\end{equation}
A \Def{subterm} of an $\Signature$-term $\TermT := \GenS \, \TermT_1 \dots \TermT_{\Arity
\App \GenS}$ is either $\TermT$ itself, or recursively a subterm of $\TermT_i$ where $i \in
[\Arity \App \GenS]$. For any $\GenS \in \Signature$, the \Def{$\GenS$-corolla} is the
$\Signature$-term $\Corolla \App \GenS := \GenS \, \LeafSimple \dots \LeafSimple$. In other
words, $\Corolla \App \GenS$ is the $\Signature$-term consisting in one single internal node
decorated by $\GenS$ and in $\Arity \App \GenS$ leaves. Let us now introduce some additional
definitions about $\Signature$-terms. Below, $\TermT$ is an $\Signature$-term.

The \Def{preorder traversal} of $\TermT$ is defined recursively as follows. If $\TermT =
\LeafSimple$, then the leaf forming $\TermT$ is visited. Otherwise, we have $\TermT = \GenS
\, \TermT_1 \TermT_2 \dots \TermT_{\Arity \App \GenS}$ where $\GenS \in \Signature$ and
$\TermT_1$, $\TermT_2$, \dots, and $\TermT_{\Arity \App \GenS}$ are $\Signature$-terms. In
this case, the root of $\TermT$ is visited first, and then, $\TermT_1$, $\TermT_2$, \dots,
and $\TermT_{\Arity \App \GenS}$ are visited from left to right according to their
respective preorder traversals. This procedure induces a total order on the leaves and
internal nodes of $\TermT$ where the first visited element is the smallest one. Note that
the symbolic notation of an $\Signature$-term already lists its leaves and internal nodes in
this order (see~\eqref{equ:example_term} and its symbolic notation given just before).

The \Def{arity} $\Arity \App \TermT$ of $\TermT$ is the number of occurrences of leaves of
$\TermT$. The leaves of $\TermT$ are numbered consecutively, starting with $1$, according to
their positions in the preorder traversal of $\TermT$. Given $i \in [\Arity \App \TermT]$,
the $i$-th leaf of $\TermT$ is \Def{extreme} if all internal nodes of $\TermT$ are visited
before the $i$-th leaf in the preorder traversal of $\TermT$.

The \Def{degree} $\Deg \App \TermT$ of $\TermT$ is the number of internal nodes of $\TermT$.
The internal nodes of $\TermT$ are numbered consecutively, starting with $1$, according to
their positions in the preorder traversal of $\TermT$. Henceforth, we identify each internal
node of $\TermT$ with the index assigned to it. When $\Deg \App \TermT \geq 1$, the
\Def{contraction} of $\TermT$ is the $\Signature$-term $\Contraction \App \TermT$ obtained
by replacing the last internal node $\Deg \App \TermT$ of $\TermT$ by a leaf. We denote by
$\InternalNodes \App \TermT$ the set $[\Deg \App \TermT]$ of internal nodes of $\TermT$. The
\Def{decoration word} of $\TermT$ is the word $\Decoration \App \TermT$ on $\Signature$ such
that for any $i \in \InternalNodes \App \TermT$, $\Decoration \App \TermT \App i$ is the
decoration of $i$ in $\TermT$. By a slight abuse of notation, let us denote by $\Arity \App
\TermT \App i$ the arity $\Arity \App \ULine{\Decoration \App \TermT \App i}$ of the
decoration of $i \in \InternalNodes \App \TermT$.

An \Def{edge} of $\TermT$ is a triple $\Par{i_1, j, i_2}$ such that $i_1, i_2 \in
\InternalNodes \App \TermT$ and $i_2$ is the $j$-th child of $i_1$, where the children of
$i_1$ are numbered from left to right, starting by $1$. For convenience, when $\Deg \App
\TermT \geq 1$, we consider that $\Par{1, 0, 1}$ is an edge of $\TermT$. This edge can be
seen as a loop on the root of $\TermT$. We denote by $\Edges \App \TermT$ the set of edges
of $\TermT$. For any internal node $i_2$ of $\TermT$, there is a unique edge of $\TermT$ of
the form $\Par{i_1, j, i_2}$, called the \Def{parent edge} of $i_2$. Under these
definitions, the \Def{parent} $\Parent \App \TermT \App i_2$ of $i_2$ is the node $i_1$, and
the \Def{local position} $\LocalPosition \App \TermT \App i_2$ of $i_2$ is the integer $j$.
With the previous convention, the internal node $1$ is the parent of itself and its local
position is $0$.

Let us give some examples of the previous definitions. The integers near each internal node
of the $\SignatureExample$-term $\TermT$ in~\eqref{equ:example_term} are the integers with
which they are identified. Moreover, we have $\Deg \App \TermT = 7$, $\Arity \App \TermT =
4$,
\begin{math}
    \Decoration \App \TermT = \GenA_3 \GenA_1 \GenA_2 \GenA_0 \GenA_3 \GenA_0 \GenA_1,
\end{math}
\begin{equation}
    \Edges \App \TermT
    = \Bra{
        \Par{1, 0, 1},
        \Par{1, 1, 2}, \Par{1, 2, 3}, \Par{3, 1, 4}, \Par{1, 3, 5}, \Par{5, 1, 6},
        \Par{5, 2, 7}
    },
\end{equation}
and the extreme leaves of $\TermT$ are the $3$-rd and the $4$-th ones. Besides, the
contraction of $\TermT$ is the $\SignatureExample$-term
\begin{math}
    \GenA_3
    \ULine{\GenA_1 \LeafSimple}
    \ULine{\GenA_2 \GenA_0 \LeafSimple}
    \ULine{\GenA_3 \GenA_0 \LeafSimple \LeafSimple}.
\end{math}

\subsection{Posets on terms} \label{subsec:posets_on_termes}
The purpose of this section is to define a partial order relation $\EasterlyWindLeq$ on the
set of $\Signature$-terms. This partial order is defined by comparing, letter by letter,
certain sequences of rational numbers obtained from $\Signature$-terms, called connection
words. We show that $\EasterlyWindLeq$ admits, as its covering relation, a rewrite rule
$\EasterlyWindCovering$ on the set of $\Signature$-terms, which consists of pruning and
grafting subterms in an appropriate way. We begin this section by introducing this rewrite
rule.

\subsubsection{A rewrite rule on terms} \label{subsubsec:rewrite_rule}
For any $i \in \N \setminus \{0\}$, let $\EasterlyWindCovering_i$ be the binary relation on
$\SetTerms \App \Signature$ defined as follows. Let $\TermT_1$ be an $\Signature$-term of
degree at least $i$ and such that its internal node $i$ is visited immediately after a leaf
in the preorder traversal of $\TermT_1$. Then, $\TermT_1 \EasterlyWindCovering_i \TermT_2$
holds if $\TermT_2$ is the $\Signature$-term obtained from $\TermT_1$ by moving to this leaf
the subterm rooted at $i$. For instance, we have
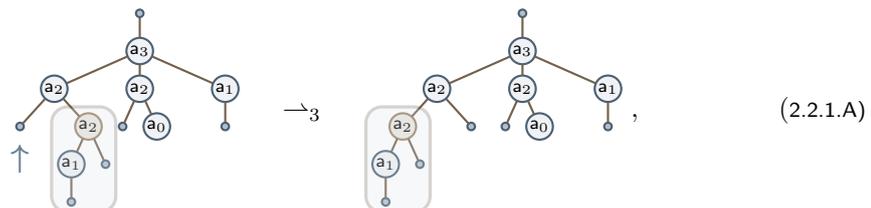
\begin{align} \label{equ:example_covering_1}
    \scalebox{1}{
        \begin{tikzpicture}[Centering,xscale=0.45,yscale=0.5]
            \node[Leaf](1)at(3.5,0){};
            \node[Node](2)at(3.5,-1){$\GenA_3$};
            \node[Node](3)at(1,-2){$\GenA_2$};
            \node[Leaf](4)at(0,-3){};
            \node[Node,MarkB](5)at(2,-3){$\GenA_2$};
            \node[Node](6)at(1.5,-4){$\GenA_1$};
            \node[Leaf](7)at(1.5,-5){};
            \node[Leaf](8)at(2.5,-4){};
            \node[Node](9)at(3.5,-2){$\GenA_2$};
            \node[Leaf](10)at(3,-3){};
            \node[Node](11)at(4,-3){$\GenA_0$};
            \node[Node](12)at(6,-2){$\GenA_1$};
            \node[Leaf](13)at(6,-3){};
            \draw[Edge](1)--(2);
            \draw[Edge](2)--(3);
            \draw[Edge](2)--(9);
            \draw[Edge](2)--(12);
            \draw[Edge](3)--(4);
            \draw[Edge](3)--(5);
            \draw[Edge](5)--(6);
            \draw[Edge](5)--(8);
            \draw[Edge](6)--(7);
            \draw[Edge](9)--(10);
            \draw[Edge](9)--(11);
            \draw[Edge](12)--(13);
            \node[FitNodes,fit=(5)(6)(7)(8)]{};
            \node[below of=4,node distance=12pt,font=\large,text=ColA]{$\bm{\uparrow}$};
        \end{tikzpicture}
    }
    & \quad \EasterlyWindCovering_3 \quad
    \scalebox{1}{
        \begin{tikzpicture}[Centering,xscale=0.45,yscale=0.5]
            \node[Leaf](1)at(3.5,0){};
            \node[Node](2)at(3.5,-1){$\GenA_3$};
            \node[Node](3)at(1,-2){$\GenA_2$};
            \node[Node,MarkB](4)at(0,-3){$\GenA_2$};
            \node[Node](5)at(-0.5,-4){$\GenA_1$};
            \node[Leaf](6)at(-0.5,-5){};
            \node[Leaf](7)at(0.5,-4){};
            \node[Leaf](8)at(2,-3){};
            \node[Node](9)at(3.5,-2){$\GenA_2$};
            \node[Leaf](10)at(3,-3){};
            \node[Node](11)at(4,-3){$\GenA_0$};
            \node[Node](12)at(6,-2){$\GenA_1$};
            \node[Leaf](13)at(6,-3){};
            \draw[Edge](1)--(2);
            \draw[Edge](2)--(3);
            \draw[Edge](2)--(9);
            \draw[Edge](2)--(12);
            \draw[Edge](3)--(4);
            \draw[Edge](3)--(8);
            \draw[Edge](4)--(5);
            \draw[Edge](4)--(7);
            \draw[Edge](5)--(6);
            \draw[Edge](9)--(10);
            \draw[Edge](9)--(11);
            \draw[Edge](12)--(13);
            \node[FitNodes,fit=(4)(5)(6)(7)]{};
        \end{tikzpicture}
    },
    \\
    \GenA_3
    \ULine{
        \GenA_2
        \LeafSimple
        \bm{\ColB{\ULine{\GenA_2 \ULine{\GenA_1 \LeafSimple} \LeafSimple}}}
    }
    \ULine{\GenA_2 \LeafSimple \GenA_0}
    \ULine{\GenA_1 \LeafSimple}
    & \quad \EasterlyWindCovering_3 \quad
    \GenA_3
    \ULine{
        \GenA_2
        \bm{\ColB{\ULine{\GenA_2 \ULine{\GenA_1 \LeafSimple} \LeafSimple}}}
        \LeafSimple
    }
    \ULine{\GenA_2 \LeafSimple \GenA_0}
    \ULine{\GenA_1 \LeafSimple}
    \notag
\end{align}
and
\begin{align} \label{equ:example_covering_2}
    \scalebox{1}{
        \begin{tikzpicture}[Centering,xscale=0.45,yscale=0.5]
            \node[Leaf](1)at(3.5,0){};
            \node[Node](2)at(3.5,-1){$\GenA_3$};
            \node[Node](3)at(1,-2){$\GenA_2$};
            \node[Leaf](4)at(0,-3){};
            \node[Node](5)at(2,-3){$\GenA_2$};
            \node[Node](6)at(1.5,-4){$\GenA_1$};
            \node[Leaf](7)at(1.5,-5){};
            \node[Leaf](8)at(2.5,-4){};
            \node[Node,MarkB](9)at(3.5,-2){$\GenA_2$};
            \node[Leaf](10)at(3,-3){};
            \node[Node](11)at(4,-3){$\GenA_0$};
            \node[Node](12)at(6,-2){$\GenA_1$};
            \node[Leaf](13)at(6,-3){};
            \draw[Edge](1)--(2);
            \draw[Edge](2)--(3);
            \draw[Edge](2)--(9);
            \draw[Edge](2)--(12);
            \draw[Edge](3)--(4);
            \draw[Edge](3)--(5);
            \draw[Edge](5)--(6);
            \draw[Edge](5)--(8);
            \draw[Edge](6)--(7);
            \draw[Edge](9)--(10);
            \draw[Edge](9)--(11);
            \draw[Edge](12)--(13);
            \node[FitNodes,fit=(9)(10)(11)]{};
            \node[below of=8,node distance=12pt,font=\large,text=ColA]{$\bm{\uparrow}$};
        \end{tikzpicture}
    }
    & \quad \EasterlyWindCovering_5 \quad
    \scalebox{1}{
        \begin{tikzpicture}[Centering,xscale=0.45,yscale=0.5]
            \node[Leaf](1)at(3,0){};
            \node[Node](2)at(3,-1){$\GenA_3$};
            \node[Node](3)at(1,-2){$\GenA_2$};
            \node[Leaf](4)at(0,-3){};
            \node[Node](5)at(2,-3){$\GenA_2$};
            \node[Node](6)at(1,-4){$\GenA_1$};
            \node[Leaf](7)at(1,-5){};
            \node[Node,MarkB](8)at(3,-4){$\GenA_2$};
            \node[Leaf](9)at(2.5,-5){};
            \node[Node](10)at(3.5,-5){$\GenA_0$};
            \node[Leaf](11)at(3,-2){};
            \node[Node](12)at(5,-2){$\GenA_1$};
            \node[Leaf](13)at(5,-3){};
            \draw[Edge](1)--(2);
            \draw[Edge](2)--(3);
            \draw[Edge](2)--(11);
            \draw[Edge](2)--(12);
            \draw[Edge](3)--(4);
            \draw[Edge](3)--(5);
            \draw[Edge](5)--(6);
            \draw[Edge](5)--(8);
            \draw[Edge](6)--(7);
            \draw[Edge](8)--(9);
            \draw[Edge](8)--(10);
            \draw[Edge](12)--(13);
            \node[FitNodes,fit=(8)(9)(10)]{};
        \end{tikzpicture}
    },
    \\
    \GenA_3
    \ULine{\GenA_2 \LeafSimple \ULine{\GenA_2 \ULine{\GenA_1 \LeafSimple} \LeafSimple}}
    \bm{\ColB{\ULine{\GenA_2 \LeafSimple \GenA_0}}}
    \ULine{\GenA_1 \LeafSimple}
    & \quad \EasterlyWindCovering_3 \quad
    \GenA_3
    \ULine{
        \GenA_2
        \LeafSimple
        \ULine{\GenA_2 \ULine{\GenA_1 \LeafSimple}
        \bm{\ColB{\ULine{\GenA_2 \LeafSimple \GenA_0}}}}
    }
    \LeafSimple
    \ULine{\GenA_1 \LeafSimple}.
    \notag
\end{align}
In these examples, in accordance with the previous definition of $\EasterlyWindCovering_i$,
each surrounded area shows the moved subterm rooted at an internal node $i$, and each arrow
shows the leaf which is visited immediately before $i$ in the preorder traversal, target of
the moved subterm. Observe that these transformations become particularly transparent when
the terms are displayed in their symbolic notation: the subterm rooted at $i$, written in
bold, is relocated on the leaf appearing immediately to its left. Remark that by setting
$\TermT_1$ as the $\SignatureExample$-term appearing in the left-hand sides
in~\eqref{equ:example_covering_1} and~\eqref{equ:example_covering_2}, there is no
$\SignatureExample$-term $\TermT_2$ such that $\TermT_1 \EasterlyWindCovering_1 \TermT_2$
nor $\TermT_1 \EasterlyWindCovering_2 \TermT_2$ because there are no leaves which are
visited before the internal nodes $1$ and $2$ in the preorder traversal of $\TermT_1$.
Moreover, there is no $\SignatureExample$-term $\TermT_2$ such that $\TermT_1
\EasterlyWindCovering_7 \TermT_2$ because an internal node, in this case the internal node
$6$, is visited just before visiting the internal node $7$ in the preorder traversal
of~$\TermT_1$.

From now on, except at some places in Section~\ref{subsubsec:fuss_catalan_lattices}, we will
exclusively use the symbolic notation for $\Signature$-terms. However, readers are invited
to convert these into their graphical notation if they feel more comfortable doing so.

It follows immediately from the definition of $\EasterlyWindCovering_i$ that
$\EasterlyWindCovering_1$ is the empty relation. Another immediate property is that for any
$i \in \N \setminus \{0\}$ and any $\Signature$-term $\TermT_1$, there is at most one
$\Signature$-term $\TermT_2$ such that $\TermT_1 \EasterlyWindCovering_i \TermT_2$. Let us
also denote by $\EasterlyWindCoveringRT_i$ the reflexive and transitive closure of
$\EasterlyWindCovering_i$.

Given two triples $x_1 := \Par{i_1, j_1, i}$ and $x_2 := \Par{i_2, j_2, i}$ of integers,
$x_1$ is \Def{dominated} by $x_2$ if the pair $\Par{i_1, -j_1}$ is lexicographically smaller
than or equal to the pair $\Par{i_2, -j_2}$ For instance, $\Par{2, 5, 4}$ is dominated by
$\Par{3, 7, 4}$ and by $\Par{2, 3, 4}$, but not by $\Par{1, 1, 4}$ neither by~$\Par{2, 6,
4}$.

The following lemma is a crucial tool that is used in several subsequent proofs.

\begin{Statement}{Lemma}{lem:covering_edge_properties}
    Let $\Signature$ be a signature and $\TermT_1$ and $\TermT_2$ be two $\Signature$-terms
    such that $\TermT_1 \EasterlyWindCovering_i \TermT_2$ for an $i \geq 2$. The following
    properties hold:
    \begin{enumerate}[label=(\roman*)]
        \item \label{item:covering_edge_properties_1}
        the decoration words of $\TermT_1$ and $\TermT_2$ are the same;
        \item \label{item:covering_edge_properties_2}
        the $\Signature$-terms $\TermT_1$ and $\TermT_2$ share all edges except for the
        parent edge of $i$;
        \item \label{item:covering_edge_properties_3}
        the parent edge of $i$ in $\TermT_1$ is dominated by the parent edge of $i$ in
        $\TermT_2$;
        \item \label{item:covering_edge_properties_4}
        for any internal nodes $i'$ and $i''$ of $\TermT_1$ and $\TermT_2$, if $i''$ is a
        descendant of $i'$ in $\TermT_1$, then $i''$ is also a descendant of $i'$ in
        $\TermT_2$.
    \end{enumerate}
\end{Statement}
\begin{Proof}
    Let $i_1 := \Parent \App \TermT_1 \App i$, $k$ be the leaf which is visited immediately
    before $i$ in the preorder traversal of $\TermT_1$, and $i_2$ be the internal node of
    $\TermT_1$ to which $k$ is attached.

    First, by definition of $\EasterlyWindCovering_i$, since $\TermT_2$ is obtained from
    $\TermT_1$ by moving on $k$ the subterm rooted at $i$, all internal nodes of $\TermT_1$
    and $\TermT_2$ are visited in the same order. Therefore, $\Decoration \App \TermT_1 =
    \Decoration \App \TermT_2$, establishing~\ref{item:covering_edge_properties_1}.

    Point~\ref{item:covering_edge_properties_2} is immediate since by definition of
    $\EasterlyWindCovering_i$, $\TermT_1$ and $\TermT_2$ differ only by the parent of their
    internal node~$i$, which is $i_1$ and $\TermT_1$ and $i_2$ in $\TermT_2$. All other
    edges are identical in both $\TermT_1$ and $\TermT_2$.

    To prove~\ref{item:covering_edge_properties_3}, we need to distinguish two cases. First,
    if both $k$ and $i$ admits the same parent $i_1 = i_2$, then, by definition of
    $\EasterlyWindCovering_i$, $\Par{i_1, j_1 - 1, i}$ is the parent edge of $i$ in
    $\TermT_2$. Therefore, the parent edge of $i$ in $\TermT_1$ is dominated by the parent
    edge of $i$ in $\TermT_2$. Otherwise, due to the fact that in $\TermT_1$, $k$ is
    attached to $i_2$, $i_1$ is the parent of $i$, and $k$ is visited immediately before $i$
    in the preorder traversal, it follows that $i_2$ is visited after $i_1$ in the preorder
    traversal of $\TermT_1$. For this reason, $i_1 < i_2$. Since $i_2$ is the parent of $i$
    in $\TermT_2$, the parent edge of $i$ in $\TermT_1$ is dominated by the parent edge of
    $i$ in $\TermT_2$.

    Finally, by hypothesis, $i$ is a descendant of $i_1$ in $\TermT_1$, and
    by~\ref{item:covering_edge_properties_2}, $i$ is a descendant of $i_2$ in $\TermT_2$.
    Moreover, by~\ref{item:covering_edge_properties_3}, $i_1 \leq i_2$. Hence, and since
    $i_2 < i$, we have either that $i_1 = i_2$ or that $i_2$ is a descendant of $i_1$.
    Therefore, in $\TermT_2$, $i$ remains a descendant of $i_1$. This property implies that,
    independently from the location of two internal nodes $i'$ and $i''$ of $\TermT_1$ and
    $\TermT_2$, if $i''$ is a descendant of $i'$ in $\TermT_1$, these two internal nodes
    enjoy the same property in~$\TermT_2$.
\end{Proof}

To illustrate Lemma~\ref{lem:covering_edge_properties}, observe that
in~\eqref{equ:example_covering_1}, the edge $\Par{2, 2, 3}$ is replaced by the dominating
edge $\Par{2, 1, 3}$, and in~\eqref{equ:example_covering_2}, the edge $\Par{1, 2, 5}$ is
replaced by the dominating edge $\Par{3, 2, 5}$. Moreover, all $\SignatureExample$-terms of
these two examples have $\GenA_3 \GenA_2 \GenA_2 \GenA_1 \GenA_2 \GenA_0 \GenA_1$ as
decoration word.

Let us denote by $\EasterlyWindCovering$ the binary relation on $\SetTerms \App \Signature$
defined as the union of all relations $\EasterlyWindCovering_i$ for $i \in \N \setminus
\{0\}$. We call $\EasterlyWindCovering$ the \Def{$\Signature$-easterly wind rewrite rule}.
Let us also denote by $\EasterlyWindCoveringRT$ the reflexive and transitive closure
of~$\EasterlyWindCovering$.

\subsubsection{Connection words} \label{subsubsec:connection_words}
The \Def{connection word} of an $\Signature$-term $\TermT$ is the word $\Connection \App
\TermT$ on $\Q$ of length $\Deg \App \TermT$ such that for any $i_2 \in \InternalNodes \App
\TermT$,
\begin{equation} \label{equ:connection}
    \Connection \App \TermT \App i
    = \Parent \App \TermT \App i + 1
    - 2^{\mbox{\small $\LocalPosition \App \TermT \App i
    - \Arity \App \TermT \App \ULine{\Parent \App \TermT \App i}$}}.
\end{equation}
In other words, $\Connection \App \TermT \App i$ expresses in binary fixed-point notation as
$\mathtt{x} \mathtt{.} \mathtt{1}^\ell$, where $\mathtt{x}$ is the binary expansion of
$\Parent \App \TermT \App i$ and $\ell$ is the number of internal nodes of $\TermT$ among
the siblings of $i$ that lie to its right.

For instance, the connection words of the $\SignatureExample$-terms $\TermT_1$ and
$\TermT_2$ of~\eqref{equ:example_covering_2} satisfy
\begin{equation}
    \Connection \App \TermT_1 =
    \Par{
        \ColA{\frac{15}{8}}, \ColA{\frac{7}{4}}, \ColA{2}, \ColA{\frac{7}{2}},
        \ColA{\frac{3}{2}}, \ColA{5}, \ColA{1}
    }
    \quad \text{and} \quad
    \Connection \App \TermT_2 =
    \Par{
        \ColA{\frac{15}{8}}, \ColA{\frac{7}{4}}, \ColA{2}, \ColA{\frac{7}{2}},
        \bm{\ColB{3}}, \ColA{5}, \ColA{1}
    }.
\end{equation}
Under the alternative interpretation of $\Connection \App \TermT$, the binary fixed-point
representation of $\Connection \App \TermT_1 \App 2$ is $\mathtt{1} \mathtt{.} \mathtt{1}^2
= \mathtt{1} \mathtt{.} \mathtt{11}$, which denotes the value $\frac{7}{4}$, and the one of
$\Connection \App \TermT_1 \App 6$ is $\mathtt{101} \mathtt{.} \mathtt{1}^0 = \mathtt{101}
\mathtt{.} \mathtt{0}$, which denotes the value $5$ as expected.

\begin{Statement}{Lemma}{lem:connection_words_interval}
    Let $\Signature$ be a signature, $\TermT$ be an $\Signature$-term, and $i$ be an
    internal node of $\TermT$. The parent of $i$ in $\TermT$ is the unique integer $i'$ such
    that $i' \leq \Connection \App i < i' + 1$.
\end{Statement}
\begin{Proof}
    Let $\Par{i', j, i}$ be the parent edge of $i$ in $\TermT$. From the
    definition~\eqref{equ:connection} of $\Connection$, the rational number $\Connection
    \App i$ is minimal when $j = \Arity \App \TermT \App i'$ and is, in this case, equal to
    $i'$. Moreover, $\Connection \App i$ is maximal when $j = 0$ and is, in this case, equal
    to $i' + 1 - 2^{\mbox{\small $- \Arity \App \TermT \App i$}}$. Since $2^{\mbox{\small $-
    \Arity \App \TermT \App i$}}$ is a positive number, the previous quantity is smaller
    than $i' + 1$. The stated property follows.
\end{Proof}

\begin{Statement}{Proposition}{prop:connection_map_injective}
    For any signature $\Signature$ and any word $w$ on $\Signature$, the map $\Connection$
    on the domain of the $\Signature$-terms having $w$ as decoration word, is injective.
\end{Statement}
\begin{Proof}
    Let $c$ be a word of rational numbers belonging to the image of the map $\Connection$ on
    the domain of $\Signature$-terms having $w$ as decoration word. Let us show that there
    exists a unique antecedent $\TermT$ of $c$ by $\Connection$. By denoting by $n$ the
    length of $w$, by Lemma~\ref{lem:connection_words_interval}, for any $i \in [n]$, there
    is a unique $i' \in [n]$ such that $i' \leq c \App i < i' + 1$. Therefore, the parent of
    the internal node $i$ of $\TermT$ is $i'$. Moreover, from~\eqref{equ:connection}, we
    have
    \begin{math}
        \LocalPosition \App \TermT \App i
        = \Arity \App \ULine{w \App i'} + \log_2 \App \ULine{i' + 1 - c \App i}.
    \end{math}
    This shows that the parent edge of $i$ in $\TermT$ is entirely specified by $c$. Since
    an $\Signature$-term is entirely specified by its decoration word and its set of edges,
    this shows the unicity of $\TermT$ and entails the statement of the proposition.
\end{Proof}

\subsubsection{Easterly wind posets} \label{subsubsec:easterly_wind_posets}
Let $\EasterlyWindLeq$ be the binary relation on $\SetTerms \App \Signature$ such that, for
any $\Signature$-terms $\TermT_1$ and $\TermT_2$, we have $\TermT_1 \EasterlyWindLeq
\TermT_2$ if $\Decoration \App \TermT_1 = \Decoration \App \TermT_2$ and, for all $i \in
\InternalNodes \App \TermT_1$, $\Connection \App \TermT_1 \App i \leq \Connection \App
\TermT_2 \App i$. Immediately from its definition, $\EasterlyWindLeq$ is a partial order
relation on $\SetTerms \App \Signature$. Let us call $\Par{\SetTerms \App \Signature,
\EasterlyWindLeq}$ the \Def{$\Signature$-easterly wind poset}.

Let us now state a series of lemmas that will be used to establish that
$\EasterlyWindCovering$ is the covering relation of the $\Signature$-easterly wind poset.

\begin{Statement}{Lemma}{lem:order_domination}
    Let $\Signature$ be a signature and $\TermT_1$ and $\TermT_2$ be two $\Signature$-terms
    of the same degree $n$ and the same decoration word. We have $\TermT_1 \EasterlyWindLeq
    \TermT_2$ if and only if for any $i \in [n]$, the parent edge of $i$ in $\TermT_1$ is
    dominated by the parent edge of $i$ in $\TermT_2$.
\end{Statement}
\begin{Proof}
    Let $\Par{i_1, j_1, i}$ (resp.\ $\Par{i_2, j_2, i}$) be the parent edge of $i$ in
    $\TermT_1$ (resp.\ $\TermT_2$). By definition of~$\EasterlyWindLeq$, the property
    $\TermT_1 \EasterlyWindLeq \TermT_2$ is equivalent to the fact that for any $i \in [n]$,
    \begin{math}
        \Connection \App \TermT_1 \App i \leq \Connection \App \TermT_2 \App i.
    \end{math}
    By~\eqref{equ:connection} and Lemma~\ref{lem:connection_words_interval}, this is
    equivalent to the fact that $i_2 > i_1$ or both $i_1 = i_2$ and $j_2 \leq j_1$. This
    says exactly that $\Par{i_1, j_1, i}$ is dominated by $\Par{i_2, j_2, i}$. The statement
    of the lemma follows.
\end{Proof}

\begin{Statement}{Lemma}{lem:rewrite_rule_partial_order}
    For any signature $\Signature$ and any $\Signature$-terms $\TermT_1$ and $\TermT_2$,
    $\TermT_1 \EasterlyWindCoveringRT \TermT_2$ implies that $\TermT_1 \EasterlyWindLeq
    \TermT_2$.
\end{Statement}
\begin{Proof}
    Assume that $\TermT_1 \EasterlyWindCovering \TermT_2$. By
    Lemma~\ref{lem:covering_edge_properties}, $\Edges \App \TermT_2$ is obtained from
    $\Edges \App \TermT_1$ by replacing an edge $\Par{i_1, j_1, i}$ by a dominating edge
    $\Par{i_2, j_2, i}$. Therefore, by Lemma~\ref{lem:order_domination}, $\TermT_1
    \EasterlyWindLeq \TermT_2$. Finally, the statement of the lemma follows from the fact
    that $\EasterlyWindLeq$ is transitive.
\end{Proof}

\begin{Statement}{Lemma}{lem:partial_order_contraction}
    For any signature $\Signature$ and any $\Signature$-terms $\TermT_1$ and $\TermT_2$ of
    degree $1$ or more, if $\TermT_1 \EasterlyWindLeq \TermT_2$ then $\Contraction \App
    \TermT_1 \EasterlyWindLeq \Contraction \App \TermT_2$.
\end{Statement}
\begin{Proof}
    Assume that $\TermT_1$ and $\TermT_2$ are two $\Signature$-terms of the same degree $n
    \geq 1$ and that $\TermT_1 \EasterlyWindLeq \TermT_2$. By definition of
    $\EasterlyWindLeq$, for all $i \in [n]$, $\Connection \App \TermT_1 \App i \leq
    \Connection \App \TermT_2 \App i$. Observe that for any $\Signature$-term $\TermT$ of
    degree $n$,
    \begin{math}
        \Edges \App \TermT
        = \Edges \App \ULine{\Contraction \App \TermT}
        \cup \Bra{\Par{\Parent \App \TermT \App n, \LocalPosition \App \TermT \App n, n}}.
    \end{math}
    Therefore, this implies that for all $i \in [n - 1]$,
    \begin{math}
        \Connection \App \ULine{\Contraction \App \TermT_1} \App i =
        \Connection \App \TermT_1 \App i \leq
        \Connection \App \TermT_2 \App i
        = \Connection \App \ULine{\Contraction \App \TermT_2} \App i.
    \end{math}
    Hence, we have $\Contraction \App \TermT_1 \EasterlyWindLeq \Contraction \App \TermT_2$
    as expected.
\end{Proof}

\begin{Statement}{Lemma}{lem:covering_contraction}
    Let $\Signature$ be a signature, and $\TermT_1$ and $\TermT_2$ be two $\Signature$-terms
    of the same degree $n \geq 1$. If there exists $i \in [n - 1]$ such that $\Contraction
    \App \TermT_1 \EasterlyWindCovering_i \Contraction \App \TermT_2$, $\Decoration \App
    \TermT_1 \App n = \Decoration \App \TermT_2 \App n$, $\Parent \App \TermT_1 \App n =
    \Parent \App \TermT_2 \App n$, and $\LocalPosition \App \TermT_1 \App n = \LocalPosition
    \App \TermT_2 \App n$, then $\TermT_1 \EasterlyWindCovering_i \TermT_2$.
\end{Statement}
\begin{Proof}
    The $\Signature$-terms $\TermT_1$ and $\TermT_2$ are both obtained by adding
    respectively to $\Contraction \App \TermT_1$ and $\Contraction \App \TermT_2$ an
    internal node $n$ decorated by the same element of $\Signature$ and through the same
    edge. Since $\Contraction \App \TermT_2$ can be obtained from $\Contraction \App
    \TermT_1$ by changing a single edge involving internal nodes smaller than $n$, as
    prescribed by the definition of $\EasterlyWindCovering_i$, it is possible to obtain
    $\TermT_2$ from $\TermT_1$ by the same changing of edge. Therefore, $\TermT_1
    \EasterlyWindCovering_i \TermT_2$.
\end{Proof}

\begin{Statement}{Lemma}{lem:partial_order_sequence_rewrites}
    Let $\Signature$ be a signature, and $\TermT_1$ and $\TermT_2$ be two $\Signature$-terms
    such that $\TermT_1 \EasterlyWindLeq \TermT_2$. There exists a sequence
    $\Par{\TermT^{(0)}, \TermT^{(1)}\dots, \TermT^{(n)}}$ of $\Signature$-terms such that
    \begin{equation}
        \TermT_1
        = \TermT^{(0)} \EasterlyWindCoveringRT_1 \TermT^{(1)}
        \EasterlyWindCoveringRT_2 \cdots \EasterlyWindCoveringRT_n \TermT^{(n)}
        = \TermT_2.
    \end{equation}
\end{Statement}
\begin{Proof}
    Let us proceed by induction on the common degree $n$ of $\TermT_1$ and $\TermT_2$. If $n
    = 0$ or $n = 1$, then $\TermT_1 = \TermT_2$ and the property holds. Otherwise, we have
    $n \geq 2$ and, since $\TermT_1 \EasterlyWindLeq \TermT_2$, by
    Lemma~\ref{lem:partial_order_contraction}, we have $\Contraction \App \TermT_1
    \EasterlyWindLeq \Contraction \App \TermT_2$. By induction hypothesis, there exists a
    sequence $\Par{\TermS^{(0)}, \TermS^{(1)}\dots, \TermS^{(n - 1)}}$ of $\Signature$-terms
    such that
    \begin{equation}
        \Contraction \App \TermT_1
        = \TermS^{(0)} \EasterlyWindCoveringRT_1 \TermS^{(1)}
        \EasterlyWindCoveringRT_2 \cdots \EasterlyWindCoveringRT_n \TermS^{(n - 1)}
        = \Contraction \App \TermT_2.
    \end{equation}
    Let, for any $i \in \HanL{n - 1}$, $\TermR^{(k)}$ be the $\Signature$-term obtained by
    adding to $\TermS^{(k)}$ an internal node $n$ decorated by $\Decoration \App \TermT_1
    \App n$ through the edge $\Par{\Parent \App \TermT_1 \App n, \LocalPosition \App
    \TermT_1 \App n, n}$. By Lemma~\ref{lem:covering_contraction},
    \begin{equation}
        \TermT_1
        = \TermR^{(0)} \EasterlyWindCoveringRT_1 \TermR^{(1)}
        \EasterlyWindCoveringRT_2 \cdots \EasterlyWindCoveringRT_n \TermR^{(n - 1)}.
    \end{equation}
    Now, since
    \begin{math}
        \Connection \App \TermT_1 \App n
        = \Connection \App \TermR^{(n - 1)} \App n
        \leq \Connection \App \TermT_2 \App n,
    \end{math}
    by Lemma~\ref{lem:order_domination}, $\TermT_2$ is obtained from $\TermR^{(n - 1)}$ by
    replacing the parent edge of $j$ in $\TermR^{(n - 1)}$ by an edge dominating it. By
    definition of $\EasterlyWindCovering$ and Lemma~\ref{lem:covering_edge_properties}, the
    parent edge of $n$ in $\TermT_2$ can be formed from the parent edge of $n$ in
    $\TermR^{(n - 1)}$ by performing a sequence of applications of the $\Signature$-easterly
    wind rewrite rule $\EasterlyWindCovering$ from $\TermR^{(n - 1)}$. Indeed, this consists
    in iteratively moving the node $n$ of $\TermR^{(n - 1)}$ as specified by the binary
    relation $\EasterlyWindCovering_n$. Therefore, we have $\TermR^{(n - 1)}
    \EasterlyWindCoveringRT_n \TermT_2$, establishing the expected property.
\end{Proof}

\begin{Statement}{Theorem}{thm:easterly_wind_posets}
    For any signature $\Signature$, the binary relation $\EasterlyWindCovering$ is the
    covering relation of the $\Signature$-easterly wind poset.
\end{Statement}
\begin{Proof}
    By Lemmas~\ref{lem:rewrite_rule_partial_order}
    and~\ref{lem:partial_order_sequence_rewrites}, the binary relations $\EasterlyWindLeq$
    and $\EasterlyWindCoveringRT$ are the same. Besides, by
    Lemma~\ref{lem:covering_edge_properties}, if $\TermT_1$ and $\TermT_2$ are two
    $\Signature$-terms satisfying $\TermT_1 \EasterlyWindCovering \TermT_2$, then $\TermT_1$
    and $\TermT_2$ differ by the parent edge of a certain internal node $i$. Therefore, we
    have $\TermT_1 \EasterlyWindCovering_i \TermT_2$. Now, by contradiction, assume that
    there is an $\Signature$-term $\TermT_3$ such that $\TermT_3 \ne \TermT_2$ and $\TermT_1
    \EasterlyWindCovering \TermT_3 \EasterlyWindCoveringRT \TermT_2$. Recall that as noticed
    in Section~\ref{subsubsec:rewrite_rule}, $\TermT_2$ is the unique $\Signature$-term such
    that $\TermT_1 \EasterlyWindCovering_i \TermT_2$. Therefore, we have $\TermT_1
    \EasterlyWindCovering_{i'} \TermT_3$ with $i' \ne i$. The fact that $\TermT_1
    \EasterlyWindCovering_{i'} \TermT_3 \EasterlyWindCoveringRT \TermT_2$ implies that
    $\TermT_1$ and $\TermT_2$ differ by the parent edge of $i'$. This yields a contradiction
    with our hypotheses. This shows that $\EasterlyWindCovering$ is the covering relation of
    the poset $\Par{\SetTerms \App \Signature, \EasterlyWindLeq}$.
\end{Proof}

Theorem~\ref{thm:easterly_wind_posets} justifies the given name for $\Par{\SetTerms \App
\Signature, \EasterlyWindLeq}$: this name of $\Signature$-``\textit{easterly wind}'' poset
is derived from the observation that the covering relation of this poset involves detaching
a subterm from the east and then attaching it to the west, as if an easterly breeze is
blowing on the tree.

For any $\TermT \in \SetTerms \App \Signature$, let
\begin{math}
    \TerminalInterval \App \TermT
    :=
    \Bra{\TermT' \in \SetTerms \App \Signature : \TermT \EasterlyWindLeq \TermT'}.
\end{math}
We call $\Par{\TerminalInterval \App \TermT, \EasterlyWindLeq}$ the
\Def{$\Signature$-easterly wind poset of $\TermT$}.
Figure~\ref{fig:example_easterly_wind_poset} shows the Hasse diagram of the
$\SignatureExample$-easterly wind poset of an $\SignatureExample$-term.
\begin{figure}[ht]
    \centering
    \def\Zoom{0.75}
\def\Scale{4.00}
\scalebox{\Zoom}{
    \begin{tikzpicture}[Centering,scale=\Scale,x={(-1cm,-1cm)},y={(1cm,-1cm)}]
        \DrawSquareGrid{1.50}{2.75}{1.00}{3.00}{0.25}
        \node[GraphLabeledVertex](1)at(1.50,1.00){
            $\GenA_3 \ULine{\GenA_3 \LeafSimple \LeafSimple \LeafSimple} \ULine{\GenA_1
            \LeafSimple} \ULine{\GenA_2 \LeafSimple \LeafSimple}$
        };
        \node[GraphLabeledVertex](2)at(1.50,3.0){
            $\GenA_3 \ULine{\GenA_3 \LeafSimple \LeafSimple \LeafSimple} \ULine{\GenA_1
            \ULine{\GenA_2 \LeafSimple \LeafSimple}} \LeafSimple$
        };
        \node[GraphLabeledVertex](3)at(2.00,1.00){
            $\GenA_3 \ULine{\GenA_3 \LeafSimple \LeafSimple \ULine{\GenA_1 \LeafSimple}}
            \LeafSimple \ULine{\GenA_2 \LeafSimple \LeafSimple}$
        };
        \node[GraphLabeledVertex](4)at(2.00,1.50){
            $\GenA_3 \ULine{\GenA_3 \LeafSimple \LeafSimple \ULine{\GenA_1 \LeafSimple}}
            \ULine{\GenA_2 \LeafSimple \LeafSimple} \LeafSimple$
        };
        \node[GraphLabeledVertex](5)at(2.00,3.00){
            $\GenA_3 \ULine{\GenA_3 \LeafSimple \LeafSimple \ULine{\GenA_1 \ULine{\GenA_2
            \LeafSimple \LeafSimple}}} \LeafSimple \LeafSimple$
        };
        \node[GraphLabeledVertex](6)at(2.50,1.00){
            $\GenA_3 \ULine{\GenA_3 \LeafSimple \ULine{\GenA_1 \LeafSimple} \LeafSimple}
            \LeafSimple \ULine{\GenA_2 \LeafSimple \LeafSimple}$
        };
        \node[GraphLabeledVertex](7)at(2.50,1.50){
            $\GenA_3 \ULine{\GenA_3 \LeafSimple \ULine{\GenA_1 \LeafSimple} \LeafSimple}
            \ULine{\GenA_2 \LeafSimple \LeafSimple} \LeafSimple$
        };
        \node[GraphLabeledVertex](8)at(2.50,2.00){
            $\GenA_3 \ULine{\GenA_3 \LeafSimple \ULine{\GenA_1 \LeafSimple} \ULine{\GenA_2
            \LeafSimple \LeafSimple}} \LeafSimple \LeafSimple$
        };
        \node[GraphLabeledVertex](9)at(2.50,3.00){
            $\GenA_3 \ULine{\GenA_3 \LeafSimple \ULine{\GenA_1 \ULine{\GenA_2 \LeafSimple
            \LeafSimple}} \LeafSimple} \LeafSimple \LeafSimple$
        };
        \node[GraphLabeledVertex](10)at(2.75,1.00){
            $\GenA_3 \ULine{\GenA_3 \ULine{\GenA_1 \LeafSimple} \LeafSimple \LeafSimple}
            \LeafSimple \ULine{\GenA_2 \LeafSimple \LeafSimple}$
        };
        \node[GraphLabeledVertex](11)at(2.75,1.50){
            $\GenA_3 \ULine{\GenA_3 \ULine{\GenA_1 \LeafSimple} \LeafSimple \LeafSimple}
            \ULine{\GenA_2 \LeafSimple \LeafSimple} \LeafSimple$
        };
        \node[GraphLabeledVertex](12)at(2.75,2.00){
            $\GenA_3 \ULine{\GenA_3 \ULine{\GenA_1 \LeafSimple} \LeafSimple \ULine{\GenA_2
            \LeafSimple \LeafSimple}} \LeafSimple \LeafSimple$
        };
        \node[GraphLabeledVertex](13)at(2.75,2.50){
            $\GenA_3 \ULine{\GenA_3 \ULine{\GenA_1 \LeafSimple} \ULine{\GenA_2 \LeafSimple
            \LeafSimple} \LeafSimple} \LeafSimple \LeafSimple$
        };
        \node[GraphLabeledVertex](14)at(2.75,3.00){
            $\GenA_3 \ULine{\GenA_3 \ULine{\GenA_1 \ULine{\GenA_2 \LeafSimple \LeafSimple}}
            \LeafSimple \LeafSimple} \LeafSimple \LeafSimple$
        };
        \draw[GraphArc](1)--(2);
        \draw[GraphArc](1)--(3);
        \draw[GraphArc](2)--(5);
        \draw[GraphArc](3)--(4);
        \draw[GraphArc](3)--(6);
        \draw[GraphArc](4)--(5);
        \draw[GraphArc](4)--(7);
        \draw[GraphArc](5)--(9);
        \draw[GraphArc](6)--(7);
        \draw[GraphArc](6)--(10);
        \draw[GraphArc](7)--(8);
        \draw[GraphArc](7)--(11);
        \draw[GraphArc](8)--(9);
        \draw[GraphArc](8)--(12);
        \draw[GraphArc](9)--(14);
        \draw[GraphArc](10)--(11);
        \draw[GraphArc](11)--(12);
        \draw[GraphArc](12)--(13);
        \draw[GraphArc](13)--(14);
    \end{tikzpicture}
}
    \caption[]{
        The Hasse diagram of the $\SignatureExample$-easterly wind poset of
        $\GenA_3 \ULine{\GenA_3 \LeafSimple \LeafSimple \LeafSimple} \ULine{\GenA_1
        \LeafSimple} \ULine{\GenA_2 \LeafSimple \LeafSimple}$.
    }
    \label{fig:example_easterly_wind_poset}
\end{figure}

\subsection{Posets on tilted terms}
In this section, we consider an idempotent map $\Tilt \App \SpecialNodes$ on the set of
$\Signature$-terms. It turns out that this map is a closure operator on the
$\Signature$-easterly wind poset, so that the set of elements closed w.r.t.\ $\Tilt \App
\SpecialNodes$ forms a subposet of the $\Signature$-easterly wind poset. This construction
will be useful in the final section of this paper, as such posets are used to construct
bases of natural Hopf algebras of free operads.

\subsubsection{Tilting map} \label{subsubsec:tilting_map}
The \Def{$\SpecialNodes$-tilting map} is the map $\Tilt \App \SpecialNodes : \SetTerms \App
\Signature \to \SetTerms \App \Signature$ defined as follows. For any $\TermT \in \SetTerms
\App \Signature$, the $\Signature$-term $\Tilt \App \SpecialNodes \App \TermT$ is obtained
from $\TermT$ by rearranging the children of every internal node $i$ with $i \in
\SpecialNodes$ so that all non-leaf children preserve their original order and precede the
children that are leaves. Let us also define the \Def{$\SpecialNodes$-reversed tilting map}
as the map $\TiltReverse \App \SpecialNodes$ in the exact same manner as the
$\SpecialNodes$-tilting map but with the difference that the children different from the
leaf appear after the children which are leaves. For instance, we have
\begin{equation} \label{equ:example_tilt_map}
    \Tilt \App {\Bra{1, 2}} \App
    \GenA_3
    \ULine{
        \GenA_2
        \LeafSimple
        \ULine{
            \GenA_2
            \LeafSimple
            \ULine{\GenA_1 \LeafSimple}
        }
    }
    \LeafSimple
    \GenA_0
    \enspace = \enspace
    \GenA_3
    \ULine{
        \GenA_2
        \ULine{
            \GenA_2
            \LeafSimple
            \ULine{\GenA_1 \LeafSimple}
        }
        \LeafSimple
    }
    \GenA_0
    \LeafSimple
\end{equation}
and
\begin{equation}
    \TiltReverse \App {\Bra{1, 2}} \App
    \GenA_3
    \ULine{
        \GenA_2
        \LeafSimple
        \ULine{
            \GenA_2
            \LeafSimple
            \ULine{\GenA_1 \LeafSimple}
        }
    }
    \LeafSimple
    \GenA_0
    \enspace = \enspace
    \GenA_3
    \LeafSimple
    \ULine{
        \GenA_2
        \LeafSimple
        \ULine{
            \GenA_2
            \LeafSimple
            \ULine{\GenA_1 \LeafSimple}
        }
    }
    \GenA_0.
\end{equation}

Given an $\Signature$-term $\TermT$ and an internal node $i$ of $\TermT$, let $\LeftBrothers
\App \TermT \App i$ be the number of internal nodes of $\TermT$ among the siblings of $i$
that lie to its left, including $i$ itself. For instance, by setting $\TermT$ as the
$\SignatureExample$-term appearing in the left-hand side of~\eqref{equ:example_tilt_map}, we
have $\LeftBrothers \App \TermT \App 1 = 0$, $\LeftBrothers \App \TermT \App 2 = 1$, and
$\LeftBrothers \App \TermT \App 5 = 2$. Observe, of course, that if $\Par{i_1, j, i}$ is an
edge of $\TermT$, then $\LeftBrothers \App \TermT \App i \leq j$.

The following lemma provides a formalization of the effect of the $\SpecialNodes$-tilting
map on an $\Signature$-term.

\begin{Statement}{Lemma}{lem:tilt_map}
    Let $\Signature$ be a signature, $\SpecialNodes$ be a set of positive integers, $\TermT$
    be an $\Signature$-term and $\Par{i', j, i}$ be an edge of $\TermT$. The following
    properties hold:
    \begin{enumerate}[label=(\roman*)]
        \item \label{item:tilt_map_1}
        if $i' \notin \SpecialNodes$, then $\Par{i', j, i}$ is an edge of $\Tilt \App
            \SpecialNodes \App \TermT$;
        \item \label{item:tilt_map_2}
        if $i' \in \SpecialNodes$, then $\Par{i', \LeftBrothers \App \TermT \App i, i}$ is
            an edge of $\Tilt \App \SpecialNodes \App \TermT$.
    \end{enumerate}
\end{Statement}
\begin{Proof}
    By definition of the map $\Tilt \App \SpecialNodes$, if $i' \notin \SpecialNodes$, then
    the children of $i'$ in $\TermT$ and in $\Tilt \App \SpecialNodes \App \TermT$ are
    arranged in the same way. This implies~\ref{item:tilt_map_1}. Besides, when $i' \in
    \SpecialNodes$, in order to obtain $\Tilt \App \SpecialNodes \App \TermT$, the children
    of $i'$ in $\TermT$ are pushed to the left taking the place of potential leaves. Since
    the children of $i'$ which are not leaves remain in the same relative order, each
    internal node $i$ of $\Tilt \App \SpecialNodes \App \TermT$ which is a children of $i'$
    appears at the $\LeftBrothers \App \TermT \App i$-th position. This
    implies~\ref{item:tilt_map_2}.
\end{Proof}

\subsubsection{A closure operator} \label{subsubsec:closure_operator}
A map $\phi : \Poset \to \Poset$ is a \Def{closure operator} (see~\cite{DP02}) of a poset
$\Par{\Poset, \Leq_\Poset}$ if the following three properties hold:
\begin{enumerate}[label=(C\arabic*)]
    \item \label{item:closure_operator_1}
    for any $x \in \Poset$, $x \Leq_\Poset \phi \App x$;
    \item \label{item:closure_operator_2}
    for any $x, x' \in \Poset$, $x \Leq_\Poset x'$ implies $\phi \App x \Leq_\Poset \phi
    \App x'$;
    \item \label{item:closure_operator_3}
    for any $x \in \Poset$, $\phi \App \ULine{\phi \App x} = \phi \App x$.
\end{enumerate}
Condition~\ref{item:closure_operator_1} says that $\phi$ is extensive,
Condition~\ref{item:closure_operator_2} says that $\phi$ is order-preserving, and
Condition~\ref{item:closure_operator_3} says that $\phi$ is idempotent.

\begin{Statement}{Lemma}{lem:left_brothers}
    Let $\Signature$ be a signature, and $\TermT_1$ and $\TermT_2$ be two $\Signature$-terms
    such that $\TermT_1 \EasterlyWindLeq \TermT_2$. If $i$ is an internal node of both
    $\TermT_1$ and $\TermT_2$ admitting the same parent in $\TermT_1$ and $\TermT_2$, then
    $\LeftBrothers \App \TermT_2 \App i \leq \LeftBrothers \App \TermT_1 \App i$.
\end{Statement}
\begin{Proof}
    Assume that $\TermT_1 \EasterlyWindCovering_{i'} \TermT_2$ where $i' \in \InternalNodes
    \App \TermT_1$. It is immediate, by definition of $\EasterlyWindCovering_{i'}$ that if
    $i'$ is not a left brother of $i$ in $\TermT_1$, then $\LeftBrothers \App \TermT_1 \App
    i = \LeftBrothers \App \TermT_2 \App i$. Assume now that $i'$ is a left brother of $i$
    in $\TermT_1$ (including the case $i' = i$). Again by definition of
    $\EasterlyWindCovering_{i'}$, if the leaf which is visited immediately before $i'$ in
    the preorder traversal of $\TermT_1$ is a child of $\Parent \App \TermT_1 \App i$, then
    $\LeftBrothers \App \TermT_2 \App i = \LeftBrothers \App \TermT_1 \App i$. Otherwise,
    $\LeftBrothers \App \TermT_2 \App i = \LeftBrothers \App \TermT_1 \App i - 1$. The facts
    that, by Theorem~\ref{thm:easterly_wind_posets}, $\EasterlyWindCovering$ is the covering
    relation of the partial order relation $\EasterlyWindLeq$, and that
    $\EasterlyWindCovering$ is the union of all $\EasterlyWindCovering_{i'}$ with $i' \geq
    1$, entail the statement of the lemma.
\end{Proof}

\begin{Statement}{Theorem}{thm:tilt_closure}
    For any signature $\Signature$ and any set $\SpecialNodes$ of positive integers, the map
    $\Tilt \App \SpecialNodes$ is a closure operator of the $\Signature$-easterly wind
    poset.
\end{Statement}
\begin{Proof}
    Let $\TermT \in \SetTerms \App \Signature$. By Lemma~\ref{lem:tilt_map}, for any edge
    $\Par{i', j_1, i}$ of $\TermT$, there is an edge $\Par{i', j_2, i}$ of $\Tilt \App
    \SpecialNodes \App \TermT$ such that $j_2 \leq j_1$. Hence, the former edge is dominated
    by the latter. Therefore, by Lemma~\ref{lem:order_domination}, $\TermT \EasterlyWindLeq
    \Tilt \App \SpecialNodes \App \TermT$, showing that $\Tilt \App \SpecialNodes$
    satisfies~\ref{item:closure_operator_1}.

    Let $\TermT_1, \TermT_2 \in \SetTerms \App \Signature$ such that $\TermT_1
    \EasterlyWindLeq \TermT_2$. Assume that $\Par{i_1, j_1, i}$ is an edge of $\TermT_1$. By
    Lemma~\ref{lem:order_domination}, there is an edge $\Par{i_2, j_2, i}$ of $\TermT_2$
    such that the former edge is dominated by the latter. Moreover, by
    Lemma~\ref{lem:tilt_map}, $\Tilt \App \SpecialNodes \App \TermT_1$ has an edge
    $\Par{i_1, j_1', i}$ with $j_1' \in \Bra{j_1, \LeftBrothers \App \TermT \App i}$. Again
    by Lemma~\ref{lem:tilt_map}, we have also that $\Par{i_2, j_2', i}$ is an edge of $\Tilt
    \App \SpecialNodes \App \TermT_2$ with $j_2' \in \Bra{j_2, \LeftBrothers \App \TermT_2
    \App i}$. Now, we have two cases depending on how $\Par{i_1, j_1, i}$ is dominated by
    $\Par{i_2, j_2, i}$.
    \begin{enumerate}
        \item If $i_1 < i_2$, then $\Par{i_1, j_1', i}$ is dominated by $\Par{i_2, j_2',
        i}$.
        \item Otherwise, we have necessarily that $i_1 = i_2$ and $j_2' \leq j_1'$. Now, by
        Lemma~\ref{lem:tilt_map}, if $i_1 \notin \SpecialNodes$, then both $j_1' = j_1$ and
        $j_2' = j_2$ hold. Otherwise, when $i_1 \in \SpecialNodes$, we have $j_1' =
        \LeftBrothers \App \TermT_1 \App i$ and $j_2' = \LeftBrothers \App \TermT_2 \App i$.
        By Lemma~\ref{lem:left_brothers}, we have in particular that $\LeftBrothers \App
        \TermT_2 \App i \leq \LeftBrothers \App \TermT_1 \App i$. It follows that in both
        sub-cases, $\Par{i_1, j_1', i}$ is dominated by $\Par{i_1, j_2', i}$.
    \end{enumerate}
    From all this, it follows that $\Par{i_1, j_1', i}$ is dominated by $\Par{i_2, j_2',
    i}$. Therefore, by Lemma~\ref{lem:order_domination}, this implies that
    \begin{math}
        \Tilt \App \SpecialNodes \App \TermT_1
        \EasterlyWindLeq
        \Tilt \App \SpecialNodes \App \TermT_2
    \end{math}
    and shows that $\Tilt \App \SpecialNodes$ satisfies~\ref{item:closure_operator_2}.

    Finally, the map $\Tilt \App \SpecialNodes$ is, immediately from its definition,
    idempotent. Therefore, \ref{item:closure_operator_3} holds.
\end{Proof}

Observe, contrary to the property highlighted by Theorem~\ref{thm:tilt_closure} for the map
$\Tilt \App \SpecialNodes$, the map $\TiltReverse \App \SpecialNodes$ is not a closure
operator of the $\Signature$-easterly wind poset. Indeed, this map is not extensive. A
counterexample involves the $\SignatureExample$-easterly wind poset, the set $\SpecialNodes
:= \Bra{1}$, and the $\SignatureExample$-term
\begin{math}
    \TermT := \GenA_2 \ULine{\GenA_1 \LeafSimple} \LeafSimple
\end{math}
because we have
\begin{math}
    \TermT
    \; \cancel{\EasterlyWindLeq} \;
    \GenA_2 \LeafSimple \ULine{\GenA_1 \LeafSimple}
    =
    \TiltReverse \App \SpecialNodes \App \TermT.
\end{math}
Moreover, the $\SpecialNodes$-reversed tilting map is not either order-preserving. A
counterexample involves the $\SignatureExample$-easterly wind poset, the set $\SpecialNodes
:= [4]$, and the $\SignatureExample$-terms
\begin{math}
    \TermT_1 :=
    \GenA_3
    \ULine{\GenA_2 \ULine{\GenA_1 \LeafSimple} \GenA_0}
    \LeafSimple
    \LeafSimple
\end{math}
and
\begin{math}
    \TermT_2 :=
    \GenA_3
    \ULine{\GenA_2 \ULine{\GenA_1 \GenA_0} \LeafSimple}
    \LeafSimple
    \LeafSimple
\end{math}
because we have $\TermT_1 \EasterlyWindLeq \TermT_2$ but
\begin{math}
    \TiltReverse \App \SpecialNodes \App \TermT_1
    =
    \GenA_3
    \LeafSimple
    \LeafSimple
    \ULine{\GenA_2 \ULine{\GenA_1 \LeafSimple} \GenA_0}
    \; \cancel{\EasterlyWindLeq} \;
    \GenA_3
    \LeafSimple
    \LeafSimple
    \ULine{\GenA_2 \LeafSimple \ULine{\GenA_1 \GenA_0}}
    =
    \TiltReverse \App \SpecialNodes \App \TermT_2.
\end{math}

\subsubsection{Kernel of the tilting map and intervals} \label{subsubsec:kernel_tilt}
Let us denote by $\Equiv_{\Tilt \App \SpecialNodes}$ the kernel of $\Tilt \App
\SpecialNodes$, that is, the equivalence relation on $\SetTerms \App \Signature$ such that
for any $\Signature$-terms $\TermT$ and $\TermT'$, $\TermT \Equiv_{\Tilt \App \SpecialNodes}
\TermT'$ holds whenever $\Tilt \App \SpecialNodes \App \TermT = \Tilt \App \SpecialNodes
\App \TermT'$. For instance,
\begin{equation}
    \GenA_4
    \ULine{\GenA_2 \LeafSimple \ULine{\GenA_1 \LeafSimple}}
    \LeafSimple
    \LeafSimple
    \ULine{\GenA_2 \LeafSimple \ULine{\GenA_2 \LeafSimple \LeafSimple}}
    \Equiv_{\Tilt \App \Bra{1, 4}}
    \GenA_4
    \LeafSimple
    \ULine{\GenA_2 \LeafSimple \ULine{\GenA_1 \LeafSimple}}
    \LeafSimple
    \ULine{\GenA_2 \ULine{\GenA_2 \LeafSimple \LeafSimple} \LeafSimple}.
\end{equation}
Observe that $\TermT \Equiv_{\Tilt \App \SpecialNodes} \TermT'$ if and only if $\TiltReverse
\App \SpecialNodes \App \TermT = \TiltReverse \App \SpecialNodes \App \TermT'$.

\begin{Statement}{Proposition}{prop:kernel_tilt_intervals}
    Let $\Signature$ be a signature and $\SpecialNodes$ be a set of positive integers. The
    $\Equiv_{\Tilt \App \SpecialNodes}$-equivalence class of an $\Signature$-term $\TermT$
    is an interval of the $\Signature$-easterly wind poset. More specifically,
    \begin{equation}
        \Han{\TermT}_{\Equiv_{\Tilt \App \SpecialNodes}}
        =
        \Han{
            \TiltReverse \App \SpecialNodes \App \TermT,
            \Tilt \App \SpecialNodes \App \TermT
        }.
    \end{equation}
\end{Statement}
\begin{Proof}
    By Theorem~\ref{thm:tilt_closure}, $\TermT \EasterlyWindLeq \Tilt \App \SpecialNodes
    \App \TermT$. Similarly, by very analogous arguments as the one used in the proof of
    this property, we have $\TiltReverse \App \SpecialNodes \App \TermT \EasterlyWindLeq
    \TermT$. Now, let $\TermT' \in \Han{\TermT}_{\Equiv_{\Tilt \App \SpecialNodes}}$. Since
    $\TermT' \Equiv_{\Tilt \App \SpecialNodes} \TermT$, we have $\Tilt \App \SpecialNodes
    \App \TermT' = \Tilt \App \SpecialNodes \App \TermT$ and $\TiltReverse \App
    \SpecialNodes \App \TermT' = \TiltReverse \App \SpecialNodes \App \TermT$. Therefore,
    from the above property, we have $\TiltReverse \App \SpecialNodes \App \TermT
    \EasterlyWindLeq \TermT' \EasterlyWindLeq \Tilt \App \SpecialNodes \App \TermT$.

    Assume now that
    \begin{math}
        \TermT' \in \Han{
            \TiltReverse \App \SpecialNodes \App \TermT,
            \Tilt \App \SpecialNodes \App \TermT
        }
    \end{math}
    and let $\Par{i, j, i'}$ be an edge of $\TermT'$. If $i \in \SpecialNodes$, then by
    Lemma~\ref{lem:order_domination}, $\Par{i, j, i'}$ is an edge of both $\TiltReverse \App
    \SpecialNodes \App \TermT$ and $\Tilt \App \SpecialNodes \App \TermT$. Otherwise, when
    $i \notin \SpecialNodes$, again by Lemma~\ref{lem:order_domination}, $\Par{i, j', i'}$
    is an edge of $\TiltReverse \App \SpecialNodes \App \TermT$ and $\Par{i, j'', i'}$ is an
    edge of $\Tilt \App \SpecialNodes \App \TermT$ with $j'' \leq j \leq j'$. Therefore, by
    definition of $\Equiv_{\Tilt \App \SpecialNodes}$, we have $\TermT' \Equiv_{\Tilt \App
    \SpecialNodes} \TermT$ so that $\TermT' \in \Han{\TermT}_{\Equiv_{\Tilt \App
    \SpecialNodes}}$.
\end{Proof}

\subsubsection{Closed elements}
An element $x$ of a poset $\Poset$ is \Def{closed} w.r.t.\ a closure operator $\phi$ of
$\Poset$ if $x$ is a fixed point of $\phi$. In this way, since by
Theorem~\ref{thm:tilt_closure}, for any set $\SpecialNodes$ of positive integers, $\Tilt
\App \SpecialNodes$ is a closure operator of $\Par{\SetTerms \App \Signature,
\EasterlyWindLeq}$, closed elements w.r.t.\ $\Tilt \App \SpecialNodes$ are well-defined and
are called \Def{$\SpecialNodes$-tilted}. For instance, the $\SignatureExample$-term
\begin{math}
    \GenA_3
    \ULine{
        \GenA_2
        \LeafSimple
        \ULine{
            \GenA_2
            \ULine{
                \GenA_1
                \LeafSimple
            }
            \LeafSimple
        }
    }
    \GenA_0
    \ULine{
        \GenA_2
        \LeafSimple
        \LeafSimple
    }
\end{math}
is $\Bra{1, 3, 6}$-tilted but is not $\Bra{2}$-tilted.

For any set $\SpecialNodes$ of positive integers, $\Leq$ is a partial order on $\Tilt \App
\SpecialNodes \App \ULine{\SetTerms \App \Signature}$. Let us call $\Par{\Tilt \App
\SpecialNodes \App \ULine{\SetTerms \App \Signature}, \EasterlyWindLeq}$ the
\Def{$\SpecialNodes$-tilted $\Signature$-easterly wind poset}. Observe that for any sets
$\SpecialNodes_1$ and $\SpecialNodes_2$ of positive integers, if $\SpecialNodes_1 \subseteq
\SpecialNodes_2$, then $\Par{\Tilt \App \SpecialNodes_2 \App \ULine{\SetTerms \App
\Signature} , \EasterlyWindLeq}$ is a subposet of $\Par{\Tilt \App \SpecialNodes_1 \App
\ULine{\SetTerms \App \Signature}, \EasterlyWindLeq}$. Of course, $\Par{\Tilt \App \emptyset
\App \ULine{\SetTerms \App \Signature}, \EasterlyWindLeq}$ is the $\Signature$-easterly wind
poset introduced in Section~\ref{subsubsec:easterly_wind_posets}. Moreover, for any $\TermT
\in \SetTerms \App \Signature$, let
\begin{math}
    \TerminalInterval \App \SpecialNodes \App \TermT
    := \ULine{\TerminalInterval \App \TermT}
    \cap
    \ULine{\Tilt \App \SpecialNodes \App \TermT}.
\end{math}
We call $\Par{\TerminalInterval \App \SpecialNodes \App \TermT, \EasterlyWindLeq}$ the
\Def{$\SpecialNodes$-tilted $\Signature$-easterly wind poset of $\TermT$}.
Figure~\ref{fig:example_tilted_easterly_wind_poset} shows the Hasse diagram of the
$\SpecialNodes$-tilted $\SignatureExample$-easterly wind poset of a term.
\begin{figure}[ht]
    \centering
    \def\Zoom{0.85}
\def\Scale{2.20}
\scalebox{\Zoom}{
    \begin{tikzpicture}[Centering,scale=\Scale,x={(-1cm,-1cm)},y={(1cm,-1cm)}]
        \DrawSquareGrid{1.50}{2.75}{1.00}{3.00}{0.25}
        \node[GraphLabeledVertex](1)at(1.50,1.00){
            $\GenA_3 \ULine{\GenA_3 \LeafSimple \LeafSimple \LeafSimple} \ULine{\GenA_1
            \LeafSimple} \ULine{\GenA_2 \LeafSimple \LeafSimple}$ };
        \node[GraphLabeledVertex](2)at(1.50,3.00){
            $\GenA_3 \ULine{\GenA_3 \LeafSimple \LeafSimple \LeafSimple} \ULine{\GenA_1
            \ULine{\GenA_2 \LeafSimple \LeafSimple}} \LeafSimple$ };
        \node[GraphLabeledVertex](3)at(2.75,1.50){
            $\GenA_3 \ULine{\GenA_3 \ULine{\GenA_1 \LeafSimple} \LeafSimple \LeafSimple}
            \ULine{\GenA_2 \LeafSimple \LeafSimple} \LeafSimple$ };
        \node[GraphLabeledVertex](4)at(2.75,2.50){
            $\GenA_3 \ULine{\GenA_3 \ULine{\GenA_1 \LeafSimple} \ULine{\GenA_2 \LeafSimple
            \LeafSimple} \LeafSimple} \LeafSimple \LeafSimple$ };
        \node[GraphLabeledVertex](5)at(2.75,3.00){
            $\GenA_3 \ULine{\GenA_3 \ULine{\GenA_1 \ULine{\GenA_2 \LeafSimple \LeafSimple}}
            \LeafSimple \LeafSimple} \LeafSimple \LeafSimple$ };
        \draw[GraphArc](1)--(2);
        \draw[GraphArc](1)--(3);
        \draw[GraphArc](2)--(5);
        \draw[GraphArc](3)--(4);
        \draw[GraphArc](4)--(5);
    \end{tikzpicture}
}
    \caption[]{
        The Hasse diagram of the $\Bra{1, 2}$-tilted $\SignatureExample$-easterly wind poset
        of $\GenA_3 \ULine{\GenA_3 \LeafSimple \LeafSimple \LeafSimple} \ULine{\GenA_1
        \LeafSimple} \ULine{\GenA_2 \LeafSimple \LeafSimple}$.
    }
    \label{fig:example_tilted_easterly_wind_poset}
 \end{figure}

\subsubsection{Scope sequences and fully tilted terms} \label{subsubsec:fully_tilted_terms}
The \Def{scope sequence} of an $\Signature$-term $\TermT$ is the word $\Scope \App \TermT$
on $\N$ of length $\Deg \App \TermT$ such that for any $i \in \InternalNodes \App \TermT$,
$\Scope \App \TermT \App i$ is the number descendants of $i$ in $\TermT$. For instance,
\begin{equation}
    \GenA_4
    \GenA_0
    \ULine{
        \GenA_2 \ULine{\GenA_1 \LeafSimple} \ULine{\GenA_2 \LeafSimple \LeafSimple}
    }
    \ULine{
        \GenA_1 \ULine{\GenA_2 \LeafSimple \LeafSimple}
    }
    = 6 0 2 0 0 1 0.
\end{equation}

\begin{Statement}{Lemma}{lem:scope_sequence_parent}
    Let $\Signature$ be a signature and $\TermT$ be an $\Signature$-term. For any internal
    node $i$ of $\TermT$ different from the root, the parent of $i$ in $\TermT$ is the
    greatest internal node $i'$ of $\TermT$ such that $i' < i$ and $i' + \Scope \App \TermT
    \App i' \geq i$.
\end{Statement}
\begin{Proof}
    By definition of the word $\Scope \App \TermT$, in $\TermT$, an internal node $i$ is a
    descendant of an internal node $i'$ if and only if $i' + 1 \leq i \leq i' + \Scope \App
    \TermT \App i'$. Therefore, each ancestor $i'$ of $i$ in $\TermT$ satisfies $i' < i$ and
    $i' + \Scope \App \TermT \App i' \geq i$. The parent of $i$ in $\TermT$ is the greatest
    internal node among the ancestors of $i$. The statement of the lemma follows.
\end{Proof}

An $\Signature$-term $\TermT$ is \Def{fully tilted} if $\TermT$ is $\N \setminus
\Bra{0}$-tilted.

\begin{Statement}{Proposition}{prop:scope_sequence_fully_tilted}
    Let $\Signature$ be a signature, and $\TermT_1$ and $\TermT_2$ be two fully tilted
    $\Signature$-terms of the same degree $n$. We have $\TermT_1 \EasterlyWindLeq \TermT_2$
    if and only if, for any $i \in [n]$, $\Scope \App \TermT_1 \App i \leq \Scope \App
    \TermT_2 \App i$.
\end{Statement}
\begin{Proof}
    Assume that $\TermT_1$ and $\TermT_2$ are two fully tilted $\Signature$-terms such that
    $\TermT_1 \EasterlyWindLeq \TermT_2$. By Lemmas~\ref{lem:covering_edge_properties}
    and~\ref{lem:partial_order_sequence_rewrites}, for any internal node $i$ of both
    $\TermT_1$ and $\TermT_2$, the number of descendants of $i$ in $\TermT_1$ is smaller
    than or equal to the number of descendants of $i$ in $\TermT_2$. This implies that
    $\Scope \App \TermT_1 \App i \leq \Scope \App \TermT_2 \App i$.

    Conversely, assume that for any $i \in [n]$, $\Scope \App \TermT_1 \App i \leq \Scope
    \App \TermT_2 \App i$. Let $i \in [n]$ and $i_1$ (resp.\ $i_2$) be the parent of $i$ in
    $\TermT_1$ (resp.\ $\TermT_2$). Lemma~\ref{lem:scope_sequence_parent} together with the
    fact that each letter of $\Scope \App \TermT_2$ is greater than or equal to the letter
    at the same position of $\Scope \App \TermT_1$ imply that $i_2 \geq i_1$. Hence, the
    parent edge of $i$ in $\TermT_1$ (resp.\ $\TermT_2$) is $\Par{i_1, j_1, i}$ (resp.\
    $\Par{i_2, j_2, i}$) for some integer $j_1$ (resp.\ $j_2$). Since $\TermT_1$ (resp.\
    $\TermT_2$) is fully tilted, $j_1 = \LeftBrothers \App \TermT_1 \App i$ (resp.\ $j_2 =
    \LeftBrothers \App \TermT_2 \App i$). This implies that $i_1 < i_2$, or both $i_1 = i_2$
    and $j_1 \geq j_2$, so that the edge $\Par{i_1, j_1, i}$ is dominated by the edge
    $\Par{i_2, j_2, i}$. By Lemma~\ref{lem:order_domination}, we have $\TermT_1
    \EasterlyWindLeq \TermT_2$, as expected.
\end{Proof}

For any fully tilted $\Signature$-term $\TermT$, the \Def{fully tilted $\Signature$-easterly
wind poset of $\TermT$} is the $\N \setminus \Bra{0}$-tilted $\Signature$-easterly wind
poset of~$\TermT$.

\section{Geometric and lattice properties} \label{sec:geometric_and_lattice_properties}
In this section, we continue to establish properties of the easterly wind posets by focusing
in particular on geometric properties and on showing that terminal intervals of these posets
are lattices.

\subsection{Geometric properties}
We begin this section by showing that the $\Signature$-easterly wind posets are
EL-shellable. To this end, we introduce an encoding of the saturated chains of these posets.
Next, we propose a realization of the $\SpecialNodes$-tilted $\Signature$-easterly wind
posets as geometric objects using connection words introduced in
Section~\ref{subsubsec:connection_words}.

\subsubsection{EL-shellability} \label{subsubsec:el_shellability}
Given two $\Signature$-terms $\TermT_1$ and $\TermT_2$ of common degree $n \geq 0$ such that
$\TermT_1 \EasterlyWindLeq \TermT_2$, a \Def{$\Par{\TermT_1, \TermT_2}$-sequence} is a word
$u$ on $[n]$ such that
\begin{equation} \label{equ:sequence_saturated_chain}
    \TermT_1
    = \TermT^{(0)} \EasterlyWindCovering_{u \App 1}
    \TermT^{(1)} \EasterlyWindCovering_{u \App 2}
    \dots
    \EasterlyWindCovering_{u \App \ULine{\Length \App u}}
    \TermT^{(\Length \App u)}
    = \TermT_2
\end{equation}
for some $\Signature$-terms $\TermT^{(0)}$, $\TermT^{(1)}$, \dots, $\TermT^{(\Length \App
u)}$. Note that, by Theorem~\ref{thm:easterly_wind_posets}, the sequence of these terms
forms a saturated chain from $\TermT_1$ to $\TermT_2$ in the $\Signature$-easterly wind
poset. We say that this saturated chain is \Def{specified} by~$u$.

\begin{Statement}{Lemma}{lem:sequence_specification}
    Let $\Signature$ be a signature, and $\TermT_1$ and $\TermT_2$ be two $\Signature$-terms
    such that $\TermT_1 \EasterlyWindLeq \TermT_2$. The set of saturated chains between
    $\TermT_1$ and $\TermT_2$ is in one-to-one correspondence with the set of
    $\Par{\TermT_1, \TermT_2}$-sequences.
\end{Statement}
\begin{Proof}
    Let $\phi$ be the map having the set of $\Par{\TermT_1, \TermT_2}$-sequences as domain
    and the set of saturated chains between $\TermT_1$ and $\TermT_2$ as codomain, sending
    any $\Par{\TermT_1, \TermT_2}$-sequence $u$ to the saturated chain $\Par{\TermT^{(0)},
    \TermT^{(1)}, \dots, \TermT^{(\ell)}}$ defined accordingly
    with~\eqref{equ:sequence_saturated_chain}. First, recall that as noticed in
    Section~\ref{subsubsec:rewrite_rule}, given an $\Signature$-term $\TermS_1$ and an
    internal node $i$ of $\TermS_1$, there is at most one $\Signature$-term $\TermS_2$ such
    that $\TermS_1 \EasterlyWindCovering_i \TermS_2$. This implies that $\phi$ is a
    well-defined map. Moreover, as a consequence of
    Lemma~\ref{lem:covering_edge_properties}, for any $\Signature$-terms $\TermS_1$ and
    $\TermS_2$, if $\TermS_1 \EasterlyWindCovering_i \TermS_2$ and $\TermS_1
    \EasterlyWindCovering_{i'} \TermS_2$ for two internal nodes $i$ and $i'$ of $\TermS_1$,
    then $i = i'$. This implies that any saturated chain $\Par{\TermT^{(0)}, \TermT^{(1)},
    \dots, \TermT^{(\ell)}}$ between $\TermT_1$ and $\TermT_2$ admits exactly one
    $\Par{\TermT_1, \TermT_2}$-sequence which is an antecedent by $\phi$. Therefore $\phi$
    is a bijection and the statement of the lemma follows.
\end{Proof}

Lemma~\ref{lem:sequence_specification} entails in particular that the notion of saturated
chain specified by a $\Par{\TermT_1, \TermT_2}$-sequence is well-defined.

\begin{Statement}{Lemma}{lem:unique_sequences_specification}
    Let $\Signature$ be a signature, and $\TermT_1$ and $\TermT_2$ be two $\Signature$-terms
    such that $\TermT_1 \EasterlyWindLeq \TermT_2$. Among all $\Par{\TermT_1,
    \TermT_2}$-sequences,
    \begin{enumerate}[label=(\roman*)]
        \item \label{item:unique_sequences_specification_1}
        there is at most one which is a weakly increasing word;
        \item \label{item:unique_sequences_specification_2}
        there is at most one which is a weakly decreasing word;
        \item \label{item:unique_sequences_specification_3}
        the one specifying the saturated chain from $\TermT_1$ to $\TermT_2$ which is
        induced by Lemma~\ref{lem:partial_order_sequence_rewrites} is a weakly increasing
        word;
        \item \label{item:unique_sequences_specification_4}
        this weakly increasing $\Par{\TermT_1, \TermT_2}$-sequence is lexicographically
        smaller than any other $\Par{\TermT_1, \TermT_2}$-sequence.
    \end{enumerate}
\end{Statement}
\begin{Proof}
    To prove the uniqueness of a weakly increasing $\Par{\TermT_1, \TermT_2}$-sequence,
    assume that $u$ and $u'$ are two such weakly increasing $\Par{\TermT_1,
    \TermT_2}$-sequences. Let $i$ be a positive integer such that $u$ and $u'$ have both the
    same number of occurrences of each letter $i'$ for any $i' < i$. Since $i = 1$ always
    satisfies this condition, such an $i$ exists. Let $v$ be the prefix of $u$ made of the
    letters smaller than $i$. By the previously stated property on $i$ and the fact that $u$
    and $u'$ are weakly increasing, $v$ is equivalently the prefix of $u'$ made of the
    letters smaller than $i$. By Lemma~\ref{lem:sequence_specification}, there is a unique
    $\Signature$-term $\TermS_1$ such that $v$ is a $\Par{\TermT_1, \TermS_1}$-sequence.
    Now, let $\TermS_2$ be the $\Signature$-term obtained from $\TermS_1$ by applying $k$
    times the rewrite rule $\EasterlyWindCovering_i$ where $k$ is the number of occurrences
    of $i$ in $u$. By Lemma~\ref{lem:covering_edge_properties}, $\TermS_2$ is obtained from
    $\TermS_1$ by moving iteratively the parent edge of $i$. Moreover, again by
    Lemma~\ref{lem:covering_edge_properties}, for all $i' \ne i$, each application of
    $\EasterlyWindCovering_{i'}$ on $\TermS_2$ does not change the edge connecting $i$ to
    its parent. Therefore, in order to obtaining $\TermT_2$ from $\TermS_2$ through both
    $\Par{\TermT_1, \TermS_1}$-sequences $u$ and $u'$, all these properties imply that the
    number of occurrences of $i$ is the same in $u$ and $u'$. This entails $u = u'$
    and~\ref{item:unique_sequences_specification_1}. In a completely similar way, this also
    shows that there is at most one $\Par{\TermT_1, \TermT_2}$-sequence which is a weakly
    decreasing word. Thus, \ref{item:unique_sequences_specification_2} holds.

    Point~\ref{item:unique_sequences_specification_3} is immediate by
    Lemma~\ref{lem:partial_order_sequence_rewrites}. Indeed, the saturated chain induced by
    the chain of the statement of the lemma which has just been cited is specified by the
    $\Par{\TermT_1, \TermT_2}$-sequence $u$ consisting in a possibly empty block of
    letter $1$, then a possibly empty block of the letter $2$, and so on. Therefore,
    \ref{item:unique_sequences_specification_3} checks out. The fact that $u$ is
    lexicographically smaller than any other $\Par{\TermT_1, \TermT_2}$-sequence holds by
    construction. Indeed, each letter of $u$ is the smallest possible in order to specify
    the right saturated chain from $\TermT_1$ to $\TermT_2$. Hence,
    \ref{item:unique_sequences_specification_4} holds.
\end{Proof}

We use the standard definitions about labelings of Hasse diagrams of posets and EL-labelings
as given in~\cite{Bjo80,BW96}, which we recall here. A \Def{labeling} of a poset
$\Par{\Poset, \Leq_\Poset}$ is a map $\lambda : \Covered_{\Poset} \to \Lambda$ where
$\Covered_{\Poset}$ is the covering relation of $\Poset$ and $\Par{\Lambda, \Leq_\Lambda}$
is a poset. Let $\bar{\lambda}$ be the map sending any saturated chain $c$ of length $k \geq
1$ of $\Poset$ to the word on $\Lambda$ of length $k - 1$ defined by
\begin{equation}
    \bar{\lambda} \App c \App i
    :=
    \lambda \App \ULine{c \App i} \App \ULine{c \App \ULine{i + 1}}
\end{equation}
for any $i \in [k - 1]$. A saturated chain of $\Poset$ is \Def{$\lambda$-increasing} (resp.\
\Def{$\lambda$-weakly decreasing}) if its image by $\bar{\lambda}$ is an increasing (resp.\
weakly decreasing) word w.r.t.\ the partial order relation $\Leq_{\Lambda}$. A saturated
chain $c$ of $\Poset$ is \Def{$\lambda$-smaller} than a saturated chain $c'$ of $\Poset$ if
$\bar{\lambda} \App c$ is smaller than $\bar{\lambda} \App c'$ for the lexicographic order
induced by $\Leq_{\Lambda}$. The labeling $\lambda$ is an \Def{EL-labeling} of $\Poset$ if
for any $x, x' \in \Poset$ satisfying $x \Leq_\Poset x'$, there is exactly one
$\lambda$-increasing saturated chain $c$ from $x$ to $x'$, and $c$ is $\lambda$-smaller than
any other saturated chains from $x$ to~$x'$.

Let us denote by $\Z^3$ the set of triples of integers endowed with the lexicographic order.
Let $\EasterlyWindEdgeLabeling : \; \EasterlyWindCovering \; \to \Z^3$ be the map defined,
for any $\Par{\TermT_1, \TermT_2} \in \; \EasterlyWindCovering$, by
\begin{equation}
    \EasterlyWindEdgeLabeling \App \TermT_1 \App \TermT_2 := \Par{i, i_1, -j_1}
\end{equation}
where $\Par{i_1, j_1, i}$ is the edge of $\TermT_1$ which is replaced by an edge $\Par{i_1',
j_1', i}$ in order to produce $\TermT_2$. This map is well-defined thanks to
Lemma~\ref{lem:covering_edge_properties}. For instance, by considering the
$\SignatureExample$-term $\TermT_1$ as the left-hand side of~\eqref{equ:example_covering_1}
and the $\SignatureExample$-term $\TermT_2$ as its right-hand side, we have
\begin{math}
    \EasterlyWindEdgeLabeling \App \TermT_1 \App \TermT_2 = \Par{3, 2, -2}.
\end{math}
With the same conventions, in~\eqref{equ:example_covering_2}, we have
\begin{math}
    \EasterlyWindEdgeLabeling \App \TermT_1 \App \TermT_2 = \Par{5, 1, -2}.
\end{math}

\begin{Statement}{Lemma}{lem:order_preserving_labeling}
    Let $\Signature$ be a signature, and $\TermT_1$, $\TermT_2$, and $\TermT_2'$ be three
    $\Signature$-terms such that $\TermT_1 \EasterlyWindCovering_{i_1} \TermT_2$ and
    $\TermT_1 \EasterlyWindCovering_{i_1'} \TermT_2'$, where $i_1$ and $i_1'$ are internal
    nodes of $\TermT_1$. We have $i_1 \leq i_1'$ if and only if
    \begin{math}
        \EasterlyWindEdgeLabeling \App \TermT_1 \App \TermT_2
    \end{math}
    is smaller than or equal to
    \begin{math}
        \EasterlyWindEdgeLabeling \App \TermT_1 \App \TermT_2'
    \end{math}
    for the lexicographic order.
\end{Statement}
\begin{Proof}
    By Lemma~\ref{lem:covering_edge_properties}, $\TermT_2$ is obtained from $\TermT_1$ by
    replacing the edge $\Par{i, j, i_1}$ by the edge $\Par{i_2, j_2, i_1}$ where $i$ and
    $i_2$ are internal nodes of both $\TermT_1$ and $\TermT_2$, and $j$ and $j_2$ are
    integers. In the same way, $\TermT_2'$ is obtained from $\TermT_1$ by replacing the edge
    $\Par{i', j', i_1'}$ by an edge $\Par{i_2', j_2', i_1'}$ where $i'$ and $i_2'$ are
    internal nodes of both $\TermT_1$ and $\TermT_2'$, and $j'$ and $j_2'$ are integers. By
    definition of $\EasterlyWindEdgeLabeling$,
    \begin{math}
        \EasterlyWindEdgeLabeling \App \TermT_1 \App \TermT_2 = \Par{i_1, i, -j}
    \end{math}
    and
    \begin{math}
        \EasterlyWindEdgeLabeling \App \TermT_1 \App \TermT_2' = \Par{i_1', i', -j'}.
    \end{math}
    The statement of the lemma follows immediately.
\end{Proof}

Observe that Lemma~\ref{lem:order_preserving_labeling} implies that if $u$ and $u'$ are two
$\Par{\TermT_1, \TermT_2}$-sequences where $\TermT_1$ and $\TermT_2$ are two
$\Signature$-terms satisfying $\TermT_1 \EasterlyWindLeq \TermT_2$ and $u$ is
lexicographically smaller than $u'$, then $\bar{\EasterlyWindEdgeLabeling} \App c$ is
smaller than $\bar{\EasterlyWindEdgeLabeling} \App c'$, where $c$ (resp.\ $c'$) is the
saturated chain specified by $u$ (resp.\ $u'$).

\begin{Statement}{Lemma}{lem:small_saturated_chain}
    Let $\Signature$ be a signature, and $\TermT_1$, $\TermT_2$, and $\TermT_3$ be three
    $\Signature$-terms such that
    \begin{math}
        \TermT_1
        \EasterlyWindCovering_{i_1} \TermT_2
        \EasterlyWindCovering_{i_2} \TermT_3,
    \end{math}
    where $i_1$ is an internal node of $\TermT_1$ and $i_2$ is an internal node of
    $\TermT_2$. The following properties hold:
    \begin{enumerate}[label=(\roman*)]
        \item \label{item:small_saturated_chain_1}
        the triples
        \begin{math}
            \EasterlyWindEdgeLabeling \App \TermT_1 \App \TermT_2
        \end{math}
        and
        \begin{math}
            \EasterlyWindEdgeLabeling \App \TermT_2 \App \TermT_3
        \end{math}
        are different;
        \item \label{item:small_saturated_chain_2}
        we have $i_1 \leq i_2$ if and only if
        \begin{math}
            \EasterlyWindEdgeLabeling \App \TermT_1 \App \TermT_2
        \end{math}
        is smaller than
        \begin{math}
            \EasterlyWindEdgeLabeling \App \TermT_2 \App \TermT_3
        \end{math}
        for the lexicographic order.
    \end{enumerate}
\end{Statement}
\begin{Proof}
    Let us first gather some properties about the fact that
    \begin{math}
        \TermT_1
        \EasterlyWindCovering_{i_1} \TermT_2
        \EasterlyWindCovering_{i_2} \TermT_3.
    \end{math}
    By Lemma~\ref{lem:covering_edge_properties}, $\TermT_2$ is obtained from $\TermT_1$ by
    replacing the edge $\Par{i, j, i_1}$ by the edge $\Par{i', j', i_1}$ where $i$ and $i'$
    are internal nodes of both $\TermT_1$ and $\TermT_2$, and $j$ and $j'$ are integers. In
    the same way, $\TermT_3$ is obtained from $\TermT_2$ by replacing the edge $\Par{i'',
    j'', i_2}$ by the edge $\Par{i''', j''', i_2}$ where $i''$ and $i'''$ are internal nodes
    of both $\TermT_2$ and $\TermT_3$, and $j''$ and $j'''$ are integers. Again by
    Lemma~\ref{lem:covering_edge_properties}, the edge $\Par{i, j, i_1}$ is dominated by the
    edge $\Par{i', j', i_1}$, and the edge $\Par{i'', j'', i_2}$ is dominated by the edge
    $\Par{i''', j''', i_2}$. Moreover, by definition of $\EasterlyWindEdgeLabeling$, we have
    \begin{math}
        \EasterlyWindEdgeLabeling \App \TermT_1 \App \TermT_2 = \Par{i_1, i, -j}
    \end{math}
    and
    \begin{math}
        \EasterlyWindEdgeLabeling \App \TermT_2 \App \TermT_3 = \Par{i_2, i'', -j''}.
    \end{math}

    First, let us assume by contradiction that
    \begin{math}
        \EasterlyWindEdgeLabeling \App \TermT_1 \App \TermT_2
        =
        \EasterlyWindEdgeLabeling \App \TermT_2 \App \TermT_3.
    \end{math}
    This implies that $\Par{i_1, i, -j} = \Par{i_2, i'', -j''}$ so that $\TermT_1$ and
    $\TermT_2$ have both the same edge $\Par{i, j, i_1} = \Par{i'', j'', i_2}$. Moreover,
    since $\TermT_1 \EasterlyWindCovering_{i_1} \TermT_2$, by
    Lemma~\ref{lem:covering_edge_properties}, $\Par{i, j, i_1}$ is not an edge of
    $\TermT_2$. This yields a contradiction and \ref{item:small_saturated_chain_1} checks
    out.

    Besides, if $i_1 < i_2$, then we have immediately that
    \begin{math}
        \EasterlyWindEdgeLabeling \App \TermT_1 \App \TermT_2
    \end{math}
    is smaller than
    \begin{math}
        \EasterlyWindEdgeLabeling \App \TermT_2 \App \TermT_3
    \end{math}
    for the lexicographic order. When $i_1 = i_2$, the edges $\Par{i', j', i_1}$ and
    $\Par{i'', j'', i_2}$ are equal. Hence, $\Par{i, j, i_1}$ differs from and is dominated
    by $\Par{i'', j'', i_2}$. For this reason,
    \begin{math}
        \EasterlyWindEdgeLabeling \App \TermT_1 \App \TermT_2
    \end{math}
    is smaller than
    \begin{math}
        \EasterlyWindEdgeLabeling \App \TermT_2 \App \TermT_3.
    \end{math}
    for the lexicographic order. Conversely, when $i_1 > i_2$, it follows directly that
    \begin{math}
        \EasterlyWindEdgeLabeling \App \TermT_1 \App \TermT_2
    \end{math}
    is greater than
    \begin{math}
        \EasterlyWindEdgeLabeling \App \TermT_2 \App \TermT_3.
    \end{math}
    Therefore, the fact that
    \begin{math}
        \EasterlyWindEdgeLabeling \App \TermT_1 \App \TermT_2
    \end{math}
    is smaller than
    \begin{math}
        \EasterlyWindEdgeLabeling \App \TermT_2 \App \TermT_3
    \end{math}
    implies that $i_1 \leq i_2$. The equivalence stated
    by~\ref{item:small_saturated_chain_2} is established.
\end{Proof}

\begin{Statement}{Theorem}{thm:el_shellability}
    For any signature $\Signature$,
    \begin{enumerate}[label=(\roman*)]
        \item \label{item:el_shellability_1}
        the labeling $\EasterlyWindEdgeLabeling$ is an EL-labeling of the
        $\Signature$-easterly wind poset;
        \item \label{item:el_shellability_2}
        there is at most one $\EasterlyWindEdgeLabeling$-weakly decreasing saturated chain
        between any pair of elements of the $\Signature$-easterly wind poset.
    \end{enumerate}
\end{Statement}
\begin{Proof}
    First of all, by Theorem~\ref{thm:easterly_wind_posets}, since $\EasterlyWindCovering$
    is the covering relation of the $\Signature$-easterly wind poset, the map
    $\EasterlyWindEdgeLabeling$ is a well-defined labeling of this poset.

    Let $\TermT_1$ and $\TermT_2$ be two $\Signature$-terms such that $\TermT_1
    \EasterlyWindLeq \TermT_2$. Assume that there exist two
    $\EasterlyWindEdgeLabeling$-weakly increasing (resp.\ $\EasterlyWindEdgeLabeling$-weakly
    decreasing) saturated chains $c$ and $c'$ between $\TermT_1$ and $\TermT_2$. By
    Lemma~\ref{lem:sequence_specification}, there exist two $\Par{\TermT_1,
    \TermT_2}$-sequences $u$ and $u'$ such that $c$ is specified by $u$ and $c'$ is
    specified by $u'$. By Lemma~\ref{lem:order_preserving_labeling}, $u$ and $u'$ are weakly
    increasing (resp.\ weakly decreasing) words. By
    Point~\ref{item:unique_sequences_specification_1} (resp.\
    Point~\ref{item:unique_sequences_specification_2}) of
    Lemma~\ref{lem:unique_sequences_specification}, we have $u = u'$. This shows that $c$
    and $c'$ are in fact the same saturated chain. Therefore, there is at most one
    $\EasterlyWindEdgeLabeling$-weakly increasing (resp.\ $\EasterlyWindEdgeLabeling$-weakly
    decreasing) saturated chain from $\TermT_1$ to $\TermT_2$. In particular, this
    proves~\ref{item:el_shellability_2}.

    Let $\TermT_1$ and $\TermT_2$ be two $\Signature$-terms such that $\TermT_1
    \EasterlyWindLeq \TermT_2$, and let $u$ be a weakly increasing $\Par{\TermT_1,
    \TermT_2}$-sequence. The existence of such $u$ is ensured by
    Point~\ref{item:unique_sequences_specification_3} of
    Lemma~\ref{lem:unique_sequences_specification}. By
    Lemma~\ref{lem:small_saturated_chain}, the saturated chain specified by $u$ is
    $\EasterlyWindEdgeLabeling$-increasing. The uniqueness of this
    $\EasterlyWindEdgeLabeling$-increasing saturated chain is shown in the previous
    paragraph of this proof. Finally, the fact that the saturated chain specified by $u$ is
    $\EasterlyWindEdgeLabeling$-smaller than any other saturated chains from $\TermT_1$ to
    $\TermT_2$ follows from Point~\ref{item:unique_sequences_specification_4} of
    Lemma~\ref{lem:unique_sequences_specification} and from
    Lemma~\ref{lem:order_preserving_labeling}. Therefore, \ref{item:el_shellability_1} is
    established.
\end{Proof}

An important consequence~\cite{BW96} of Point~\ref{item:el_shellability_2} of
Theorem~\ref{thm:el_shellability} is that the codomain of the Möbius function $\Mobius$ of
the $\Signature$-easterly wind poset is the set $\Bra{-1, 0, 1}$.

Besides, since by Theorem~\ref{thm:tilt_closure}, for any set $\SpecialNodes$ of positive
integers, $\Tilt \App \SpecialNodes$ is a closure operator of the $\Signature$-easterly wind
poset, by the Crapo's Closure Theorem~\cite{Cra66}, the Möbius function $\Mobius \App
\SpecialNodes$ of the $\SpecialNodes$-tilted $\Signature$-easterly wind poset satisfies
\begin{equation}
    \Mobius \App \SpecialNodes \App \TermT_1 \App \TermT_2
    =
    \sum_{\TermT_2' \in \SetTerms \App \Signature}
    \Iverson{\Tilt \App \SpecialNodes \App \TermT_2' = \TermT_2}
    \;
    \Mobius \App \TermT_1 \App \TermT_2'
\end{equation}
for all $\SpecialNodes$-tilted $\Signature$-terms $\TermT_1$ and $\TermT_2$.

\subsubsection{Geometric realization}
Let $\Poset$ be an interval of the $\SpecialNodes$-tilted $\Signature$-easterly wind poset.
Since two $\Signature$-terms are comparable only if they have the same degree, let us denote
by $n$ the common degree of the $\Signature$-terms of $\Poset$. The \Def{geometric
realization} $\GeometricRealization \App \Poset$ of $\Poset$ is the embedding of the Hasse
diagram of $\Poset$ in the space $\R^n$ such that each $\TermT \in \Poset$ gives rise to a
vertex of coordinates $\Connection \App \TermT$ and each pair $\Par{\TermT_1, \TermT_2}$ of
$\Signature$-terms of $\Poset$ gives rise to an edge, provided that $\TermT_1$ is covered by
$\TermT_2$ in $\Poset$.

Moreover, when $\SpecialNodes = \emptyset$, since by
Lemma~\ref{lem:covering_edge_properties} and Theorem~\ref{thm:easterly_wind_posets}, the
connection sequences of two $\Signature$-terms which are in relation for
$\EasterlyWindCovering$ differ in exactly one component, every edge of
$\GeometricRealization \App \Poset$ is parallel to a line passing through the origin and a
point of $\R^n$ the form $\Par{0, \dots, 0, 1, 0, \dots, 0}$. For this reason, the geometric
realizations of $\Signature$-easterly wind posets are cubic~\cite{CG22,Com23}. In general,
$\GeometricRealization \App \Poset$ is not cubic when $\Poset$ is an interval of an
$\SpecialNodes$-tilted $\Signature$-easterly wind poset with $\SpecialNodes \ne \emptyset$.
Figure~\ref{fig:example_geometric_realization} shows examples of geometric realizations of
such intervals.
\begin{figure}[ht]
    \centering
    \def\Zoom{0.4}
\def\Scale{5.0}
\scalebox{\Zoom}{
    \begin{tikzpicture}[Centering,scale=\Scale,
    x={(0,-1cm)},y={(-1cm,-1cm)}, z={(1cm,-1cm)}]
        \DrawCubeGrid{1.50}{1.75}{2.00}{2.75}{1.00}{3.00}{0.25}
        \node[GraphVertex](1)at(1.50,2.00,1.00){};
        \node[GraphVertex](2)at(1.50,2.00,3.00){};
        \node[GraphVertex](3)at(1.50,2.50,1.00){};
        \node[GraphVertex](4)at(1.50,2.50,2.00){};
        \node[GraphVertex](5)at(1.50,2.50,3.00){};
        \node[GraphVertex](6)at(1.50,2.75,1.00){};
        \node[GraphVertex](7)at(1.50,2.75,2.00){};
        \node[GraphVertex](8)at(1.50,2.75,2.50){};
        \node[GraphVertex](9)at(1.50,2.75,3.00){};
        \node[GraphVertex](10)at(1.75,2.00,1.00){};
        \node[GraphVertex](11)at(1.75,2.00,1.50){};
        \node[GraphVertex](12)at(1.75,2.00,3.00){};
        \node[GraphVertex](13)at(1.75,2.50,1.00){};
        \node[GraphVertex](14)at(1.75,2.50,1.50){};
        \node[GraphVertex](15)at(1.75,2.50,2.00){};
        \node[GraphVertex](16)at(1.75,2.50,3.00){};
        \node[GraphVertex](17)at(1.75,2.75,1.00){};
        \node[GraphVertex](18)at(1.75,2.75,1.50){};
        \node[GraphVertex](19)at(1.75,2.75,2.00){};
        \node[GraphVertex](20)at(1.75,2.75,2.50){};
        \node[GraphVertex](21)at(1.75,2.75,3.00){};
        \draw[GraphArc](1)--(2);
        \draw[GraphArc](1)--(3);
        \draw[GraphArc](1)--(10);
        \draw[GraphArc](2)--(5);
        \draw[GraphArc](2)--(12);
        \draw[GraphArc](3)--(4);
        \draw[GraphArc](3)--(6);
        \draw[GraphArc](3)--(13);
        \draw[GraphArc](4)--(5);
        \draw[GraphArc](4)--(7);
        \draw[GraphArc](4)--(15);
        \draw[GraphArc](5)--(9);
        \draw[GraphArc](5)--(16);
        \draw[GraphArc](6)--(7);
        \draw[GraphArc](6)--(17);
        \draw[GraphArc](7)--(8);
        \draw[GraphArc](7)--(19);
        \draw[GraphArc](8)--(9);
        \draw[GraphArc](8)--(20);
        \draw[GraphArc](9)--(21);
        \draw[GraphArc](10)--(11);
        \draw[GraphArc](10)--(13);
        \draw[GraphArc](11)--(12);
        \draw[GraphArc](11)--(14);
        \draw[GraphArc](12)--(16);
        \draw[GraphArc](13)--(14);
        \draw[GraphArc](13)--(17);
        \draw[GraphArc](14)--(15);
        \draw[GraphArc](14)--(18);
        \draw[GraphArc](15)--(16);
        \draw[GraphArc](15)--(19);
        \draw[GraphArc](16)--(21);
        \draw[GraphArc](17)--(18);
        \draw[GraphArc](18)--(19);
        \draw[GraphArc](19)--(20);
        \draw[GraphArc](20)--(21);
    \end{tikzpicture}
}
    \hspace{6em}
    \def\Zoom{0.4}
\def\Scale{5.0}
\scalebox{\Zoom}{
    \begin{tikzpicture}[Centering,scale=\Scale, x={(-1cm,-1cm)},y={(1cm,-1cm)}]
        \DrawSquareGrid{2.00}{2.75}{1.50}{3.00}{0.25}
        \node[GraphVertex](1)at(2.00,1.50){};
        \node[GraphVertex](2)at(2.00,3.00){};
        \node[GraphVertex](3)at(2.50,1.50){};
        \node[GraphVertex](4)at(2.50,2.00){};
        \node[GraphVertex](5)at(2.50,3.00){};
        \node[GraphVertex](6)at(2.75,1.50){};
        \node[GraphVertex](7)at(2.75,2.00){};
        \node[GraphVertex](8)at(2.75,2.50){};
        \node[GraphVertex](9)at(2.75,3.00){};
        \draw[GraphArc](1)--(2);
        \draw[GraphArc](1)--(3);
        \draw[GraphArc](2)--(5);
        \draw[GraphArc](3)--(4);
        \draw[GraphArc](3)--(6);
        \draw[GraphArc](4)--(5);
        \draw[GraphArc](4)--(7);
        \draw[GraphArc](5)--(9);
        \draw[GraphArc](6)--(7);
        \draw[GraphArc](7)--(8);
        \draw[GraphArc](8)--(9);
    \end{tikzpicture}
}
    \caption[]{
        The geometric realization of the $\SignatureExample$-easterly wind poset of $\TermT$
        on the left, and the geometric realization of the $\Bra{1}$-tilted
        $\SignatureExample$-easterly wind poset of $\Tilt \App \Bra{1} \App \TermT$ on the
        right, for $\TermT := \GenA_3 \LeafSimple \ULine{\GenA_3 \LeafSimple \LeafSimple
        \ULine{\GenA_1 \LeafSimple}} \ULine{\GenA_1 \LeafSimple}$. The top (resp.\ bottom)
        point has $\Par{\frac{15}{18}, \frac{3}{2}, 2, 1}$ (resp.\ $\Par{\frac{15}{18},
        \frac{7}{4}, \frac{11}{4}, 3}$) as coordinates on the left, and the top (resp.\
        bottom) point has $\Par{\frac{15}{18}, \frac{7}{4}, 2, \frac{3}{2}}$ (resp.\
        $\Par{\frac{15}{18}, \frac{7}{4}, \frac{11}{4}, 3}$) on the right.
    }
    \label{fig:example_geometric_realization}
\end{figure}

\subsection{Lattice properties}
The purpose of this section is to show that all terminal intervals of $\Signature$-easterly
wind posets are lattices. We also show that this property holds for $\SpecialNodes$-tilted
$\Signature$-easterly wind posets.

\subsubsection{Join-semilattice structure} \label{subsubsec:join_semilattice}
For any words $w_1$ and $w_2$ on $\Q$ of the same length $n$, let $w_3 := w_1 \JJoin w_2$ be
the sequence of length $n$ such that for any $i \in [n]$, $w_3 \App i = \max \Bra{w_1 \App
i, w_2 \App i}$.

Let $\TermT$ be an $\Signature$-term and let on the set $\TerminalInterval \App \TermT$ the
binary operation $\JJoin$ defined as follows. For any $\Signature$-terms $\TermT_1$ and
$\TermT_2$ of $\TerminalInterval \App \TermT$, let $\TermT_1 \JJoin \TermT_2$ be the
$\Signature$-term $\TermT_3$ having $\Decoration \App \TermT$ as decoration word and such
that its connection word is $\Connection \App \TermT_1 \JJoin \Connection \App \TermT_2$.
For instance, for
\begin{math}
    \TermT
    := \GenA_3 \LeafSimple \ULine{\GenA_3 \LeafSimple \LeafSimple \GenA_0} \ULine{\GenA_3
    \LeafSimple \LeafSimple \LeafSimple},
\end{math}
we have
\begin{equation}
    \GenA_3 \LeafSimple \ULine{
        \GenA_3 \LeafSimple \GenA_0 \LeafSimple
    }
    \ULine{\GenA_3 \LeafSimple \LeafSimple \LeafSimple}
    \JJoin
    \GenA_3 \ULine{
        \GenA_3 \LeafSimple \LeafSimple \GenA_0
    }
    \ULine{\GenA_3 \LeafSimple \LeafSimple \LeafSimple} \LeafSimple
    = \GenA_3 \ULine{
        \GenA_3 \LeafSimple \GenA_0 \LeafSimple
    }
    \ULine{\GenA_3 \LeafSimple \LeafSimple \LeafSimple} \LeafSimple.
\end{equation}
Note that by Lemma~\ref{lem:covering_edge_properties}, since $\TerminalInterval \App \TermT
\EasterlyWindLeq \TermT_1$ and $\TerminalInterval \App \TermT \EasterlyWindLeq \TermT_2$, we
have $\Decoration \App \TermT_1 = \Decoration \App \TermT = \Decoration \App \TermT_2$.

\begin{Statement}{Lemma}{lem:extreme_leaves_preservation}
    Let $\Signature$ be a signature and $\TermT_1$ and $\TermT_2$ be two $\Signature$-terms
    such that $\TermT_1 \EasterlyWindLeq \TermT_2$. If $\TermT_1$ has an internal node $i$
    such that its $j$-th child is an extreme leaf, then in $\TermT_2$, the $j$-th child of
    $i$ is an extreme leaf.
\end{Statement}
\begin{Proof}
    Assume that $\TermT_1 \EasterlyWindCovering \TermT_2$ and that $\TermT_1$ has an
    internal node $i$ such that its $j$-th child is an extreme leaf. Let us call $k$ this
    leaf. Since $k$ is extreme in $\TermT_1$, there is no internal node of $\TermT_1$ which
    is visited after $k$ in the preorder traversal of $\TermT_1$. For this reason,
    $\TermT_2$ cannot be obtained by replacing $k$ by any other subterm through the
    $\Signature$-easterly wind rewrite rule. This shows that $k$ remains an extreme leaf in
    $\TermT_2$, so that the $j$-th child of $i$ is an extreme leaf. Now, the fact that, by
    Theorem~\ref{thm:easterly_wind_posets}, $\EasterlyWindCovering$ is the covering relation
    of the partial order relation $\EasterlyWindLeq$ entails the statement of the lemma.
\end{Proof}

\begin{Statement}{Proposition}{prop:join_operation}
    For any signature $\Signature$ and any $\Signature$-term $\TermT$, the operation
    $\JJoin$ is well-defined on the $\Signature$-easterly wind poset of $\TermT$.
\end{Statement}
\begin{Proof}
    We have to show that for any $\Signature$-terms $\TermT_1$ and $\TermT_2$ such that
    $\TermT \EasterlyWindLeq \TermT_1$ and $\TermT \EasterlyWindLeq \TermT_2$, the word $c_1
    \JJoin c_2$ is the connection word of an $\Signature$-term of $\TerminalInterval \App
    \TermT$, having $\Decoration \App \TermT$ as decoration word, where $c_1 := \Connection
    \App \TermT_1$ and $c_2 := \Connection \App \TermT_2$. Let us prove this property by
    induction on $n$, the degree of $\TermT$. This is immediately true when $n = 0$. Assume
    that $n \geq 1$. By Lemma~\ref{lem:partial_order_contraction}, we have $\Contraction
    \App \TermT \EasterlyWindLeq \Contraction \App \TermT_1$ and $\Contraction \App \TermT
    \EasterlyWindLeq \Contraction \App \TermT_2$. Let us set $c_1' := \Connection \App
    \ULine{\Contraction \App \TermT_1}$, $c_2' := \Connection \App \ULine{\Contraction \App
    \TermT_2}$, and $c' := c_1' \JJoin c_2'$. By induction hypothesis, $c'$ is the
    connection word of an $\Signature$-term $\TermT'$ of $\TerminalInterval \App
    \ULine{\Contraction \App \TermT}$ whose decoration word is $\Decoration \App
    \ULine{\Contraction \App \TermT}$. Without loss of generality, assume that $c_1 \App n
    \geq c_2 \App n$ and let $\Par{i_1, j_1, n}$ be the parent edge of $n$ in $\TermT_1$.
    Let $\TermT''$ be the $\Signature$-term obtained by adding to $\TermT'$ the edge
    $\Par{i_1, j_1, n}$ so that the added internal node $n$ is decorated by $\Decoration
    \App \TermT \App n$. Note that since $n$ is the last visited internal node of $\TermT_1$
    in the preorder traversal of $\TermT_1$, in $\Contraction \App \TermT_1$, the $j_1$-th
    child of the internal node $i_1$ is an extreme leaf. Therefore, by
    Lemma~\ref{lem:extreme_leaves_preservation}, in $\TermT'$, the $j_1$-child of $\TermT'$
    is a leaf. This ensures that it is possible to build $\TermT''$ as stated. By
    construction, $\Connection \App \TermT'' = \Connection \App \TermT' \Conc c_1 \App n$.
    Now, since $\Contraction \App \TermT \EasterlyWindLeq \TermT'$ and $\Connection \App
    \TermT \App n \leq \Connection \App \TermT_1 \App n$, we have $\TermT \EasterlyWindLeq
    \TermT''$. Moreover, by construction, $\Decoration \App \TermT'' = \Decoration \App
    \TermT' \Conc \Decoration \App \TermT \App n = \Decoration \App \TermT$. This shows the
    stated property.
\end{Proof}

\subsubsection{Lattice structure}
Let us now state one of the most important results of this section.

\begin{Statement}{Theorem}{thm:easterly_wind_lattices}
    For any signature $\Signature$ and any $\Signature$-term $\TermT$, the subposet
    $\TerminalInterval \App \TermT$ of the $\Signature$-easterly wind poset is a lattice.
    Moreover, this lattice admits $\JJoin$ as join operation.
\end{Statement}
\begin{Proof}
    Let us first prove that $\JJoin$ is the join operation of the poset
    $\Par{\TerminalInterval \App \TermT, \EasterlyWindLeq}$. First of all, by
    Proposition~\ref{prop:join_operation}, the operation $\JJoin$ is well-defined on
    $\TerminalInterval \App \TermT$. Let us show that for any $\Signature$-terms $\TermT_1$
    and $\TermT_2$ such that $\TermT \EasterlyWindLeq \TermT_1$ and $\TermT \EasterlyWindLeq
    \TermT_2$, $\TermT' := \TermT_1 \JJoin \TermT_2$ is the unique minimal element of the
    set $\TerminalInterval \App \TermT_1 \cap \TerminalInterval \App \TermT_2$. For this,
    let $\TermT''$ be an $\Signature$-term such that $\TermT_1 \EasterlyWindLeq \TermT''$
    and $\TermT_2 \EasterlyWindLeq \TermT''$. By definition of $\EasterlyWindLeq$, for any
    $i \in [\Deg \App \TermT]$, $\Connection \App \TermT'' \App i \geq \Connection \App
    \TermT_1 \App i$ and $\Connection \App \TermT'' \App i \geq \Connection \App \TermT_2
    \App i$ so that $\Connection \App \TermT'' \App i \geq \max \Bra{\Connection \App
    \TermT_1 \App i, \Connection \App \TermT_2 \App i}$. Since, by definition of $\JJoin$,
    $\Connection \App \TermT' \App i = \max \Bra{\Connection \App \TermT_1 \App i,
    \Connection \App \TermT_2 \App i}$, we have $\TermT' \EasterlyWindLeq \TermT''$. This
    shows that  $\JJoin$ is the join operation of $\Par{\TerminalInterval \App \TermT,
    \EasterlyWindLeq}$. Finally, since any join-semilattice with a unique minimal element is
    a lattice~\cite{Sta11}, the stated property holds.
\end{Proof}

A lattice $\Lattice$ with meet operation $\Meet$ and join operation $\JJoin$ is \Def{join
semi-distributive}~\cite{FJN95} if for any $x, x_1, x_2 \in \Lattice$, $x_1 \JJoin x = x
\JJoin x_2$ implies $x \JJoin \Par{x_1 \Meet x_2} = x_1 \JJoin x$.

For some $\Signature$-terms $\TermT$, the $\Signature$-easterly wind posets of $\TermT$ are
not join semi-distributive lattices. Indeed, let us consider the $\SignatureExample$-terms
\begin{enumerate}[label=$\bullet$,before=\begin{multicols}{2},after=\end{multicols}]
    \item
    \begin{math}
        \TermT
        := \GenA_4 \ULine{\GenA_1 \LeafSimple} \ULine{\GenA_1 \LeafSimple} \ULine{\GenA_1
        \LeafSimple} \ULine{\GenA_1 \LeafSimple};
    \end{math}
    \item
    \begin{math}
        \TermS
        := \GenA_4 \ULine{\GenA_1 \LeafSimple} \ULine{\GenA_1 \ULine{\GenA_1 \ULine{\GenA_1
        \LeafSimple}}} \LeafSimple \LeafSimple;
    \end{math}
    \item
    \begin{math}
        \TermS_1
        := \GenA_4 \ULine{\GenA_1 \ULine{\GenA_1 \ULine{\GenA_1 \LeafSimple}}} \LeafSimple
        \ULine{\GenA_1 \LeafSimple} \LeafSimple;
    \end{math}
    \item
    \begin{math}
        \TermS_2 := \GenA_4 \ULine{\GenA_1 \ULine{\GenA_1 \LeafSimple}} \ULine{\GenA_1
        \ULine{\GenA_1 \LeafSimple}} \LeafSimple \LeafSimple;
    \end{math}
    \item
    \begin{math}
        \TermS'
        :=
        \GenA_4 \ULine{\GenA_1 \ULine{\GenA_1 \ULine{\GenA_1 \ULine{\GenA_1 \LeafSimple}}}}
        \LeafSimple \LeafSimple \LeafSimple.
    \end{math}
\end{enumerate}
It is easy to check that $\TermS$, $\TermS_1$, and $\TermS_2$ belong to the
$\SignatureExample$-easterly wind poset of $\TermT$ and that we have
\begin{math}
    \TermS_1 \JJoin \TermS = \TermS' = \TermS \JJoin \TermS_2.
\end{math}
Now, since we have
\begin{math}
    \TermS \JJoin \Par{\TermS_1 \Meet \TermS_2}
    = \TermS \JJoin \TermT = \TermS \ne \TermS',
\end{math}
this yields a contradiction with the required relation to be join semi-distributive.

An equivalence relation $\Equiv$ on a lattice $\Lattice$ is a \Def{lattice
congruence}~\cite{CS98,Rea04} of $\Lattice$ if each $\Equiv$-equivalence class is an
interval of $\Lattice$ and both the maps sending each $x \in \Lattice$ to the smallest or
greatest element of $[x]_{\Equiv}$ are order-preserving. Observe that despite the fact that
by Proposition~\ref{prop:kernel_tilt_intervals}, for any set $\SpecialNodes$ of positive
integers, each $\Equiv_{\Tilt \App \SpecialNodes}$-equivalence class is an interval of the
$\Signature$-easterly wind poset, by the remark stated at the end of
Section~\ref{subsubsec:closure_operator}, $\TiltReverse \App \SpecialNodes$ is not
order-preserving. For this reason, for any $\Signature$-term $\TermT$, the restriction of
the equivalence relation $\Equiv_{\Tilt \App \SpecialNodes}$ on the lattice
$\TerminalInterval \App \TermT$ is not in general a lattice congruence of this lattice.

\subsubsection{Lattice structure on tilted terms}
Let $\TermT$ be an $\Signature$-term and $\SpecialNodes$ be a set of positive integers. Let
on the set $\TerminalInterval \App \SpecialNodes \App \TermT$ the binary operation $\JJoin
\App \SpecialNodes$ defined for any $\SpecialNodes$-tilted $\Signature$-terms $\TermT_1$ and
$\TermT_2$ by
\begin{math}
    \TermT_1 \JJoin \App \SpecialNodes \; \TermT_2
    := \Tilt \App \SpecialNodes \App \ULine{\TermT_1 \JJoin \TermT_2}
\end{math}
where $\JJoin$ is the operation defined in Section~\ref{subsubsec:join_semilattice}.

\begin{Statement}{Proposition}{prop:tilted_easterly_wind_lattices}
    For any signature $\Signature$, any set $\SpecialNodes$ of positive integers, and any
    $\SpecialNodes$-tilted $\Signature$-term $\TermT$, the subposet
    $\TerminalInterval \App \SpecialNodes \App \TermT$ of the $\Signature$-easterly wind
    poset is a lattice. Moreover, this lattice admits $\JJoin \App \SpecialNodes$ as join
    operation.
\end{Statement}
\begin{Proof}
    This is a consequence of the fact that by Theorem~\ref{thm:tilt_closure}, $\Tilt \App
    \SpecialNodes$ is a closure operator of $\TerminalInterval \App \TermT$ and the fact
    that, by Theorem~\ref{thm:easterly_wind_lattices}, $\TerminalInterval \App \TermT$ is a
    lattice. Indeed, as exposed in~\cite{DP02}, if $\phi$ is a closure operator of a lattice
    $\Lattice$, then $\phi \App \Lattice$ is also a lattice. This lattice has the same meet
    operation as the one of $\Lattice$, and admits the join operation $\JJoin'$ satisfying
    $x_1 \JJoin' x_2 = \phi \App \ULine{x_1 \JJoin x_2}$ where $\JJoin$ is the join
    operation of $\Lattice$, for any $x_1, x_2 \in \Lattice$.
\end{Proof}

Remark that even if, as provided by Proposition~\ref{prop:tilted_easterly_wind_lattices},
$\TerminalInterval \App \SpecialNodes \App \TermT$ is a lattice, this lattice is not a
sublattice of $\TerminalInterval \App \TermT$. Indeed, in the $\SignatureExample$-easterly
wind poset of
\begin{math}
    \GenA_2
    \ULine{
        \GenA_2
        \ULine{
            \GenA_3
            \ULine{
                \GenA_2 \LeafSimple \LeafSimple
            }
            \ULine{
                \GenA_2 \LeafSimple \LeafSimple
            }
            \LeafSimple
        }
        \ULine{
            \GenA_2 \LeafSimple \LeafSimple
        }
    }
    \LeafSimple,
\end{math}
we have for instance
\begin{equation}
    \GenA_2
    \ULine{
        \GenA_2
        \ULine{
            \GenA_3
            \ULine{
                \GenA_2 \LeafSimple \LeafSimple
            }
            \ULine{
                \GenA_2 \LeafSimple \LeafSimple
            }
            \ULine{
                \GenA_2 \LeafSimple \LeafSimple
            }
        }
        \LeafSimple
    }
    \LeafSimple
    \JJoin
    \GenA_2
    \ULine{
        \GenA_2
        \ULine{
            \GenA_3
            \ULine{
                \GenA_2
                \LeafSimple
                \ULine{
                    \GenA_2
                    \LeafSimple
                    \LeafSimple
                }
            }
            \LeafSimple
            \LeafSimple
        }
        \ULine{
            \GenA_2 \LeafSimple \LeafSimple
        }
    }
    \LeafSimple
    =
    \GenA_2
    \ULine{
        \GenA_2
        \ULine{
            \GenA_3
            \ULine{
                \GenA_2
                \LeafSimple
                \ULine{
                    \GenA_2
                    \LeafSimple
                    \LeafSimple
                }
            }
            \LeafSimple
            \ULine{
                \GenA_2
                \LeafSimple
                \LeafSimple
            }
        }
        \LeafSimple
    }
    \LeafSimple,
\end{equation}
but even if the two operands are $\Bra{1, 3}$-tilted, the result is not.

\section{Easterly wind lattices of forests} \label{sec:easterly_wind_lattices_forests}
Easterly-wind posets are sufficiently large to contain, as subposets, several notable
structures. In particular, we introduce a notion of forests regarded as specific
$\Signature$-terms and study their posets. We then obtain, on the one hand, a family of
lattices on Fuss–Catalan objects and, on the other hand, an alternative construction of the
Tamari lattices. The section concludes with lattice structures on leaning forests and
several key notions that will later be used to describe the natural Hopf algebras of free
nonsymmetric operads.

\subsection{Forests and maximal intervals} \label{subsubsec:intervals_forests}
We introduce the notion of $\Signature$-forests, which are particular $\Signature$-terms. We
also study certain maximal intervals involving such forests in the easterly wind posets.

\subsubsection{Forests} \label{subsubsec:forests}
By seeing $\N$ as the signature such that for any $n \in \N$, $\Arity \App n := n$, let
$\Signature_\N$ be the signature $\Signature \sqcup \N$. An \Def{$\Signature$-forest} is an
$\Signature_\N$-term $\ForestF$ of the form $\ForestF = n \, \TermT_1 \dots \TermT_n$ where
$n \in \N$ and for any $i \in [n]$, $\TermT_i$ is an $\Signature$-term. For instance,
\begin{math}
    4 \,
    \GenA_0
    \ULine{\GenA_2 \ULine{\GenA_1 \LeafSimple} \LeafSimple}
    \LeafSimple
    \ULine{\GenA_1 \ULine{\GenA_3 \LeafSimple \LeafSimple \LeafSimple}}
\end{math}
is an $\SignatureExample$-forest. Moreover, the \Def{concatenation} of two
$\Signature$-forests $n \, \TermT_1 \dots \TermT_n$ and $n' \, \TermT_1' \dots \TermT_{n'}'$
is the $\Signature$-forest
\begin{equation}
    n \, \TermT_1 \dots \TermT_n \; \Conc \; n' \, \TermT_1' \dots \TermT_{n'}'
    := \ULine{n + n'} \, \TermT_1 \dots \TermT_n \TermT_1' \dots \TermT_{n'}'.
\end{equation}

An $\Signature$-forest $\ForestF = n \, \TermT_1 \dots \TermT_n$ is \Def{balanced} if $\Deg
\App \ForestF = n + 1$. For instance, $0$ and $3 \, \ULine{\GenA_2 \LeafSimple \GenA_0}
\LeafSimple \ULine{\GenA_3 \LeafSimple \LeafSimple \LeafSimple}$ are balanced
$\SignatureExample$-forests. On the contrary, $2 \, \ULine{\GenA_2 \GenA_0 \LeafSimple}
\GenA_0$ is an $\SignatureExample$-forest which is not balanced. The \Def{size} of a
balanced $\Signature$-forest $\ForestF$ is the decoration of the root of $\ForestF$
(or, equivalently, the degree of $\ForestF$ minus one). Observe that the concatenation of
two balanced $\Signature$-forests is a balanced $\Signature$-forest.

\subsubsection{Maximal intervals}
For any $\Signature$-forest $\ForestF$ of arity $1$ or more and any $\Signature$-term
$\TermT$, let $\ForestF \CompositionFirst \TermT$ be the $\Signature$-forest obtained by
replacing the leftmost leaf of $\ForestF$ by $\TermT$. For instance,
\begin{equation}
    \ColA{
        3 \,
        \ULine{\GenA_2 \GenA_0 \GenA_0}
        \ULine{
            \GenA_2
            \ULine{\GenA_1 \LeafSimple}
            \LeafSimple
        }
        \LeafSimple
    }
    \; \CompositionFirst \;
    \bm{\ColB{\GenA_3 \LeafSimple \ULine{\GenA_1 \LeafSimple} \LeafSimple}}
    =
    \ColA{
        3 \,
        \ULine{\GenA_2 \GenA_0 \GenA_0}
        \ULine{
            \GenA_2
            \ULine{
                \GenA_1
                \ULine{
                    \bm{\ColB{\GenA_3 \LeafSimple \ULine{\GenA_1 \LeafSimple} \LeafSimple}}
                }
            }
            \LeafSimple
        }
        \LeafSimple
    }.
\end{equation}
Now, given a word $w$ on $\Signature$ of length $n \in \N$, let the $\Signature$-forests
\begin{equation}
    \ForestMin \App w
    :=
    n \,
    \ULine{\Corolla \App \ULine{w \App 1}}
    \ULine{\Corolla \App \ULine{w \App 2}}
    \dots \ULine{\Corolla \App \ULine{w \App n}}
\end{equation}
and
\begin{equation}
    \ForestMax \App w
    :=
    \Par{
        \dots
        \Par{
            \Par{
                \ULine{\Corolla \App n}
                \CompositionFirst
                \Corolla \App \ULine{w \App 1}
            }
            \CompositionFirst
            \Corolla \App \ULine{w \App 2}
        }
        \dots
    }
    \CompositionFirst
    \Corolla \App \ULine{w \App n}.
\end{equation}
For instance, for $w := \GenA_2 \GenA_1 \GenA_0 \GenA_0 \GenA_0 \GenA_3 \GenA_2 \in
\SignatureExample^*$, we have
\begin{equation}
    \ForestMin \App w
    =
    7
    \ULine{\GenA_2 \LeafSimple \LeafSimple}
    \ULine{\GenA_1 \LeafSimple}
    \GenA_0 \GenA_0 \GenA_0
    \ULine{\GenA_3 \LeafSimple \LeafSimple \LeafSimple}
    \ULine{\GenA_2 \LeafSimple \LeafSimple}
\end{equation}
and
\begin{equation}
    \ForestMax \App w
    =
    7
    \ULine{\GenA_2 \ULine{\GenA_1 \GenA_0} \GenA_0}
    \GenA_0
    \ULine{\GenA_3 \ULine{\GenA_2 \LeafSimple \LeafSimple} \LeafSimple \LeafSimple}
    \LeafSimple \LeafSimple \LeafSimple \LeafSimple.
\end{equation}
Note that $\ForestMin \App \epsilon = 0 = \ForestMax \App \epsilon$. It is straightforward
to prove that both $\ForestMin \App w$ and $\ForestMax \App w$ are well-defined balanced
$\Signature$-forests of size~$\Length \App w$.

\begin{Statement}{Theorem}{thm:maximal_intervals}
    For any signature $\Signature$, any set $\SpecialNodes$ of positive integers, and any
    word $w$ on $\Signature$ of length $n \in \N$,
    \begin{enumerate}[label=(\roman*)]
        \item \label{item:maximal_intervals_1}
        in the $\SpecialNodes$-tilted $\Signature_\N$-easterly wind poset,
        $\ForestMin \App w \EasterlyWindLeq \ForestMax \App w$;
        \item \label{item:maximal_intervals_2}
        the $\Signature$-forest $\ForestMin \App w$ (resp.\ $\ForestMax \App
        w$) is a minimal (resp.\ maximal) element of the $\SpecialNodes$-tilted
        $\Signature_\N$-easterly wind poset.
    \end{enumerate}
\end{Statement}
\begin{Proof}
    Observe first that neither $\ForestMin \App w$ nor $\ForestMax \App w$ depend
    on the set $\SpecialNodes$. Indeed, directly from their definitions, it follows that
    these two $\Signature$-forests are $\SpecialNodes$-tilted for any set $\SpecialNodes$ of
    positive integers. For this reason, in this proof, we consider simply that
    $\SpecialNodes = \emptyset$.

    To prove~\ref{item:maximal_intervals_1}, let us proceed by induction on $n$. When $n =
    0$, $\ForestMin = 0 = \ForestMax$ so that the property is satisfied. Assume
    that the property holds for any word $w$ on $\Signature$ of length $n \geq 0$ and let $a
    \in \Signature$. By definition of $\ForestMin$, we have
    \begin{math}
        \Edges \App \ULine{\ForestMin \App \ULine{w \Conc a}}
        = \Edges \App \ULine{\ForestMin \App w} \cup \Bra{\Par{1, n + 1, n + 2}}.
    \end{math}
    Moreover, by definition of $\ForestMax$, we have
    \begin{math}
        \Edges \App \ULine{\ForestMax \App \ULine{w \Conc a}}
        = \Edges \App \ULine{\ForestMax \App w} \cup \Bra{\Par{i, j, n + 2}}
    \end{math}
    where $i$ is the internal node of $\ForestMax \App w$ which is the parent of the
    leftmost leaf and $j$ is the position of this leaf in its siblings. By induction
    hypothesis, we have $\ForestMin \App w \EasterlyWindLeq \ForestMax \App w$.
    Thus, by Lemma~\ref{lem:order_domination}, for any $i' \in [n + 1]$, the parent edge of
    $i'$ in $\ForestMin \App w$ is dominated by that of $i'$ in $\ForestMax \App w$.
    Moreover, the edge $\Par{1, n + 1, n + 2}$ is dominated by the edge $\Par{i, j, n + 2}$.
    Indeed, otherwise, we would have $i = 1$ and $j > n + 1$, which is absurd since the
    internal node $1$ of both $\ForestMin \App \ULine{w \Conc a}$ and $\ForestMax
    \App \ULine{w \Conc a}$ has arity $n + 1$. Therefore, again by
    Lemma~\ref{lem:order_domination}, we have $\ForestMin \App \ULine{w \Conc a}
    \EasterlyWindLeq \ForestMax \App \ULine{w \Conc a}$ as expected.

    To prove~\ref{item:maximal_intervals_2}, assume first that $\TermT$ is an
    $\Signature$-forest such that $\TermT \EasterlyWindLeq \ForestMin \App w$. Hence,
    by Lemma~\ref{lem:order_domination}, for any $i \in [2, n]$, the parent edge
    $\Par{\Parent \App \TermT \App i, \LocalPosition \App \TermT \App i, i}$ of $i$ in
    $\TermT$ is dominated by the parent edge $\Par{1, i - 1, i}$ of $i$ in $\ForestMin
    \App w$. By definition of the notion of edge domination, we have necessarily $\Parent
    \App \TermT \App i = 1$ and $\LocalPosition \App \TermT \App i \geq i - 1$. It follows
    that $\LocalPosition \App \TermT \App i = i - 1$ so that $\TermT = \ForestMin \App w$.
    This shows that $\ForestMin \App w$ is a minimal element of the $\Signature_\N$-easterly
    wind poset. Finally, observe that in $\ForestMax \App w$, all leaves are visited after
    all internal internal nodes in its preorder traversal. Therefore, there is no
    $\Signature$-forest $\TermT$ such that $\ForestMax \EasterlyWindCovering \TermT$.
    Therefore, by Lemma~\ref{lem:partial_order_sequence_rewrites}, this implies that
    $\ForestMax \App w$ is a maximal element of the $\Signature_\N$-easterly wind poset.
\end{Proof}

By Theorem~\ref{thm:maximal_intervals} and
Proposition~\ref{prop:tilted_easterly_wind_lattices}, the interval $\Han{\ForestMin \App w,
\ForestMax \App w}$ of the $\SpecialNodes$-tilted $\Signature_\N$-easterly wind poset is a
maximal interval and a lattice. Let us call it the \Def{balanced forest
$\SpecialNodes$-tilted $\Signature$-easterly wind poset of~$w$}. When $\SpecialNodes =
\emptyset$, this poset is the \Def{balanced forest $\Signature$-easterly wind poset of~$w$}.
For instance, by replacing all decorations $\GenA_3$ of the roots of the terms of
Figure~\ref{fig:example_easterly_wind_poset} by the decoration $3 \in \N$, the resulting
Hasse diagram of this figure is the balanced forest $\emptyset$-tilted
$\SignatureExample$-easterly wind poset of $\GenA_3 \GenA_1 \GenA_2$. With the same change,
Figure~\ref{fig:example_tilted_easterly_wind_poset} shows the Hasse diagram of the balanced
forest $\Bra{1, 2}$-tilted $\SignatureExample$-easterly wind lattice of~$\GenA_3 \GenA_1
\GenA_2$.

\subsection{Catalan lattices and Fuss-Catalan lattices}
The purpose of this section is to build, as a particular balanced forest
$\SpecialNodes$-tilted easterly wind posets of some words, partial order structures on the
combinatorial set of Catalan or Fuss-Catalan objects. We begin by constructing such
structures by introducing a nontrivial bijection between the set of some balanced forests
and the set of terms whose internal node all have a fixed arity. Using this bijection, we
prove that our first family of posets admits as underlying set the set of Fuss-Catalan
objects. Independently, we consider a balanced forest easterly wind poset on fully tilted
terms and show that this construction yields an alternative description of the well-known
Tamari partial order.

\subsubsection{Fuss-Catalan lattices} \label{subsubsec:fuss_catalan_lattices}
\begin{figure}[ht]
    \centering
    \def\Zoom{0.5}
\def\Scale{3}
\scalebox{\Zoom}{
    \begin{tikzpicture}[Centering,scale=\Scale,
    x={(0,-5cm)},y={(-1cm,-1cm)},z={(1cm,-1cm)}]
        \DrawCubeGrid{1.75}{2.00}{1.50}{3.00}{1.00}{4.00}{0.25}
        \node[GraphLabeledVertex](1)at(1.75,1.50,1.00){
            $4 \ULine{1 \LeafSimple} \ULine{1 \LeafSimple} \ULine{1 \LeafSimple} \ULine{1
            \LeafSimple}$
        };
        \node[GraphLabeledVertex](2)at(1.75,1.50,4.00){
            $4 \ULine{1 \LeafSimple} \ULine{1 \LeafSimple} \ULine{1 \ULine{1 \LeafSimple}}
            \LeafSimple$
        };
        \node[GraphLabeledVertex](3)at(1.75,3.00,1.00){
            $4 \ULine{1 \LeafSimple} \ULine{1 \ULine{1 \LeafSimple}} \LeafSimple \ULine{1
            \LeafSimple}$
        };
        \node[GraphLabeledVertex](4)at(1.75,3.00,1.50){
            $4 \ULine{1 \LeafSimple} \ULine{1 \ULine{1 \LeafSimple}} \ULine{1 \LeafSimple}
            \LeafSimple$
        };
        \node[GraphLabeledVertex](5)at(1.75,3.00,4.00){
            $4 \ULine{1 \LeafSimple} \ULine{1 \ULine{1 \ULine{1 \LeafSimple}}} \LeafSimple
            \LeafSimple$
        };
        \node[GraphLabeledVertex](6)at(2.00,1.50,1.00){
            $4 \ULine{1 \ULine{1 \LeafSimple}} \LeafSimple \ULine{1 \LeafSimple} \ULine{1
            \LeafSimple}$
        };
        \node[GraphLabeledVertex](7)at(2.00,1.50,4.00){
            $4 \ULine{1 \ULine{1 \LeafSimple}} \LeafSimple \ULine{1 \ULine{1 \LeafSimple}}
            \LeafSimple$
        };
        \node[GraphLabeledVertex](8)at(2.00,1.75,1.00){
            $4 \ULine{1 \ULine{1 \LeafSimple}} \ULine{1 \LeafSimple} \LeafSimple \ULine{1
            \LeafSimple}$
        };
        \node[GraphLabeledVertex](9)at(2.00,1.75,1.50){
            $4 \ULine{1 \ULine{1 \LeafSimple}} \ULine{1 \LeafSimple} \ULine{1 \LeafSimple}
            \LeafSimple$
        };
        \node[GraphLabeledVertex](10)at(2.00,1.75,4.00){
            $4 \ULine{1 \ULine{1 \LeafSimple}} \ULine{1 \ULine{1 \LeafSimple}} \LeafSimple
            \LeafSimple$
        };
        \node[GraphLabeledVertex](11)at(2.00,3.00,1.00){
            $4 \ULine{1 \ULine{1 \ULine{1 \LeafSimple}}} \LeafSimple \LeafSimple \ULine{1
            \LeafSimple}$
        };
        \node[GraphLabeledVertex](12)at(2.00,3.00,1.50){
            $4 \ULine{1 \ULine{1 \ULine{1 \LeafSimple}}} \LeafSimple \ULine{1 \LeafSimple}
            \LeafSimple$
        };
        \node[GraphLabeledVertex](13)at(2.00,3.00,1.75){
            $4 \ULine{1 \ULine{1 \ULine{1 \LeafSimple}}} \ULine{1 \LeafSimple} \LeafSimple
            \LeafSimple$
        };
        \node[GraphLabeledVertex](14)at(2.00,3.00,4.00){
            $4 \ULine{1 \ULine{1 \ULine{1 \ULine{1 \LeafSimple}}}} \LeafSimple \LeafSimple
            \LeafSimple$
        };
        \draw[GraphArc](1)--(2);
        \draw[GraphArc](1)--(3);
        \draw[GraphArc](1)--(6);
        \draw[GraphArc](2)--(5);
        \draw[GraphArc](2)--(7);
        \draw[GraphArc](3)--(4);
        \draw[GraphArc](3)--(11);
        \draw[GraphArc](4)--(5);
        \draw[GraphArc](4)--(12);
        \draw[GraphArc](5)--(14);
        \draw[GraphArc](6)--(7);
        \draw[GraphArc](6)--(8);
        \draw[GraphArc](7)--(10);
        \draw[GraphArc](8)--(9);
        \draw[GraphArc](8)--(11);
        \draw[GraphArc](9)--(10);
        \draw[GraphArc](9)--(12);
        \draw[GraphArc](10)--(14);
        \draw[GraphArc](11)--(12);
        \draw[GraphArc](12)--(13);
        \draw[GraphArc](13)--(14);
    \end{tikzpicture}
}
    \hspace{0.25em}
    \def\Zoom{0.5}
\def\Scale{2.75}
\scalebox{\Zoom}{
    \begin{tikzpicture}[Centering,scale=\Scale,x={(-1cm,-1cm)},y={(1cm,-1cm)}]
        \DrawSquareGrid{1.50}{2.50}{1.00}{3.50}{0.50}
        \node[GraphLabeledVertex](1)at(1.50,1.00){
            $3 \ULine{2 \LeafSimple \LeafSimple} \ULine{2 \LeafSimple \LeafSimple} \ULine{2
            \LeafSimple \LeafSimple}$
        };
        \node[GraphLabeledVertex](2)at(1.50,3.00){
            $3 \ULine{2 \LeafSimple \LeafSimple} \ULine{2 \LeafSimple \ULine{2 \LeafSimple
            \LeafSimple}} \LeafSimple$
        };
        \node[GraphLabeledVertex](3)at(1.50,3.50){
            $3 \ULine{2 \LeafSimple \LeafSimple} \ULine{2 \ULine{2 \LeafSimple \LeafSimple}
            \LeafSimple} \LeafSimple$
        };
        \node[GraphLabeledVertex](4)at(2.00,1.00){
            $3 \ULine{2 \LeafSimple \ULine{2 \LeafSimple \LeafSimple}} \LeafSimple \ULine{2
            \LeafSimple \LeafSimple}$
        };
        \node[GraphLabeledVertex](5)at(2.00,1.50){
            $3 \ULine{2 \LeafSimple \ULine{2 \LeafSimple \LeafSimple}} \ULine{2 \LeafSimple
            \LeafSimple} \LeafSimple$
        };
        \node[GraphLabeledVertex](6)at(2.00,3.00){
            $3 \ULine{2 \LeafSimple \ULine{2 \LeafSimple \ULine{2 \LeafSimple \LeafSimple}}}
            \LeafSimple \LeafSimple$
        };
        \node[GraphLabeledVertex](7)at(2.00,3.50){
            $3 \ULine{2 \LeafSimple \ULine{2 \ULine{2 \LeafSimple \LeafSimple} \LeafSimple}}
            \LeafSimple \LeafSimple$
        };
        \node[GraphLabeledVertex](8)at(2.50,1.00){
            $3 \ULine{2 \ULine{2 \LeafSimple \LeafSimple} \LeafSimple} \LeafSimple \ULine{2
            \LeafSimple \LeafSimple}$
        };
        \node[GraphLabeledVertex](9)at(2.50,1.50){
            $3 \ULine{2 \ULine{2 \LeafSimple \LeafSimple} \LeafSimple} \ULine{2 \LeafSimple
            \LeafSimple} \LeafSimple$
        };
        \node[GraphLabeledVertex](10)at(2.50,2.00){
            $3 \ULine{2 \ULine{2 \LeafSimple \LeafSimple} \ULine{2 \LeafSimple \LeafSimple}}
            \LeafSimple \LeafSimple$
        };
        \node[GraphLabeledVertex](11)at(2.50,3.00){
            $3 \ULine{2 \ULine{2 \LeafSimple \ULine{2 \LeafSimple \LeafSimple}} \LeafSimple}
            \LeafSimple \LeafSimple$
        };
        \node[GraphLabeledVertex](12)at(2.50,3.50){
            $3 \ULine{2 \ULine{2 \ULine{2 \LeafSimple \LeafSimple} \LeafSimple} \LeafSimple}
            \LeafSimple \LeafSimple$
        };
        \draw[GraphArc](1)--(2);
        \draw[GraphArc](1)--(4);
        \draw[GraphArc](2)--(3);
        \draw[GraphArc](2)--(6);
        \draw[GraphArc](3)--(7);
        \draw[GraphArc](4)--(5);
        \draw[GraphArc](4)--(8);
        \draw[GraphArc](5)--(6);
        \draw[GraphArc](5)--(9);
        \draw[GraphArc](6)--(7);
        \draw[GraphArc](6)--(11);
        \draw[GraphArc](7)--(12);
        \draw[GraphArc](8)--(9);
        \draw[GraphArc](9)--(10);
        \draw[GraphArc](10)--(11);
        \draw[GraphArc](11)--(12);
    \end{tikzpicture}
}
    \caption[]{
        The Hasse diagram of the $1$-Fuss-Catalan easterly wind poset of order $4$ on the
        left, and the Hasse diagram of the $2$-Fuss-Catalan easterly wind poset of order $3$
        on the right. Observe that these graphs are not regular (that is, not all vertices
        have the same degree).
    }
    \label{fig:example_catalan_poset}
\end{figure}
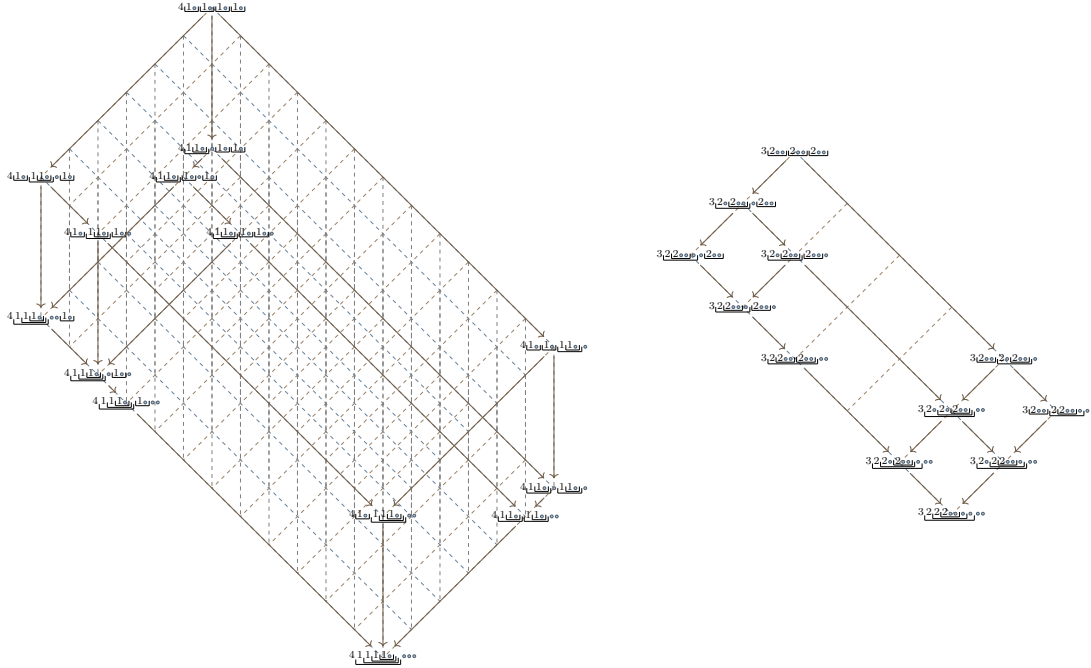

For any $n, m \in \N$, the \Def{$m$-Fuss-Catalan easterly wind poset of order $n$} is the
balanced forest $\N$-easterly wind poset of the word $m^n$, as defined in
Section~\ref{subsubsec:intervals_forests}. Figure~\ref{fig:example_catalan_poset} shows the
Hasse diagrams of some such posets.

\begin{Statement}{Lemma}{lem:fuss_catalan_posets_objects}
    For any $m, n \in \N$, $\TerminalInterval \App \ULine{\ForestMin \App m^n}$ is the
    set of $\Signature$-forests of the form $\ForestF = n \, \TermT_1 \dots \TermT_n$ such
    that $\ForestF$ is balanced and for any $\ell \in [n]$, $\Deg \App \TermT_{n - \ell + 1}
    + \dots + \Deg \App \TermT_n \leq \ell$.
\end{Statement}
\begin{Proof}
    Assume first that $\ForestF = n \, \TermT_1 \dots \TermT_n$ and $\ForestF' = n \,
    \TermT'_1 \dots \TermT'_n$ are two $\N$-forests such that $\ForestF$ satisfies the
    condition of the statement and $\ForestF \EasterlyWindCovering \ForestF'$. From the
    definition of $\EasterlyWindCovering$, we have either that $\Deg \App \TermT_j = \Deg
    \App \TermT'_j$ for all $j \in [n]$, or that there exists $j \in [n - 1]$ such that
    $\Deg \App \TermT_j = \Deg \App \TermT'_{j + 1} = 0$ and $\Deg \App \TermT_{j + 1} =
    \Deg \App \TermT'_j \ne 0$. In both cases, $\ForestF'$ satisfies also the condition of
    the statement. Since $\ForestMin \App m^n$ satisfies the condition of the statement and,
    by Theorem~\ref{thm:easterly_wind_posets}, $\EasterlyWindCovering$ is the covering
    relation of $\EasterlyWindLeq$, this shows that for any $\N$-forest $\ForestF$ such that
    $\ForestMin \App m^n \EasterlyWindLeq \ForestF$, $\ForestF$ satisfies the condition of
    the statement.

    Conversely, assume that $\ForestF = n \, \TermT_1 \dots \TermT_n$ is an $\N$-forest
    satisfying the condition of the statement and let $i$ be an internal node of $\ForestF$.
    First, if $\Parent \App \ForestF \App i \geq 2$, since $\Parent \App
    \ULine{\ForestMin \App m^n} \App i = 1$, the parent edge of $i$ in $\ForestF$
    dominates the parent edge of $i$ in $\ForestMin \App m^n$. Otherwise, when $\Parent
    \App \ForestF \App i = 1$, let us set $j := \LocalPosition \App \ForestF \App i$. Since
    $\ForestF$ satisfies the condition of the statement,
    \begin{math}
        \Deg \App \TermT_j + \dots + \Deg \App \TermT_n \leq n - j + 1.
    \end{math}
    Moreover, as the subterm of $\ForestF$ rooted at $i$ is $\TermT_j$ and the set of
    internal nodes of $\ForestF$ which are greater than of equal to $i$ is $\Bra{i, \dots, n
    + 1}$, we have
    \begin{math}
        \Deg \App \TermT_j + \dots + \Deg \App \TermT_n = n - i + 2.
    \end{math}
    This implies that $j \leq i - 1$, so that the parent edge of $i$ in $\ForestF$ dominates
    the parent edge $\Par{1, i - 1, i}$ of $i$ in $\ForestMin \App m^n$. Hence, by
    Lemma~\ref{lem:order_domination}, we have in both cases $\ForestMin \App m^n
    \EasterlyWindLeq \ForestF$, as expected.
\end{Proof}

For any $m \in \N$, an \Def{$m$-tree} is an element of $\SetTerms \App \Bra{m + 1}$.
Besides, an \Def{$m$-binary tree} is an element of $\SetTerms \App \ULine{\SetTerms \App
\Bra{m}}$ where, here, $\SetTerms \App \Bra{m}$ is seen as a signature whose all elements
are of arity $2$. In other words, an $m$-binary tree is a binary tree whose internal nodes
are decorated by $m - 1$-trees. Given an $m - 1$-tree $\TermS$ of degree $k \geq 1$, the
\Def{right comb of $\TermS$} is the $m$-binary tree such that its root is decorated by
$\TermS$, the first child of the root is a leaf, and the second child is a right comb binary
tree consisting in $k - 1$ internal nodes decorated by $\LeafSimple$. For instance, the
right comb of the $2$-tree $\TermS := 3 \, \ULine{\Corolla \App 3} \LeafSimple
\ULine{\Corolla \App 3}$ is the $3$-binary tree
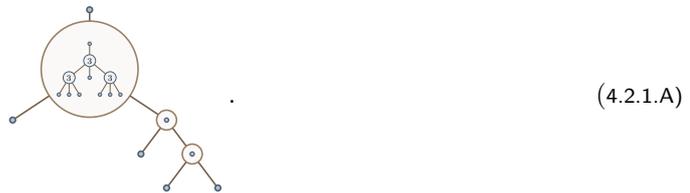
\begin{equation}
    \def\Zoom{0.6}
    \scalebox{.75}{
        \begin{tikzpicture}[Centering,xscale=0.9,yscale=0.6]
            \node[Leaf](1)at(2,1.25){};
            \node[Node,MarkB,fill=ColB!5,inner sep=1pt](2)at(2,-0.5){
                \scalebox{\Zoom}{
                    \begin{tikzpicture}[Centering,xscale=0.3,yscale=0.5]
                        \node[Leaf](1A)at(3,0){};
                        \node[Node](2A)at(3,-1){$3$};
                        \node[Node](3A)at(1,-2){$3$};
                        \node[Leaf](4A)at(0,-3){};
                        \node[Leaf](5A)at(1,-3){};
                        \node[Leaf](6A)at(2,-3){};
                        \node[Leaf](7A)at(3,-2){};
                        \node[Node](8A)at(5,-2){$3$};
                        \node[Leaf](9A)at(4,-3){};
                        \node[Leaf](10A)at(5,-3){};
                        \node[Leaf](11A)at(6,-3){};
                        \draw[Edge](1A)--(2A);
                        \draw[Edge](2A)--(3A);
                        \draw[Edge](2A)--(7A);
                        \draw[Edge](2A)--(8A);
                        \draw[Edge](3A)--(4A);
                        \draw[Edge](3A)--(5A);
                        \draw[Edge](3A)--(6A);
                        \draw[Edge](8A)--(9A);
                        \draw[Edge](8A)--(10A);
                        \draw[Edge](8A)--(11A);
                    \end{tikzpicture}
                }
            };
            \node[Leaf](3)at(0.5,-2){};
            \node[Node,MarkB,fill=ColB!5,inner sep=1pt](6)at(3.5,-2){$\LeafSimple$};
            \node[Leaf](7)at(3,-3){};
            \node[Node,MarkB,fill=ColB!5,inner sep=1pt](8)at(4,-3){$\LeafSimple$};
            \node[Leaf](9)at(3.5,-4){};
            \node[Leaf](16)at(4.5,-4){};
            \draw[Edge](1)--(2);
            \draw[Edge](2)--(3);
            \draw[Edge](2)--(6);
            \draw[Edge](6)--(7);
            \draw[Edge](6)--(8);
            \draw[Edge](8)--(9);
            \draw[Edge](8)--(16);
        \end{tikzpicture}
    }.
\end{equation}

Let the map $\ToBinaryTree$ from the set of $m$-trees to the set of $m$-binary trees defined
recursively, for any $m$-tree $\TermT$, as follows. First, if $\TermT = \LeafSimple$, then
$\ToBinaryTree \App \TermT = \LeafSimple$. Otherwise, let $\TermS$ be the $m - 1$-tree
obtained by keeping the root of $\TermT$ and by deleting recursively all first subterms of
the kept internal nodes. In this process, let us denote by $\TermT_1$, \dots, $\TermT_k$ the
forgotten subterms from left to right, with $k := \Deg \App \TermS$. Now, let $\TermR$ be
the right comb of $\TermS$. With these definitions, $\ToBinaryTree \App \TermT$ is obtained
by replacing, for any $i \in [k]$, the $i$-th leaf of $\TermR$ by $\ToBinaryTree \App
\TermT_i$. The last leaf of $\TermR$ is left as is. For instance, for $m := 3$, we have
\begin{equation} \label{equ:example_to_binary_tree}
    \ToBinaryTree \App
    \scalebox{.75}{
        \begin{tikzpicture}[Centering,xscale=0.3,yscale=0.7]
            \node[Leaf](1)at(8,0){};
            \node[Node](2)at(8,-1){$4$};
            \node[Node](3)at(1.5,-2){$4$};
            \node[Leaf](4)at(0,-3){};
            \node[Leaf](5)at(1,-3){};
            \node[Leaf](6)at(2,-3){};
            \node[Leaf](7)at(3,-3){};
            \node[Node](8)at(5.5,-2){$4$};
            \node[Leaf](9)at(4,-3){};
            \node[Leaf](10)at(5,-3){};
            \node[Leaf](11)at(6,-3){};
            \node[Leaf](12)at(7,-3){};
            \node[Node](13)at(14.5,-2){$4$};
            \node[Node](14)at(12.5,-3){$4$};
            \node[Node](15)at(10.5,-4){$4$};
            \node[Leaf](16)at(9,-5){};
            \node[Leaf](17)at(10,-5){};
            \node[Leaf](18)at(11,-5){};
            \node[Leaf](19)at(12,-5){};
            \node[Leaf](20)at(12,-4){};
            \node[Leaf](21)at(13,-4){};
            \node[Node](22)at(14.5,-4){$4$};
            \node[Leaf](23)at(13,-5){};
            \node[Leaf](24)at(14,-5){};
            \node[Leaf](25)at(15,-5){};
            \node[Leaf](26)at(16,-5){};
            \node[Leaf](27)at(14,-3){};
            \node[Leaf](28)at(15,-3){};
            \node[Leaf](29)at(16,-3){};
            \node[Leaf](30)at(10.5,-2){};
            \draw[Edge](1)--(2);
            \draw[Edge](2)--(3);
            \draw[Edge](2)--(8);
            \draw[Edge](2)--(13);
            \draw[Edge](3)--(4);
            \draw[Edge](3)--(5);
            \draw[Edge](3)--(6);
            \draw[Edge](3)--(7);
            \draw[Edge](8)--(9);
            \draw[Edge](8)--(10);
            \draw[Edge](8)--(11);
            \draw[Edge](8)--(12);
            \draw[Edge](13)--(14);
            \draw[Edge](13)--(27);
            \draw[Edge](13)--(28);
            \draw[Edge](13)--(29);
            \draw[Edge](14)--(15);
            \draw[Edge](14)--(20);
            \draw[Edge](14)--(21);
            \draw[Edge](14)--(22);
            \draw[Edge](15)--(16);
            \draw[Edge](15)--(17);
            \draw[Edge](15)--(18);
            \draw[Edge](15)--(19);
            \draw[Edge](22)--(23);
            \draw[Edge](22)--(24);
            \draw[Edge](22)--(25);
            \draw[Edge](22)--(26);
            \draw[Edge](2)--(30);
        \end{tikzpicture}
    }
    \enspace = \enspace
    \def\Zoom{0.6}
    \scalebox{.75}{
        \begin{tikzpicture}[Centering,xscale=1.2,yscale=0.8]
            \node[Leaf](1)at(2,0.75){};
            \node[Node,MarkB,fill=ColB!5,inner sep=1pt](2)at(2,-0.5){
                \scalebox{\Zoom}{
                    \begin{tikzpicture}[Centering,xscale=0.3,yscale=0.5]
                        \node[Leaf](1A)at(3,0){};
                        \node[Node](2A)at(3,-1){$3$};
                        \node[Node](3A)at(1,-2){$3$};
                        \node[Leaf](4A)at(0,-3){};
                        \node[Leaf](5A)at(1,-3){};
                        \node[Leaf](6A)at(2,-3){};
                        \node[Leaf](7A)at(3,-2){};
                        \node[Node](8A)at(5,-2){$3$};
                        \node[Leaf](9A)at(4,-3){};
                        \node[Leaf](10A)at(5,-3){};
                        \node[Leaf](11A)at(6,-3){};
                        \draw[Edge](1A)--(2A);
                        \draw[Edge](2A)--(3A);
                        \draw[Edge](2A)--(7A);
                        \draw[Edge](2A)--(8A);
                        \draw[Edge](3A)--(4A);
                        \draw[Edge](3A)--(5A);
                        \draw[Edge](3A)--(6A);
                        \draw[Edge](8A)--(9A);
                        \draw[Edge](8A)--(10A);
                        \draw[Edge](8A)--(11A);
                    \end{tikzpicture}
                }
            };
            \node[Node,MarkB,fill=ColB!5,inner sep=1pt](3)at(0.5,-2){
                \scalebox{\Zoom}{
                    \begin{tikzpicture}[Centering,xscale=0.3,yscale=0.5]
                        \node[Leaf](1B)at(1,0){};
                        \node[Node](2B)at(1,-1){$3$};
                        \node[Leaf](3B)at(0,-2){};
                        \node[Leaf](4B)at(1,-2){};
                        \node[Leaf](5B)at(2,-2){};
                        \draw[Edge](1B)--(2B);
                        \draw[Edge](2B)--(3B);
                        \draw[Edge](2B)--(4B);
                        \draw[Edge](2B)--(5B);
                    \end{tikzpicture}
                }
            };
            \node[Leaf](4)at(0,-3){};
            \node[Leaf](5)at(1,-3){};
            \node[Node,MarkB,fill=ColB!5,inner sep=1pt](6)at(3.5,-2){$\LeafSimple$};
            \node[Leaf](7)at(3,-3){};
            \node[Node,MarkB,fill=ColB!5,inner sep=1pt](8)at(4,-3){$\LeafSimple$};
            \node[Node,MarkB,fill=ColB!5,inner sep=1pt](9)at(3.5,-4){
                \scalebox{\Zoom}{
                    \begin{tikzpicture}[Centering,xscale=0.4,yscale=0.5]
                        \node[Leaf](1C)at(2,0){};
                        \node[Node](2C)at(2,-1){$3$};
                        \node[Leaf](3C)at(1,-2){};
                        \node[Leaf](4C)at(2,-2){};
                        \node[Node](5C)at(3,-2){$3$};
                        \node[Leaf](6C)at(2,-3){};
                        \node[Leaf](7C)at(3,-3){};
                        \node[Leaf](8C)at(4,-3){};
                        \draw[Edge](1C)--(2C);
                        \draw[Edge](2C)--(3C);
                        \draw[Edge](2C)--(4C);
                        \draw[Edge](2C)--(5C);
                        \draw[Edge](5C)--(6C);
                        \draw[Edge](5C)--(7C);
                        \draw[Edge](5C)--(8C);
                    \end{tikzpicture}
                }
            };
            \node[Node,MarkB,fill=ColB!5,inner sep=1pt](10)at(2.5,-5){
                \scalebox{\Zoom}{
                    \begin{tikzpicture}[Centering,xscale=0.3,yscale=0.5]
                        \node[Leaf](1D)at(1,0){};
                        \node[Node](2D)at(1,-1){$3$};
                        \node[Leaf](3D)at(0,-2){};
                        \node[Leaf](4D)at(1,-2){};
                        \node[Leaf](5D)at(2,-2){};
                        \draw[Edge](1D)--(2D);
                        \draw[Edge](2D)--(3D);
                        \draw[Edge](2D)--(4D);
                        \draw[Edge](2D)--(5D);
                    \end{tikzpicture}
                }
            };
            \node[Leaf](11)at(2,-6){};
            \node[Leaf](12)at(3,-6){};
            \node[Node,MarkB,fill=ColB!5,inner sep=1pt](13)at(4.5,-5){$\LeafSimple$};
            \node[Leaf](14)at(4,-6){};
            \node[Leaf](15)at(5,-6){};
            \node[Leaf](16)at(4.5,-4){};
            \draw[Edge](1)--(2);
            \draw[Edge](2)--(3);
            \draw[Edge](2)--(6);
            \draw[Edge](3)--(4);
            \draw[Edge](3)--(5);
            \draw[Edge](6)--(7);
            \draw[Edge](6)--(8);
            \draw[Edge](8)--(9);
            \draw[Edge](8)--(16);
            \draw[Edge](9)--(10);
            \draw[Edge](9)--(13);
            \draw[Edge](10)--(11);
            \draw[Edge](10)--(12);
            \draw[Edge](13)--(14);
            \draw[Edge](13)--(15);
        \end{tikzpicture}
    }.
\end{equation}

Finally, let $\BinaryTreesImage \App m$ be the set of $m$-binary trees $\TermR$ such that
from any internal node of $\TermR$ decorated by an $m - 1$-tree $\TermS$ of degree $k \geq
1$, there is a right branch consisting in $k - 1$ internal nodes decorated by $\LeafSimple$,
and each internal node decorated by $\LeafSimple$ in $\TermR$ is a part of such right
branch.

\begin{Statement}{Lemma}{lem:to_binary_tree_bijection}
    For any $m \in \N$, from the domain consisting in the set of $m$-trees and on the
    codomain $\BinaryTreesImage \App m$, the map $\ToBinaryTree$ is a bijection.
\end{Statement}
\begin{Proof}
    Let $\phi : \BinaryTreesImage \App m \to \SetTerms \App \Bra{m + 1}$ be the map defined
    recursively, for any $\TermR \in \BinaryTreesImage \App m$ of degree $n$, as follows.
    First, if $n = 0$, then $\TermR = \LeafSimple$. In this case, set $\phi \App \TermR :=
    \LeafSimple$. Otherwise, we have $n \geq 1$ and, from the description of
    $\BinaryTreesImage \App m$, $\TermR$ is the right comb $\TermR'$ of an $m - 1$-tree
    $\TermS$ of degree $k \geq 1$ such that for any $i \in [k]$, the $i$-th leaf of
    $\TermR'$ is attached to a subterm $\TermR_i$ of $\TermR$. Since for any $i \in [k]$,
    $\TermR_i$ belongs to $\BinaryTreesImage \App m$, the $m$-tree $\TermT_i := \phi \App
    \TermR_i$ is, by induction, well-defined. Let also $\TermT'$ be the $m$-tree obtained by
    adding to each internal node of $\TermS$ a leaf as first child. We define $\phi \App
    \TermR$ as the $m$-tree obtained by replacing, for any $i \in [k]$, the first leaf of
    the internal node $i$ of $\TermT'$ by $\TermT_i$. It follows by induction on $n$ that
    $\phi$ is a well-defined map. Again by induction on $n$, it is straightforward to show
    that $\phi$ is the inverse map of $\ToBinaryTree$.
\end{Proof}

The \Def{inorder traversal} of an $m$-binary tree $\TermR$ is defined recursively as
follows. It $\TermR = \LeafSimple$, then the inorder traversal of $\TermR$ is empty.
Otherwise, we have $\TermR = \TermS \, \TermR_1 \TermR_2$ where $\TermS$ is an $m - 1$-tree,
and $\TermR_1$ and $\TermR_2$ are two $m$-binary trees. In this case, $\TermR_1$ is visited
according to the inorder traversal, then the root of $\TermR$, and finally, $\TermR_2$ is
visited according to the inorder traversal. This procedure induces a total order on the
internal nodes of $\TermR$ where the first visited internal node is the smallest one. Now,
let the map $\ToForest$ from the set of $m$-binary trees to the set of balanced $\N$-forests
such that, for any $m$-binary tree $\TermR$ of degree $n$, $\ToForest \App \TermR$ is the
$\N$-forest $n \, \TermS_1 \dots \TermS_n$ where for any $i \in [n]$, $\TermS_i$ is the
decoration of the $i$-th visited internal node of $\TermR$ w.r.t.\ the inorder traversal of
$\TermR$. For instance, by considering the $3$-binary tree $\TermR$ of the right-hand side
of~\eqref{equ:example_to_binary_tree},
\begin{equation}
    \ToForest \App \TermR
    =
    7 \,
    \ULine{\Corolla \App 3}
    \ULine{
        3 \,
        \ULine{\Corolla \App 3}
        \LeafSimple
        \ULine{\Corolla \App 3}
    }
    \LeafSimple
    \ULine{\Corolla \App 3}
    \ULine{3 \, \LeafSimple \LeafSimple \ULine{\Corolla \App 3}}
    \LeafSimple
    \LeafSimple.
\end{equation}

\begin{Statement}{Lemma}{lem:to_forest_bijection}
    For any $m \in \N$, from the domain $\BinaryTreesImage \App m$ and on the codomain
    $\bigcup_{n \in \N} \TerminalInterval \App \Par{\ForestMin \App m^n}$, the map
    $\ToForest $ is a bijection.
\end{Statement}
\begin{Proof}
    Let $\phi : \bigcup_{n \in \N} \TerminalInterval \App \Par{\ForestMin \App m^n} \to
    \BinaryTreesImage \App m$ be the map defined recursively, for any $\ForestF \in
    \TerminalInterval \App \Par{\ForestMin \App m^n}$, $n \in \N$, as follows. First,
    if $n = 0$, then $\ForestF = 0$. In this case, set $\phi \App \ForestF := \LeafSimple$.
    Otherwise, we have $n \geq 1$ and by Lemma~\ref{lem:fuss_catalan_posets_objects}, it
    follows by induction on $n$ that $\ForestF$ decomposes as
    \begin{equation}
        \ForestF =
        \ForestF_1 \; \Conc \;
        1 \, \TermS
        \; \Conc \; k - 1 \, \underbrace{\LeafSimple \dots \LeafSimple}_{k - 1}
        \; \Conc \; \ForestF_2
    \end{equation}
    where $\ForestF_1$ and $\ForestF_2$ are two $\N$-forests and $\TermS$ is an $m - 1$-tree
    $\TermS$ of degree $k \geq 1$. Let us consider this decomposition when the size $\ell
    \geq 0$ of $\ForestF_1$ is minimal. Since $\ForestF_1 \Conc \ForestF_2$ satisfies the
    conditions described in Lemma~\ref{lem:fuss_catalan_posets_objects}, the $m$-binary tree
    $\TermR := \phi \App \ULine{\ForestF_1 \Conc \ForestF_2}$ is, by induction,
    well-defined. Let also $\TermR'$ be the right comb of $\TermS$. We define $\phi \App
    \ForestF$ as the $m$-binary tree obtained by inserting the root of $\TermR'$ onto the
    unique edge of $\TermR$ such that, w.r.t.\ the inorder traversal, the internal nodes
    coming from $\ForestF_1$ are visited first, then the ones of $\TermR'$ are visited, and
    finally the ones coming from $\ForestF_2$ are visited. It follows by induction on $n$
    that $\phi$ is a well-defined map. Again by induction on $n$, it is straightforward to
    show that $\phi$ is the inverse map of~$\ToForest$.
\end{Proof}

\begin{Statement}{Theorem}{thm:fuss_catalan_posets}
    For any $m, n \in \N$, the underlying set of the $m$-Fuss-Catalan easterly wind poset of
    order $n$ is in one-to-one correspondence with the set of $m$-trees of degree $n$. The
    map $\ToForest \circ \ToBinaryTree$ is such a one-to-one correspondence.
\end{Statement}
\begin{Proof}
    By Lemma~\ref{lem:to_binary_tree_bijection}, $\ToBinaryTree$ is a bijection between the
    set of $m$-trees and $\BinaryTreesImage \App m$. Moreover, By
    Lemma~\ref{lem:to_forest_bijection}, $\ToForest$ is a bijection between
    $\BinaryTreesImage \App m$ and $\bigcup_{n \in \N} \TerminalInterval \App
    \ULine{\ForestMin \App m^n}$. Since these two bijections preserve the degree, the
    composition $\ToForest \circ \ToBinaryTree$ satisfies the property described in the
    statement of the theorem.
\end{Proof}

By Theorem~\ref{thm:fuss_catalan_posets}, the $m$-Fuss-Catalan easterly wind posets involve
the combinatorial family of Fuss-Catalan objects. Hence, the cardinality of such posets of
order $n$ is
\begin{equation}
    \frac{1}{mn + 1} \binom{mn + n}{n}.
\end{equation}
Many other posets involving this family of objects exist~\cite{BPR12,CG22}, and our posets
differ from those presented in the cited works.

\subsubsection{Rooted tree lattices} \label{subsubsec:rooted_tree_lattices}
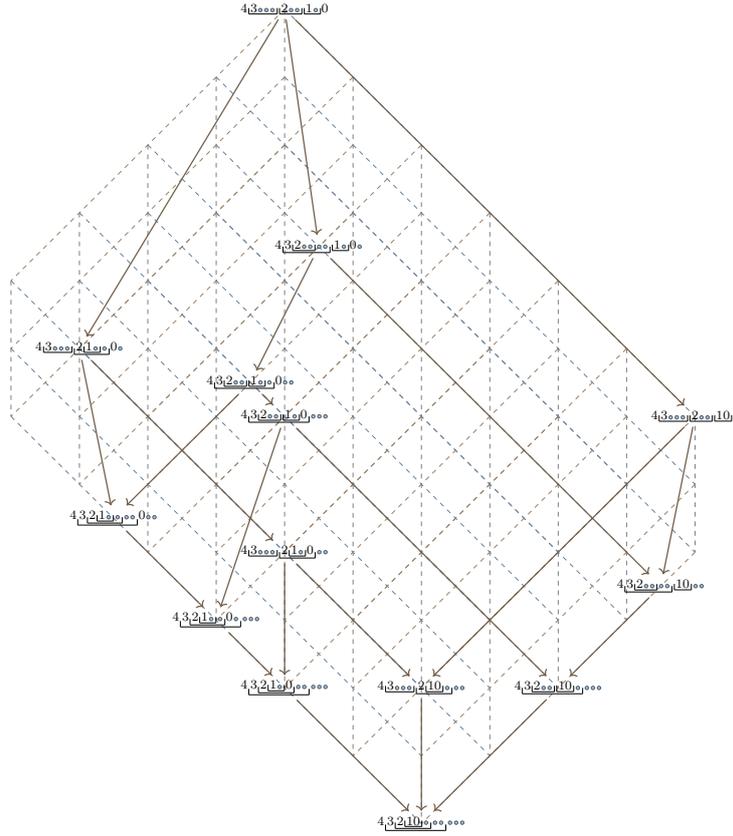
\begin{figure}[ht]
    \centering
    \def\Zoom{0.6}
\def\Scale{3.0}
\scalebox{\Zoom}{
    \begin{tikzpicture}[Centering,scale=\Scale,
    x={(0,-1cm)},y={(-1cm,-1cm)}, z={(1cm,-1cm)}]
        \DrawCubeGrid{1.75}{2.75}{1.50}{3.50}{1.00}{4.00}{0.50}
        \node[GraphLabeledVertex](1)at(1.75,1.50,1.00){
            $4 \ULine{3 \LeafSimple \LeafSimple \LeafSimple} \ULine{2 \LeafSimple
            \LeafSimple} \ULine{1 \LeafSimple} 0$
        };
        \node[GraphLabeledVertex](2)at(1.75,1.50,4.00){
            $4 \ULine{3 \LeafSimple \LeafSimple \LeafSimple} \ULine{2 \LeafSimple
            \LeafSimple} \ULine{1 0} \LeafSimple$
        };
        \node[GraphLabeledVertex](3)at(1.75,3.50,1.50){
            $4 \ULine{3 \LeafSimple \LeafSimple \LeafSimple} \ULine{2 \ULine{1 \LeafSimple}
            \LeafSimple} 0 \LeafSimple$
        };
        \node[GraphLabeledVertex](4)at(1.75,3.50,3.00){
            $4 \ULine{3 \LeafSimple \LeafSimple \LeafSimple} \ULine{2 \ULine{1 \LeafSimple}
            0} \LeafSimple \LeafSimple$
        };
        \node[GraphLabeledVertex](5)at(1.75,3.50,4.00){
            $4 \ULine{3 \LeafSimple \LeafSimple \LeafSimple} \ULine{2 \ULine{1 0}
            \LeafSimple} \LeafSimple \LeafSimple$
        };
        \node[GraphLabeledVertex](6)at(2.75,1.75,1.50){
            $4 \ULine{3 \ULine{2 \LeafSimple \LeafSimple} \LeafSimple \LeafSimple} \ULine{1
            \LeafSimple} 0 \LeafSimple$
        };
        \node[GraphLabeledVertex](7)at(2.75,1.75,4.00){
            $4 \ULine{3 \ULine{2 \LeafSimple \LeafSimple} \LeafSimple \LeafSimple} \ULine{1
            0} \LeafSimple \LeafSimple$
        };
        \node[GraphLabeledVertex](8)at(2.75,2.50,1.75){
            $4 \ULine{3 \ULine{2 \LeafSimple \LeafSimple} \ULine{1 \LeafSimple} \LeafSimple}
            0 \LeafSimple \LeafSimple$
        };
        \node[GraphLabeledVertex](9)at(2.75,2.50,2.00){
            $4 \ULine{3 \ULine{2 \LeafSimple \LeafSimple} \ULine{1 \LeafSimple} 0}
            \LeafSimple \LeafSimple \LeafSimple$
        };
        \node[GraphLabeledVertex](10)at(2.75,2.50,4.00){
            $4 \ULine{3 \ULine{2 \LeafSimple \LeafSimple} \ULine{1 0} \LeafSimple}
            \LeafSimple \LeafSimple \LeafSimple$
        };
        \node[GraphLabeledVertex](11)at(2.75,3.50,1.75){
            $4 \ULine{3 \ULine{2 \ULine{1 \LeafSimple} \LeafSimple} \LeafSimple \LeafSimple}
            0 \LeafSimple \LeafSimple$
        };
        \node[GraphLabeledVertex](12)at(2.75,3.50,2.50){
            $4 \ULine{3 \ULine{2 \ULine{1 \LeafSimple} \LeafSimple} 0 \LeafSimple}
            \LeafSimple \LeafSimple \LeafSimple$
        };
        \node[GraphLabeledVertex](13)at(2.75,3.50,3.00){
            $4 \ULine{3 \ULine{2 \ULine{1 \LeafSimple} 0} \LeafSimple \LeafSimple}
            \LeafSimple \LeafSimple \LeafSimple$
        };
        \node[GraphLabeledVertex](14)at(2.75,3.50,4.00){
            $4 \ULine{3 \ULine{2 \ULine{1 0} \LeafSimple} \LeafSimple \LeafSimple}
            \LeafSimple \LeafSimple \LeafSimple$
        };
        \draw[GraphArc](1)--(2);
        \draw[GraphArc](1)--(3);
        \draw[GraphArc](1)--(6);
        \draw[GraphArc](2)--(5);
        \draw[GraphArc](2)--(7);
        \draw[GraphArc](3)--(4);
        \draw[GraphArc](3)--(11);
        \draw[GraphArc](4)--(5);
        \draw[GraphArc](4)--(13);
        \draw[GraphArc](5)--(14);
        \draw[GraphArc](6)--(7);
        \draw[GraphArc](6)--(8);
        \draw[GraphArc](7)--(10);
        \draw[GraphArc](8)--(9);
        \draw[GraphArc](8)--(11);
        \draw[GraphArc](9)--(10);
        \draw[GraphArc](9)--(12);
        \draw[GraphArc](10)--(14);
        \draw[GraphArc](11)--(12);
        \draw[GraphArc](12)--(13);
        \draw[GraphArc](13)--(14);
    \end{tikzpicture}
}
    \caption[]{
        The Hasse diagram of the rooted tree easterly wind poset of order $4$.
    }
    \label{fig:example_rooted_tree_lattice}
\end{figure}

For any $n \in \N$, let $\DecreasingWord \App n$ be the word on $\N$ of length $n$ such that
for any $i \in [n]$, $\DecreasingWord \App n \App i = n - i$. The \Def{rooted tree easterly
wind poset of order $n$} is the fully tilted $\N$-easterly wind poset of $\DecreasingWord
\App n$, as defined in Section~\ref{subsubsec:fully_tilted_terms}.
Figure~\ref{fig:example_rooted_tree_lattice} shows the Hasse diagrams of such poset.
By definition, all $\N$-terms of the rooted tree easterly wind posets are balanced forests.
Observe that the minimal element of the rooted tree easterly wind poset of order $n$ is
\begin{math}
    \DecreasingTree \App n
    := \ForestMin \App \ULine{\DecreasingWord \App n}
    = n \ULine{\Corolla \App \ULine{n - 1}} \dots \ULine{\Corolla \App 0}.
\end{math}
For instance,
\begin{math}
    \DecreasingTree \App 4
    = 4 \,
    \ULine{3 \LeafSimple \LeafSimple \LeafSimple}
    \ULine{2 \LeafSimple \LeafSimple}
    \ULine{1 \LeafSimple}
    0.
\end{math}

A \Def{rooted tree} is recursively defined as a node together with a possibly empty list of
rooted trees, each of which is attached to the node as a child. The \Def{size} of a rooted
tree $\TermR$ is the number of nodes of $\TermR$. The \Def{underlying rooted tree} of an
$\Signature$-term $\TermT$ is the rooted tree $\ToRootedTree \App \TermT$ obtained by
removing the leaves and their adjacent edges in $\TermT$, as well as the decorations of its
internal nodes. For instance,
\begin{equation}
    \ToRootedTree
    \App
    \ULine{
        \GenA_4
        \LeafSimple
        \ULine{\GenA_2 \ULine{\GenA_1 \LeafSimple} \GenA_0}
        \LeafSimple
        \ULine{\GenA_3 \ULine{\GenA_2 \LeafSimple \LeafSimple} \LeafSimple \LeafSimple}
    }
    =
    \scalebox{.65}{
        \begin{tikzpicture}[Centering,xscale=0.6,yscale=0.60]
            \node[Node](1)at(1.25,0){};
            \node[Node](2)at(0.5,-1){};
            \node[Node](3)at(0,-2){};
            \node[Node](4)at(1,-2){};
            \node[Node](5)at(2,-1){};
            \node[Node](6)at(2,-2){};
            \draw[Edge](1)--(2);
            \draw[Edge](1)--(5);
            \draw[Edge](2)--(3);
            \draw[Edge](2)--(4);
            \draw[Edge](5)--(6);
        \end{tikzpicture}
    }.
\end{equation}
Moreover, the \Def{fully tilted term} $\ToFullyTiltedTerm \App \TermR$ of a rooted tree
$\TermR$ is the fully tilted $\N$-term obtained by labelling from $n - 1$ to $0$ the nodes
of $\TermR$ w.r.t.\ the preorder traversal and then, by grafting to each node some leaves as
rightmost children such that each node labeled by $k$ has $k$ children. This ensures that
$\ToFullyTiltedTerm \App \TermR$ is fully tilted. For instance,
\begin{equation}
    \ToFullyTiltedTerm
    \App
    \scalebox{.65}{
        \begin{tikzpicture}[Centering,xscale=0.6,yscale=0.60]
            \node[Node](1)at(1.25,0){};
            \node[Node](2)at(0.5,-1){};
            \node[Node](3)at(0,-2){};
            \node[Node](4)at(1,-2){};
            \node[Node](5)at(2,-1){};
            \node[Node](6)at(2,-2){};
            \draw[Edge](1)--(2);
            \draw[Edge](1)--(5);
            \draw[Edge](2)--(3);
            \draw[Edge](2)--(4);
            \draw[Edge](5)--(6);
        \end{tikzpicture}
    }
    =
    5 \,
    \ULine{
        4 \,
        \ULine{3 \, \LeafSimple \LeafSimple \LeafSimple}
        \ULine{2 \, \LeafSimple \LeafSimple}
        \LeafSimple
        \LeafSimple
    }
    \ULine{1 \, 0}
    \LeafSimple
    \LeafSimple
    \LeafSimple.
\end{equation}

\begin{Statement}{Lemma}{lem:rooted_trees_bijection}
    For any $n \in \N$, the map $\ToRootedTree$ is a  one-to-one correspondence from the
    underlying set of the rooted tree easterly wind poset of order $n$ and the set of
    rooted trees of size~$n$.
\end{Statement}
\begin{Proof}
    Let $T$ be the underlying set of the rooted tree easterly wind poset of order $n$. The
    fact that the map $\ToRootedTree$ on the domain $T$ is injective is a consequence of the
    fact that all $\N$-terms of $T$ are fully tilted. Let us show that each rooted tree
    $\TermR$ admits an antecedent in $T$ for the map $\ToRootedTree$. By setting $\TermT :=
    \ToFullyTiltedTerm \App \TermR$, from the definitions of the maps $\ToRootedTree$ and
    $\ToFullyTiltedTerm$, it follows immediately that $\ToRootedTree \App \TermT = \TermR$.
    It remains to prove that $\TermT$ belongs to $T$. By construction of the $\N$-term
    $\ForestMin \App \ULine{\DecreasingWord \App n}$, for any internal node $i$ of $\TermT$,
    the parent edge of $i$ in $\ForestMin \App \ULine{\DecreasingWord \App n}$ is dominated
    by that of $i$ in $\TermT$. Hence, by Lemma~\ref{lem:order_domination}, $\ForestMin \App
    \ULine{\DecreasingWord \App n} \EasterlyWindLeq \TermT$. This shows that $\TermT \in T$
    and implies the statement of the lemma.
\end{Proof}

The \Def{scope sequence} of a rooted tree $\TermR$ of size $n$ is the word $\Scope \App
\TermR$ of length $n$ such that for any $i \in [n]$, $\Scope \App \TermR \App i$ is the
number descendants of the $i$-th visited node in the preorder traversal of $\TermR$. For
instance,
\begin{equation}
    \Scope
    \App
    \scalebox{.65}{
        \begin{tikzpicture}[Centering,xscale=0.6,yscale=0.60]
            \node[Node](1)at(1.25,0){};
            \node[Node](2)at(0.5,-1){};
            \node[Node](3)at(0,-2){};
            \node[Node](4)at(1,-2){};
            \node[Node](5)at(2,-1){};
            \node[Node](6)at(2,-2){};
            \draw[Edge](1)--(2);
            \draw[Edge](1)--(5);
            \draw[Edge](2)--(3);
            \draw[Edge](2)--(4);
            \draw[Edge](5)--(6);
        \end{tikzpicture}
    }
    =
    5 2 0 0 1 0.
\end{equation}
Observe that the scope sequence of an $\Signature$-term $\TermT$ as defined in
Section~\ref{subsubsec:fully_tilted_terms} and the scope sequence of the rooted tree
$\ToRootedTree \App \TermT$ coincide, that is $\Scope \App \ULine{\ToRootedTree \App \TermT}
= \Scope \App \TermT$.

The \Def{Tamari partial order}~\cite{Tam62} is a partial order $\TamariLeq$ defined on the
family of Catalan objects of a given size. We consider here the following description of
this order involving rooted trees~\cite{Knu04}. Given two rooted trees $\TermR_1$ and
$\TermR_2$ of the same size $n$, $\TermR_1 \TamariLeq \TermR_2$ holds if and only if for any
$i \in [n]$, $\Scope \App \TermR_1 \App i \leq \Scope \App \TermR_2 \App i$.

\begin{Statement}{Proposition}{prop:rooted_tree_lattices}
    For any $n \in \N$, the map $\ToRootedTree$ is a poset isomorphism between the rooted
    tree easterly wind poset of order $n$ and the Tamari poset of order~$n$.
\end{Statement}
\begin{Proof}
    Let $T$ be the underlying set of the rooted tree easterly wind poset of order $n$. First
    of all, by Lemma~\ref{lem:rooted_trees_bijection}, the map $\ToRootedTree$ is a
    bijection between $T$ and the underlying set of the Tamari poset of order $n$. It
    remains to prove that $\ToRootedTree$ is an order embedding, that is, for any $\TermT_1,
    \TermT_2 \in T$, $\TermT_1 \EasterlyWindLeq \TermT_2$ if and only if $\ToRootedTree \App
    \TermT_1 \TamariLeq \ToRootedTree \App \TermT_2$. This property is a consequence of
    Proposition~\ref{prop:scope_sequence_fully_tilted} and the fact that, as noticed above,
    for any $\N$-term $\TermT$, $\Scope \App \TermT = \Scope \App \ULine{\ToRootedTree \App
    \TermT}$.
\end{Proof}

\subsection{Leaning forest lattices} \label{subsec:leaning_forests}
Leaning forest $\Signature$-easterly wind lattices will be used to construct bases of
natural Hopf algebras of free nonsymmetric operads in the next section. Here we introduce
these lattices and establish some of their properties. We also define two concatenation
operations and a shuffle operation on leaning forests. These operations will subsequently be
employed to describe the product of natural Hopf algebras of free nonsymmetric operads on
alternative bases.

\subsubsection{Leaning forests} \label{subsubsec:leaning_forests}
An $\Signature$-forest $\ForestF = n \, \TermT_1 \dots \TermT_n$ is \Def{leaning} if
$\ForestF$ is balanced and is $\Bra{1}$-tilted. For instance, $0$ and $4 \, \ULine{\GenA_2
\LeafSimple \ULine{\GenA_1 \GenA_0}} \GenA_0 \LeafSimple \LeafSimple$ are leaning
$\SignatureExample$-forests. On the contrary, $4 \, \GenA_0 \LeafSimple \ULine{\GenA_2
\ULine{\GenA_2 \LeafSimple \LeafSimple} \ULine{\GenA_1 \LeafSimple}} \LeafSimple$ is a
balanced $\SignatureExample$-forest which is not leaning since it is not $\Bra{1}$-tilted,
and $3 \, \ULine{\GenA_1 \LeafSimple} \LeafSimple \LeafSimple$ is a $\Bra{1}$-tilted
$\Signature$-forest which is not leaning since it is not balanced. Let us denote by
$\SetLeaningForests \App \Signature$ the set of leaning $\Signature$-forests. The
\Def{length} of a leaning $\Signature$-forest $\ForestF$ is the number of children subterms
of the root of $\ForestF$ which are not leaves. For instance, the length of $4 \,
\ULine{\GenA_2 \LeafSimple \ULine{\GenA_1 \GenA_0}} \GenA_0 \LeafSimple \LeafSimple$ is~$2$.

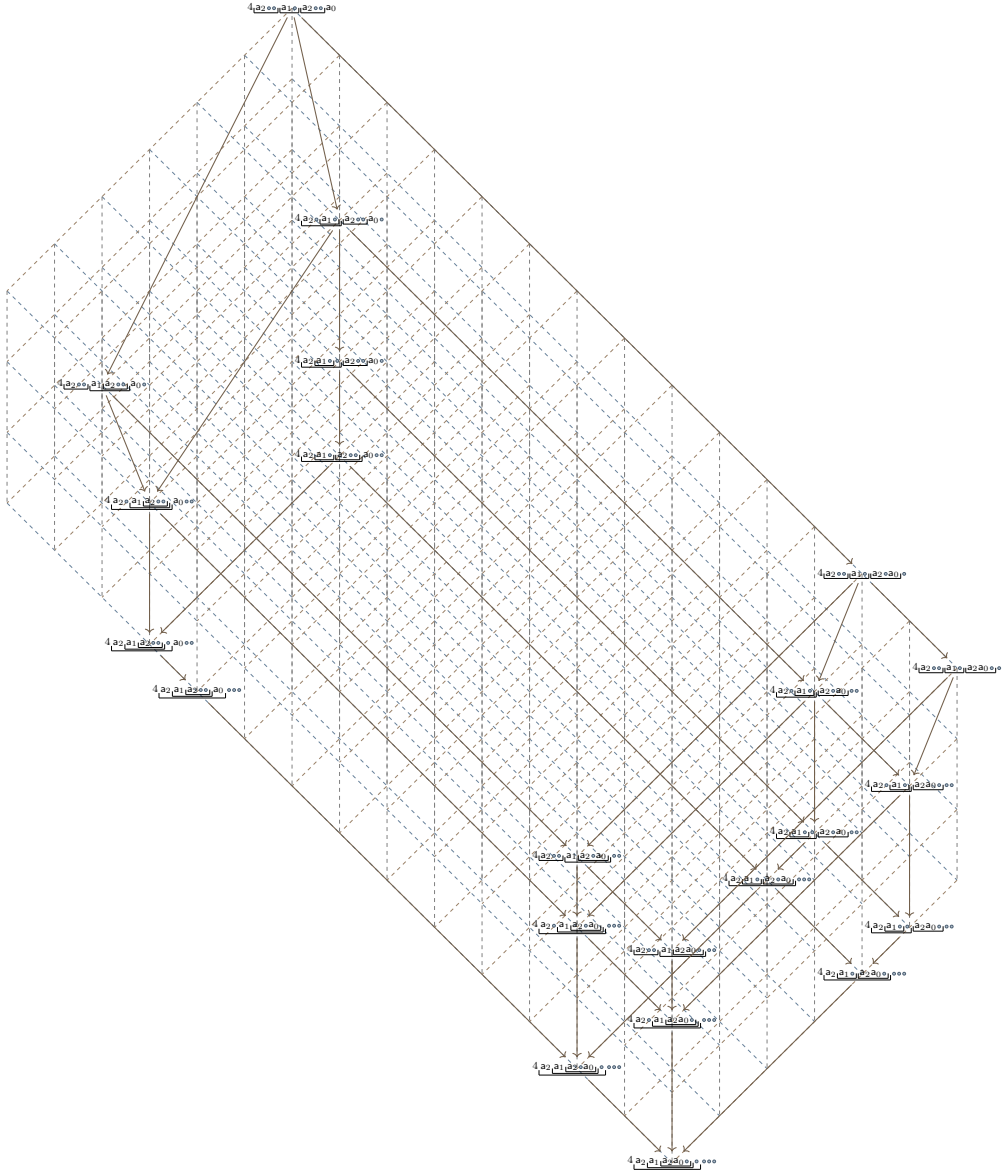
\begin{figure}[ht]
    \centering
    \def\Zoom{0.5}
\def\Scale{5.0}
\scalebox{\Zoom}{
    \begin{tikzpicture}[Centering,scale=\Scale,
    x={(0,-1.5cm)}, y={(-1cm,-1cm)}, z={(1cm,-1cm)}]
        \DrawCubeGrid{1.75}{2.50}{1.50}{3.00}{1.00}{4.50}{0.25}
        \node[GraphLabeledVertex](1)at(1.75,1.50,1.00){
            $4 \ULine{\GenA_2 \LeafSimple \LeafSimple} \ULine{\GenA_1 \LeafSimple}
            \ULine{\GenA_2 \LeafSimple \LeafSimple} \GenA_0$
        };
        \node[GraphLabeledVertex](2)at(1.75,1.50,4.00){
            $4 \ULine{\GenA_2 \LeafSimple \LeafSimple} \ULine{\GenA_1 \LeafSimple}
            \ULine{\GenA_2 \LeafSimple \GenA_0} \LeafSimple$
        };
        \node[GraphLabeledVertex](3)at(1.75,1.50,4.50){
            $4 \ULine{\GenA_2 \LeafSimple \LeafSimple} \ULine{\GenA_1 \LeafSimple}
            \ULine{\GenA_2 \GenA_0 \LeafSimple} \LeafSimple$
        };
        \node[GraphLabeledVertex](4)at(1.75,3.00,1.50){
            $4 \ULine{\GenA_2 \LeafSimple \LeafSimple} \ULine{\GenA_1 \ULine{\GenA_2
            \LeafSimple \LeafSimple}} \GenA_0 \LeafSimple$
        };
        \node[GraphLabeledVertex](5)at(1.75,3.00,4.00){
            $4 \ULine{\GenA_2 \LeafSimple \LeafSimple} \ULine{\GenA_1 \ULine{\GenA_2
            \LeafSimple \GenA_0}} \LeafSimple \LeafSimple$
        };
        \node[GraphLabeledVertex](6)at(1.75,3.00,4.50){
            $4 \ULine{\GenA_2 \LeafSimple \LeafSimple} \ULine{\GenA_1 \ULine{\GenA_2
            \GenA_0 \LeafSimple}} \LeafSimple \LeafSimple$
        };
        \node[GraphLabeledVertex](7)at(2.00,1.75,1.50){
            $4 \ULine{\GenA_2 \LeafSimple \ULine{\GenA_1 \LeafSimple}} \ULine{\GenA_2
            \LeafSimple \LeafSimple} \GenA_0 \LeafSimple$
        };
        \node[GraphLabeledVertex](8)at(2.00,1.75,4.00){
            $4 \ULine{\GenA_2 \LeafSimple \ULine{\GenA_1 \LeafSimple}} \ULine{\GenA_2
            \LeafSimple \GenA_0} \LeafSimple \LeafSimple$
        };
        \node[GraphLabeledVertex](9)at(2.00,1.75,4.50){
            $4 \ULine{\GenA_2 \LeafSimple \ULine{\GenA_1 \LeafSimple}} \ULine{\GenA_2
            \GenA_0 \LeafSimple} \LeafSimple \LeafSimple$
        };
        \node[GraphLabeledVertex](10)at(2.00,3.00,1.75){
            $4 \ULine{\GenA_2 \LeafSimple \ULine{\GenA_1 \ULine{\GenA_2 \LeafSimple
            \LeafSimple}}} \GenA_0 \LeafSimple \LeafSimple$
        };
        \node[GraphLabeledVertex](11)at(2.00,3.00,4.00){
            $4 \ULine{\GenA_2 \LeafSimple \ULine{\GenA_1 \ULine{\GenA_2 \LeafSimple
            \GenA_0}}} \LeafSimple \LeafSimple \LeafSimple$
        };
        \node[GraphLabeledVertex](12)at(2.00,3.00,4.50){
            $4 \ULine{\GenA_2 \LeafSimple \ULine{\GenA_1 \ULine{\GenA_2 \GenA_0
            \LeafSimple}}} \LeafSimple \LeafSimple \LeafSimple$
        };
        \node[GraphLabeledVertex](13)at(2.50,1.75,1.50){
            $4 \ULine{\GenA_2 \ULine{\GenA_1 \LeafSimple} \LeafSimple} \ULine{\GenA_2
            \LeafSimple \LeafSimple} \GenA_0 \LeafSimple$
        };
        \node[GraphLabeledVertex](14)at(2.50,1.75,4.00){
            $4 \ULine{\GenA_2 \ULine{\GenA_1 \LeafSimple} \LeafSimple} \ULine{\GenA_2
            \LeafSimple \GenA_0} \LeafSimple \LeafSimple$
        };
        \node[GraphLabeledVertex](15)at(2.50,1.75,4.50){
            $4 \ULine{\GenA_2 \ULine{\GenA_1 \LeafSimple} \LeafSimple} \ULine{\GenA_2
            \GenA_0 \LeafSimple} \LeafSimple \LeafSimple$
        };
        \node[GraphLabeledVertex](16)at(2.50,2.00,1.75){
            $4 \ULine{\GenA_2 \ULine{\GenA_1 \LeafSimple} \ULine{\GenA_2 \LeafSimple
            \LeafSimple}} \GenA_0 \LeafSimple \LeafSimple$
        };
        \node[GraphLabeledVertex](17)at(2.50,2.00,4.00){
            $4 \ULine{\GenA_2 \ULine{\GenA_1 \LeafSimple} \ULine{\GenA_2 \LeafSimple
            \GenA_0}} \LeafSimple \LeafSimple \LeafSimple$
        };
        \node[GraphLabeledVertex](18)at(2.50,2.00,4.50){
            $4 \ULine{\GenA_2 \ULine{\GenA_1 \LeafSimple} \ULine{\GenA_2 \GenA_0
            \LeafSimple}} \LeafSimple \LeafSimple \LeafSimple$
        };
        \node[GraphLabeledVertex](19)at(2.50,3.00,1.75){
            $4 \ULine{\GenA_2 \ULine{\GenA_1 \ULine{\GenA_2 \LeafSimple \LeafSimple}}
            \LeafSimple} \GenA_0 \LeafSimple \LeafSimple$
        };
        \node[GraphLabeledVertex](20)at(2.50,3.00,2.00){
            $4 \ULine{\GenA_2 \ULine{\GenA_1 \ULine{\GenA_2 \LeafSimple \LeafSimple}}
            \GenA_0} \LeafSimple \LeafSimple \LeafSimple$
        };
        \node[GraphLabeledVertex](21)at(2.50,3.00,4.00){
            $4 \ULine{\GenA_2 \ULine{\GenA_1 \ULine{\GenA_2 \LeafSimple \GenA_0}}
            \LeafSimple} \LeafSimple \LeafSimple \LeafSimple$
        };
        \node[GraphLabeledVertex](22)at(2.50,3.00,4.50){
            $4 \ULine{\GenA_2 \ULine{\GenA_1 \ULine{\GenA_2 \GenA_0 \LeafSimple}}
            \LeafSimple} \LeafSimple \LeafSimple \LeafSimple$
        };
        \draw[GraphArc](1)--(2);
        \draw[GraphArc](1)--(4);
        \draw[GraphArc](1)--(7);
        \draw[GraphArc](2)--(3);
        \draw[GraphArc](2)--(5);
        \draw[GraphArc](2)--(8);
        \draw[GraphArc](3)--(6);
        \draw[GraphArc](3)--(9);
        \draw[GraphArc](4)--(5);
        \draw[GraphArc](4)--(10);
        \draw[GraphArc](5)--(6);
        \draw[GraphArc](5)--(11);
        \draw[GraphArc](6)--(12);
        \draw[GraphArc](7)--(8);
        \draw[GraphArc](7)--(10);
        \draw[GraphArc](7)--(13);
        \draw[GraphArc](8)--(9);
        \draw[GraphArc](8)--(11);
        \draw[GraphArc](8)--(14);
        \draw[GraphArc](9)--(12);
        \draw[GraphArc](9)--(15);
        \draw[GraphArc](10)--(11);
        \draw[GraphArc](10)--(19);
        \draw[GraphArc](11)--(12);
        \draw[GraphArc](11)--(21);
        \draw[GraphArc](12)--(22);
        \draw[GraphArc](13)--(14);
        \draw[GraphArc](13)--(16);
        \draw[GraphArc](14)--(15);
        \draw[GraphArc](14)--(17);
        \draw[GraphArc](15)--(18);
        \draw[GraphArc](16)--(17);
        \draw[GraphArc](16)--(19);
        \draw[GraphArc](17)--(18);
        \draw[GraphArc](17)--(21);
        \draw[GraphArc](18)--(22);
        \draw[GraphArc](19)--(20);
        \draw[GraphArc](20)--(21);
        \draw[GraphArc](21)--(22);
    \end{tikzpicture}
}
    \caption[]{
        The Hasse diagram of the leaning forest $\SignatureExample$-easterly wind poset of
        $\GenA_4 \GenA_2 \GenA_1 \GenA_2 \GenA_0$.
    }
    \label{fig:example_forest_poset}
\end{figure}

For any $w \in \Signature^*$, the \Def{leaning forest $\Signature$-easterly wind poset of
$w$} is the balanced forest $\Bra{1}$-tilted $\Signature$-easterly wind poset of $w$, as
defined in Section~\ref{subsubsec:intervals_forests}. Figure~\ref{fig:example_forest_poset}
shows the Hasse diagram of such poset.

\begin{Statement}{Proposition}{prop:all_leaning_forests}
    Let $\Signature$ be a signature and let $w$ be a word on $\Signature$ of length $n \in
    \N$. The leaning forest $\Signature$-easterly wind poset of $w$ contains all leaning
    $\Signature$-forests having $n \Conc w$ as decoration word.
\end{Statement}
\begin{Proof}
    Let us prove that for any leaning forest $\ForestF$ having $n \Conc w$ as decoration
    word, $\ForestMin \App w \EasterlyWindLeq \ForestF$. For this, let $i$ be an
    internal node of $\ForestF$ different from the root, and let $\Par{i', j, i}$ be the
    parent edge of $i$ in $\ForestF$. From the definition of $\ForestMin \App w$, the
    parent edge of the internal node $i$ in $\ForestMin \App w$ is $\Par{1, i - 1, i}$.
    Now, if $i' \geq 2$, then the parent edge of $i$ in $\ForestF$ dominates that of $i$
    in $\ForestMin \App w$. Otherwise, we have $i' = 1$. Let us prove in this case that
    $j \leq i - 1$. Indeed, assume by contradiction that $j \geq i$. Under this assumption,
    among the first $j - 1$ children of the root of $\ForestF$, since they comprise only $i
    - 2$ internal nodes, there must be at least one leaf. Therefore, $\ForestF$ would not be
    $\Bra{1}$-tilted, contradicting our hypotheses. Hence, even in this case, the parent
    edge of $i$ in $\ForestF$ dominates that of $i$ in $\ForestMin \App w$. It follows
    now from Lemma~\ref{lem:order_domination} that $\ForestMin \App w \EasterlyWindLeq
    \ForestF$.
\end{Proof}

Proposition~\ref{prop:all_leaning_forests} provides the notable property that any leaning
$\Signature$-forest belong to a leaning forest $\Signature$-easterly wind poset.

\subsubsection{Restrictions} \label{subsubsec:restrictions}
Let $\ForestF$ be a leaning $\Signature$-forest of size $n$. Given a subset $I$ of $[n]$,
let $\bar{I} := \Bra{i + 1 : i \in I}$. The \Def{restriction} $\ForestF \App I $of
$\ForestF$ on $I$ is the leaning $\Signature$-forest obtained by keeping only the internal
nodes from the set $\{1\} \cup \Bar{I}$ of $\ForestF$ and their adjacent edges, and by
setting $\# I$ as the decoration of the root. For instance, by considering the leaning
$\SignatureExample$-forest
\begin{equation} \label{equ:example_restriction}
    \ForestF :=
    5 \,
    \GenA_0
    \ULine{\GenA_2 \LeafSimple \ULine{\GenA_1 \GenA_0}}
    \ULine{\GenA_3 \LeafSimple \LeafSimple \LeafSimple}
    \LeafSimple
    \LeafSimple,
\end{equation}
we have
\begin{math}
    \ForestF \App \Bra{1, 4}
    = 2 \, \GenA_0 \GenA_0,
\end{math}
\begin{math}
    \ForestF \App \Bra{2, 3}
    = 2 \, \ULine{\GenA_2 \LeafSimple \ULine{\GenA_1 \LeafSimple}} \LeafSimple,
\end{math}
and
\begin{math}
    \ForestF \App \Bra{2, 4, 5}
    = 3 \,
    \ULine{\GenA_2 \LeafSimple \LeafSimple}
    \GenA_0
    \ULine{\GenA_3 \LeafSimple \LeafSimple \LeafSimple}.
\end{math}

When $I$ is an interval of $[n]$, let $\theta_I : \bar{I} \to \Han{\# I}$ be the map defined
for any $i \in \bar{I}$ by $\theta_I \App i := i - \min I + 1$. Observe that for any $i \in
\bar{I}$, the internal node $i$ of $\ForestF$ gives rise to the internal node $\theta_I \App
i$ in $\ForestF \App I$. This map will be used to lighten the notation during the proof of
the following lemma.

\begin{Statement}{Lemma}{lem:restriction_order}
    Let $\Signature$ be a signature, and let $\ForestF_1$ and $\ForestF_2$ be two leaning
    $\Signature$-forests of the same size $n$. If $\ForestF_1 \EasterlyWindLeq \ForestF_2$
    and $I$ is an interval of $[n]$, then $\ForestF_1 \App I \EasterlyWindLeq \ForestF_2
    \App I$.
\end{Statement}
\begin{Proof}
    Assume that $\ForestF_1 \EasterlyWindLeq \ForestF_2$ and that $I$ is an interval of
    $[n]$. As the case where $I$ is empty is immediate, we assume that $I$ is nonempty. Let
    $i$ be an internal node of both $\ForestF_1$ and $\ForestF_2$ belonging to $\bar{I}$.
    Let also $i_1$ (resp.\ $i_2$) be the parent of $i$ in $\ForestF_1$ (resp.\
    $\ForestF_2$). We have several cases to explore depending on whether $i_1$ and $i_2$
    belong to $\bar{I}$.
    \begin{enumerate}
        \item Assume first that $i_1 \in \bar{I}$. By Lemma~\ref{lem:order_domination}, the
        parent edge of $i$ in $\ForestF_1$ is dominated by that of $i$ in $\ForestF_2$. In
        particular, this implies that $i_1 \leq i_2 < i$ so that, since $I$ is an interval,
        $i_2 \in \bar{I}$. For this reason, we have also
        \begin{math}
            \LocalPosition \App \ULine{\ForestF_2 \App I} \App \ULine{\theta_I \App i}
            =
            \LocalPosition \App \ForestF_2 \App i.
        \end{math}
        Similarly, since $i_1 \in \bar{I}$,
        \begin{math}
            \LocalPosition \App \ULine{\ForestF_1 \App I} \App \ULine{\theta_I \App i}
            =
            \LocalPosition \App \ForestF_1 \App i.
        \end{math}
        These properties imply that the parent edge of $\theta_I \App i$ in $\ForestF_1 \App
        I$ is dominated by that of $\theta_I \App i$ in~$\ForestF_2 \App I$.
        \item Assume now that $i_1 \notin \bar{I}$ and $i_2 \notin \bar{I}$. By definition
        of the restriction operation, the parent of $i$ is the internal node $1$ in both
        $\ForestF_1 \App I$ and $\ForestF_2 \App I$. Moreover, since both $\ForestF_1 \App
        I$ and $\ForestF_2 \App I$ are leaning,
        \begin{math}
            \LocalPosition \App \ULine{\ForestF_1 \App I} \App \ULine{\theta_I \App i}
            \geq
            \LocalPosition \App \ULine{\ForestF_2 \App I} \App \ULine{\theta_I \App i}.
        \end{math}
        Therefore, the parent edge of $\theta_I \App i$ in $\ForestF_1 \App I$ is dominated
        by that of $\theta_I \App i$ in $\ForestF_2 \App I$.
        \item In the remaining case, $i_1 \notin \bar{I}$ and $i_2 \in \bar{I}$. By
        definition of the restriction operation, the parent of $i$ is the internal node $1$
        in $\ForestF_1$. Moreover, since $i_2 \in \bar{I}$, $\theta_I \App i_2 \geq 2$. This
        shows that the parent edge of $\theta_I \App i$ in $\ForestF_1 \App I$ is dominated
        by that of $\theta_I \App i$ in $\ForestF_2 \App I$.
    \end{enumerate}
    We have shown that for any internal node $\theta_I \App i$ of both $\ForestF_1 \App I$
    and $\ForestF_2 \App I$, the parent edge of $\theta_I \App i$ in $\ForestF_1 \App I$ is
    dominated by that of $\theta_I \App i$ in $\ForestF_2 \App I$. Therefore, by
    Lemma~\ref{lem:order_domination}, $\ForestF_1 \App I \EasterlyWindLeq \ForestF_2 \App
    I$.
\end{Proof}

Let us introduce two specific restrictions. Given a leaning $\Signature$-forest $\ForestF$
of size $n$, for any $k \in \HanL{n}$, let $\Top \App k \App \ForestF := \ForestF \App
\Han{1, k}$ and $\Bottom \App k \App \ForestF := \ForestF \App \Han{k + 1, n}$. We call the
first restriction the \Def{$k$-top restriction} of $\ForestF$ and the second, the
\Def{$k$-bottom restriction} of $\ForestF$. For instance, by considering the leaning
$\SignatureExample$-forest $\ForestF$ of~\eqref{equ:example_restriction}, we have
\begin{enumerate}[label=$\bullet$,before=\begin{multicols}{2},after=\end{multicols}]
    \item
    \begin{math}
        \Par{\Top \App 0 \App \ForestF, \Bottom \App 0 \App \ForestF}
        =
        \Par{0, \ForestF};
    \end{math}
    \item
    \begin{math}
        \Par{
            \Top \App 1 \App \ForestF,
            \Bottom \App 1 \App \ForestF
        }
        =
        \Par{
            1 \, \GenA_0,
            4 \,
            \ULine{\GenA_2 \LeafSimple \ULine{\GenA_1 \GenA_0}}
            \ULine{\GenA_3 \LeafSimple \LeafSimple \LeafSimple}
            \LeafSimple
            \LeafSimple
        };
    \end{math}
    \item
    \begin{math}
        \Par{
            \Top \App 2 \App \ForestF,
            \Bottom \App 2 \App \ForestF
        }
        =
        \Par{
            2 \, \GenA_0 \ULine{\GenA_2 \LeafSimple \LeafSimple},
            3 \,
            \ULine{\GenA_1 \GenA_0}
            \ULine{\GenA_3 \LeafSimple \LeafSimple \LeafSimple}
            \LeafSimple
        };
    \end{math}
    \item
    \begin{math}
        \Par{
            \Top \App 3 \App \ForestF,
            \Bottom \App 3 \App \ForestF
        }
        =
        \Par{
            3 \, \GenA_0 \ULine{\GenA_2 \LeafSimple \ULine{\GenA_1 \LeafSimple}},
            2 \,
            \GenA_0
            \ULine{\GenA_3 \LeafSimple \LeafSimple \LeafSimple}
        };
    \end{math}
    \item
    \begin{math}
        \Par{
            \Top \App 4 \App \ForestF,
            \Bottom \App 4 \App \ForestF
        }
        \! = \!
        \Par{
            4 \, \GenA_0 \ULine{\GenA_2 \LeafSimple \ULine{\GenA_1 \GenA_0}}
            \LeafSimple
            \LeafSimple
            \LeafSimple,
            1 \, \ULine{\GenA_3 \LeafSimple \LeafSimple \LeafSimple}
        };
    \end{math}
    \item
    \begin{math}
        \Par{
            \Top \App 5 \App \ForestF
            \Bottom \App 5 \App \ForestF
        }
        =
        \Par{\ForestF, 0}.
    \end{math}
\end{enumerate}

\subsubsection{Over and under operations} \label{subsubsec:over_under}
Let $\Over$ be the \Def{over} operation on leaning $\Signature$-forests defined, for any
leaning $\Signature$-forests $\ForestF_1$ and $\ForestF_2$, by $\ForestF_1 \Over \ForestF_2
:= \Tilt \App \Bra{1} \App \ULine{\ForestF_1 \Conc \ForestF_2}$ where $\Tilt$ is the tilting
map defined in Section~\ref{subsubsec:tilting_map} and $\Conc$ is the concatenation
operation defined in Section~\ref{subsubsec:forests}. For instance, on leaning
$\SignatureExample$-forests,
\begin{equation}
    \ColA{
        3 \,
        \ULine{\GenA_1 \LeafSimple} \ULine{\GenA_2 \LeafSimple
        \ULine{\GenA_2 \LeafSimple \LeafSimple}} \LeafSimple
    }
    \; \Over \;
    \bm{
        \ColB{
            5 \,
            \ULine{\GenA_3 \LeafSimple \LeafSimple \LeafSimple}
            \ULine{\GenA_2 \LeafSimple \LeafSimple}
            \ULine{\GenA_1 \ULine{\GenA_2 \LeafSimple \LeafSimple}}
            \ULine{\GenA_1 \LeafSimple}
            \LeafSimple
        }
    }
    =
    8 \,
    \ColA{
        \ULine{\GenA_1 \LeafSimple} \ULine{\GenA_2 \LeafSimple
        \ULine{\GenA_2 \LeafSimple \LeafSimple}}
    }
    \bm{
        \ColB{
            \ULine{\GenA_3 \LeafSimple \LeafSimple \LeafSimple}
            \ULine{\GenA_2 \LeafSimple \LeafSimple}
            \ULine{\GenA_1 \ULine{\GenA_2 \LeafSimple \LeafSimple}}
            \ULine{\GenA_1 \LeafSimple}
            \LeafSimple
        }
    }
    \ColA{\LeafSimple}.
\end{equation}

Let $\ForestF_1$ and $\ForestF_2$ be two leaning $\Signature$-forests such that $\ForestF_1$
is of size $n_1$ and of length $\ell_1$. From the definition of the over operation, each
edge $\Par{i', j, i}$ of $\ForestF_1 \Over \ForestF_2$ obeys to the following rules:
\begin{enumerate}
    \item if $i' < i \leq n_1 + 1$ (that is, both $i'$ and $i$ come from internal nodes of
    $\ForestF_1$), then $\Par{i', j, i}$ is an edge of~$\ForestF_1$;
    \item if $n_1 + 2 \leq i' < i$ (that is, both $i'$ and $i$ come from internal nodes of
    $\ForestF_2$), then $\Par{i' - n_1, j, i - n_1}$ is an edge of~$\ForestF_2$;
    \item otherwise, we have $i' \leq n_1 + 1$ (that is, $i'$ comes from an internal node of
    $\ForestF_1$) and $n_1 + 2 \leq i$ (that is, $i$ comes from an internal node of
    $\ForestF_2$). In this case, we have $i' = 1$ (that is, $i'$ is the root of
    $\ForestF_1$) and that $\Par{1, j - \ell_1, i - n_1}$ is an edge of~$\ForestF_2$.
\end{enumerate}

Similarly, let $\Under$ be the \Def{under} operation on leaning $\Signature$-forests such
that, for any leaning $\Signature$-forests $\ForestF_1$ and $\ForestF_2$, by denoting by
$n_1$ the size of $\ForestF_1$, by $n_2$ the size of $\ForestF_2$, and by $r$ the number of
extreme leaves of $\ForestF_1$, $\ForestF_1 \Under \ForestF_2$ is the leaning
$\Signature$-forest built by grafting, for any $j \in \Han{n_2 + r}$, the $j$-th subterm of
the root of $\ForestF_2 \Conc \ULine{\Corolla \App r}$ onto the $j$-th extreme leaf of
$\ForestF_1 \Conc \ULine{\Corolla \App n_2}$. For instance, on leaning
$\SignatureExample$-forests,
\begin{equation}
    \ColA{
        3 \,
        \ULine{\GenA_1 \LeafSimple} \ULine{\GenA_2 \LeafSimple
        \ULine{\GenA_2 \LeafSimple \LeafSimple}} \LeafSimple
    }
    \; \Under \;
    \bm{
        \ColB{
            5 \,
            \ULine{\GenA_3 \LeafSimple \LeafSimple \LeafSimple}
            \ULine{\GenA_2 \LeafSimple \LeafSimple}
            \ULine{\GenA_1 \ULine{\GenA_2 \LeafSimple \LeafSimple}}
            \ULine{\GenA_1 \LeafSimple}
            \LeafSimple
        }
    }
    =
    8 \,
    \ColA{
        \ULine{\GenA_1 \LeafSimple}
        \ULine{
            \GenA_2
            \LeafSimple
            \ULine{
                \GenA_2
                \bm{\ColB{\ULine{\GenA_3 \LeafSimple \LeafSimple \LeafSimple}}}
                \bm{\ColB{\ULine{\GenA_2 \LeafSimple \LeafSimple}}}
            }
        }
        \bm{\ColB{\ULine{\GenA_1 \ULine{\GenA_2 \LeafSimple \LeafSimple}}}}
    }
    \bm{\ColB{\ULine{\GenA_1 \LeafSimple}}}
    \LeafSimple
    \LeafSimple
    \LeafSimple
    \LeafSimple.
\end{equation}

Let $\ForestF_1$ and $\ForestF_2$ be two leaning $\Signature$-forests such that $\ForestF_1$
is of size $n_1$ and $\ForestF_2$ is of size $n_2$. From the definition of the under
operation, each edge $\Par{i', j, i}$ of $\ForestF_1 \Under \ForestF_2$ obeys to the
following rules:
\begin{enumerate}
    \item if $i' < i \leq n_1 + 1$ (that is, both $i'$ and $i$ come from internal nodes of
    $\ForestF_1$), then $\Par{i', j, i}$ is an edge of~$\ForestF_1$;
    \item if $n_1 + 2 \leq i' < i$ (that is, both $i'$ and $i$ come from internal nodes of
    $\ForestF_2$), then $\Par{i' - n_1, j, i - n_1}$ is an edge of~$\ForestF_2$;
    \item otherwise, we have $i' \leq n_1 + 1$ (that is, $i'$ comes from an internal node of
    $\ForestF_1$) and $n_1 + 2 \leq i$ (that is, $i$ comes from an internal node of
    $\ForestF_2$). In this case, $\Par{1, j', i - n_1}$ is an edge of~$\ForestF_2$ where in
    $\ForestF_1 \Conc \ULine{\Corolla \App n_2}$, the child of $i'$ at position $j$ is the
    $j'$-th extreme leaf of this $\Signature$-forest.
\end{enumerate}

\begin{Statement}{Lemma}{lem:over_under_order_compatibility}
    Let $\Signature$ be a signature and let $\ForestF_1$, $\ForestF_1'$, $\ForestF_2$, and
    $\ForestF_2'$ be leaning $\Signature$-forests. If $\ForestF_1 \EasterlyWindLeq
    \ForestF_1'$ and $\ForestF_2 \EasterlyWindLeq \ForestF_2'$, then
    \begin{enumerate}[label=(\roman*),before=\begin{multicols}{2},after=\end{multicols}]
        \item \label{item:over_under_order_compatibility_1}
        $\ForestF_1 \Over \ForestF_2 \EasterlyWindLeq \ForestF_1' \Over \ForestF_2'$;
        \item \label{item:over_under_order_compatibility_2}
        $\ForestF_1 \Under \ForestF_2 \EasterlyWindLeq \ForestF_1' \Under \ForestF_2'$.
    \end{enumerate}
\end{Statement}
\begin{Proof}
    Le us assume that $\ForestF_1 \EasterlyWindLeq \ForestF_1'$ and $\ForestF_2
    \EasterlyWindLeq \ForestF_2'$. By Lemma~\ref{lem:order_domination}, for any internal
    node $i$ of $\ForestF_1$ (resp.\ $\ForestF_2$), the parent edge of $i$ in $\ForestF_1$
    (resp.\ $\ForestF_2$) is dominated by that of $i$ in $\ForestF_1'$ (resp.\
    $\ForestF_2'$).

    Now, by the previous description of the edges of a leaning $\Signature$-forest obtained
    from the over operation applied on two leaning $\Signature$-forests, together with the
    fact that for any leaning $\Signature$-forests $\ForestG$ and $\ForestG'$, $\ForestG
    \EasterlyWindLeq \ForestG'$ implies that the length of $\ForestG$ is greater than or
    equal to the length of $\ForestG'$, for any internal node $i$ of $\ForestF_1 \Over
    \ForestF_2$ the parent edge of $i$ in $\ForestF_1 \Over \ForestF_2$ is dominated by that
    of $i$ in $\ForestF_1' \Over \ForestF_2'$. Therefore, by
    Lemma~\ref{lem:order_domination}, \ref{item:over_under_order_compatibility_1} holds.

    Similarly, by the previous description of the edges of a leaning $\Signature$-forest
    obtained from the under operation applied on two leaning $\Signature$-forests, together
    with the fact that for any leaning $\Signature$-forests $\ForestG$ and $\ForestG'$,
    $\ForestG \EasterlyWindLeq \ForestG'$ implies that the number of extreme leaves of
    $\ForestG$ is smaller than or equal to the number of extreme leaves of $\ForestG'$, for
    any internal node $i$ of $\ForestF_1 \Over \ForestF_2$ the parent edge of $i$ in
    $\ForestF_1 \Under \ForestF_2$ is dominated by that of $i$ in $\ForestF_1' \Under
    \ForestF_2'$. Therefore, by Lemma~\ref{lem:order_domination},
    \ref{item:over_under_order_compatibility_2} holds.
\end{Proof}

For any leaning $\Signature$-forest $\ForestF$ of size $n$ and any $k \in \HanL{0, n}$, let
\begin{equation}
    \Over \App k \App \ForestF
    :=
    \Top \App k \App \ForestF \Over \Bottom \App k \App \ForestF
\end{equation}
and
\begin{equation}
    \Under \App k \App \ForestF
    :=
    \Top \App k \App \ForestF \Under \Bottom \App k \App \ForestF.
\end{equation}

\begin{Statement}{Lemma}{lem:over_under_interval}
    Let $\Signature$ be a signature and let $\ForestF$ be a leaning $\Signature$-forest of
    size $n$. For any $k \in \HanL{0, n}$,
    \begin{equation}
        \Over \App k \App \ForestF
        \EasterlyWindLeq \ForestF \EasterlyWindLeq
        \Under \App k \App \ForestF.
    \end{equation}
\end{Statement}
\begin{Proof}
    From the definition of $\Over \App k$ (resp.\ $\Under \App k$) and the previous
    description of the edges of a leaning $\Signature$-forest obtained from the over (resp.\
    under) operation applied on two leaning $\Signature$-forests, the leaning
    $\Signature$-forests $\ForestF$ and $\Over \App k \App \ForestF$ (resp.\ $\Under \App k
    \App \ForestF$) have the same edges, except possibly the parent edges of the internal
    nodes $i$ with $i \geq k + 2$, which are of the form $\Par{i', j, i}$ in $\Over \App k
    \App \ForestF$ (resp.\ $\Under \App k \App \ForestF$) where $j$ is an integer, and of
    the form $\Par{i'', j', i}$ in $\ForestF$ with $i'' > i'$ (resp.\ $i'' < i'$), or both
    $i'' = i'$ and $j' = j$. By Lemma~\ref{lem:order_domination}, this implies that $\Over
    \App k \App \ForestF \EasterlyWindLeq \ForestF$ (resp.\ $\ForestF \EasterlyWindLeq
    \Under \App k \App \ForestF$).
\end{Proof}

\begin{Statement}{Lemma}{lem:over_under_restriction}
    Let $\Signature$ be a signature, and let $\ForestF_1$ and $\ForestF_2$ be two leaning
    $\Signature$-forests. By denoting by $k$ the size of $\ForestF_1$, for any leaning
    $\Signature$-forest $\ForestF$, the two following properties hold:
    \begin{enumerate}[label=(\roman*)]
        \item \label{item:over_under_restriction_1}
        $\ForestF_1 \Over \ForestF_2 \EasterlyWindLeq \ForestF$ if and only if $\ForestF_1
        \EasterlyWindLeq \Top \App k \App \ForestF$ and $\ForestF_2 \EasterlyWindLeq \Bottom
        \App k \App \ForestF$;
        \item \label{item:over_under_restriction_2}
        $\ForestF \EasterlyWindLeq \ForestF_1 \Under \ForestF_2$ if and only if $\Top \App k
        \App \ForestF \EasterlyWindLeq \ForestF_1$ and $\Bottom \App k \App \ForestF
        \EasterlyWindLeq \ForestF_2$.
    \end{enumerate}
\end{Statement}
\begin{Proof}
    Assume first that $\ForestF_1 \Over \ForestF_2 \EasterlyWindLeq \ForestF$ (resp.\
    $\ForestF \EasterlyWindLeq \ForestF_1 \Under \ForestF_2$). By
    Lemma~\ref{lem:restriction_order}, we have
    \begin{math}
        \Top \App k \App \ULine{\ForestF_1 \Over \ForestF_2}
        \EasterlyWindLeq
        \Top \App k \App \ForestF
    \end{math}
    (resp.\
    \begin{math}
        \Top \App k \App \ForestF
        \EasterlyWindLeq
        \Top \App k \App \ULine{\ForestF_1 \Under \ForestF_2}
    \end{math})
    and
    \begin{math}
        \Bottom \App k \App \ULine{\ForestF_1 \Over \ForestF_2}
        \EasterlyWindLeq
        \Bottom \App k \App \ForestF
    \end{math}
    (resp.\
    \begin{math}
        \Bottom \App k \App \ForestF
        \EasterlyWindLeq
        \Bottom \App k \App \ULine{\ForestF_1 \Under \ForestF_2}
    \end{math}).
    Now, by definition of the over (resp.\ under) operation, we have $\Top \App k \App
    \ULine{\ForestF_1 \Over \ForestF_2} = \ForestF_1$ (resp.\ $\Top \App k \App
    \ULine{\ForestF_1 \Under \ForestF_2} = \ForestF_1$) and $\Bottom \App k \App
    \ULine{\ForestF_1 \Over \ForestF_2} = \ForestF_2$ (resp.\ $\Bottom \App k \App
    \ULine{\ForestF_1 \Under \ForestF_2} = \ForestF_2$). Hence, $\ForestF_1 \EasterlyWindLeq
    \Top \App k \App \ForestF$ (resp.\ $\Top \App k \App \ForestF \EasterlyWindLeq
    \ForestF_1$) and $\ForestF_2 \EasterlyWindLeq \Bottom \App k \App \ForestF$ (resp.\
    $\Bottom \App k \App \ForestF \EasterlyWindLeq \ForestF_2$). This proves the direct
    implication of~\ref{item:over_under_restriction_1} (resp.\
    \ref{item:over_under_restriction_2}).

    Assume conversely that $\ForestF_1 \EasterlyWindLeq \Top \App k \App \ForestF$ (resp.\
    $\Top \App k \App \ForestF \EasterlyWindLeq \ForestF_1$) and $\ForestF_2
    \EasterlyWindLeq \Bottom \App k \App \ForestF$ (resp.\ $\Bottom \App k \App \ForestF
    \EasterlyWindLeq \ForestF_2$). By Lemma~\ref{lem:over_under_interval}, we have $\Over
    \App k \App \ForestF \EasterlyWindLeq \ForestF$ (resp.\ $\ForestF \EasterlyWindLeq
    \Under \App k \App \ForestF$). Given this, by
    Lemma~\ref{lem:over_under_order_compatibility}, we have $\ForestF_1 \Over \ForestF_2
    \EasterlyWindLeq \Over \App k \App \ForestF$ (resp.\ $\Under \App k \App \ForestF
    \EasterlyWindLeq \ForestF_1 \Under \ForestF_2$). Therefore, by using the transitivity of
    $\EasterlyWindLeq$, this implies that $\ForestF_1 \Over \ForestF_2 \EasterlyWindLeq
    \ForestF$ (resp.\ $\ForestF \EasterlyWindLeq \ForestF_1 \Under \ForestF_2$). The
    converse of~\ref{item:over_under_restriction_1} (resp.\
    \ref{item:over_under_restriction_2}) has been established.
\end{Proof}

\subsubsection{Shuffle product} \label{subsubsec:shuffle_leaning_forest}
Let $\ForestF_1$ and $\ForestF_2$ be two leaning $\Signature$-forests of respective sizes
$n_1$ and $n_2$. The \Def{shifted shuffle} of $\ForestF_1$ and $\ForestF_2$ is the set
$\ForestF_1 \cshuffle \ForestF_2$ of leaning $\Signature$-forests $\ForestF$ of size $n_1 +
n_2$ such that $\Top \App n_1 \App \ForestF = \ForestF_1$ and $\Bottom \App n_1 \App
\ForestF = \ForestF_2$.

For instance, let the leaning $\SignatureExample$-forests
\begin{math}
    \ForestF_1
    :=
    \ColA{4 \, \ULine{\GenA_1 \LeafSimple} \ULine{\GenA_3 \LeafSimple \ULine{\GenA_3
    \LeafSimple \GenA_0 \LeafSimple} \LeafSimple} \LeafSimple \LeafSimple}
\end{math}
and
\begin{math}
    \ForestF_2
    :=
    \ColB{3 \, \ULine{\GenA_2 \LeafSimple \ULine{\GenA_1 \LeafSimple}} \GenA_0 \LeafSimple}.
\end{math}
The set $\ForestF_1 \cshuffle \ForestF_2$ contains exactly the four leaning
$\SignatureExample$-forests
\begin{enumerate}[label=$\bullet$,before=\begin{multicols}{2},after=\end{multicols}]
    \item
    \begin{math}
        7 \,
        \ColA{
            \ULine{\GenA_1 \LeafSimple}
            \ULine{
                \GenA_3 \LeafSimple \ULine{\GenA_3 \LeafSimple \GenA_0 \LeafSimple}
                \LeafSimple
            }
        }
        \ColB{\ULine{\GenA_2 \LeafSimple \ULine{\GenA_1
        \LeafSimple}}} \ColB{\GenA_0} \LeafSimple \LeafSimple \LeafSimple;
    \end{math}
    \item
    \begin{math}
        7 \,
            \ColA{\ULine{\GenA_1 \LeafSimple}}
            \ULine{
                \ColA{\GenA_3 \LeafSimple}
                \ColA{\ULine{\GenA_3 \LeafSimple \GenA_0 \LeafSimple}}
                \ColB{\ULine{\GenA_2 \LeafSimple \ULine{\GenA_1 \LeafSimple}}}
            }
        \ColB{\GenA_0} \LeafSimple \LeafSimple \LeafSimple \LeafSimple;
    \end{math}
    \item
    \begin{math}
        7 \,
        \ColA{
            \ULine{\GenA_1 \LeafSimple}
            \ULine{
                \GenA_3 \LeafSimple \ULine{\GenA_3 \LeafSimple
                \GenA_0
                \ColB{\ULine{\GenA_2 \LeafSimple \ULine{\GenA_1 \LeafSimple}}}}
                \LeafSimple
            }
        }
        \ColB{\GenA_0} \LeafSimple \LeafSimple \LeafSimple \LeafSimple;
    \end{math}
    \item
    \begin{math}
        7 \,
        \ColA{\ULine{\GenA_1 \LeafSimple}}
        \ULine{
            \ColA{\GenA_3 \LeafSimple}
            \ColA{
                \ULine{
                    \GenA_3
                    \LeafSimple
                    \GenA_0
                    \ColB{\ULine{\GenA_2 \LeafSimple \ULine{\GenA_1 \LeafSimple}}}
                }
            }
            \ColB{\GenA_0}
        }
        \LeafSimple \LeafSimple \LeafSimple \LeafSimple \LeafSimple.
    \end{math}
\end{enumerate}
Observe that among these four $\SignatureExample$-forests, the first (resp.\ last) one is
$\ForestF_1 \Over \ForestF_2$ (resp.\ $\ForestF_1 \Under \ForestF_2$).

\begin{Statement}{Proposition}{prop:leaning_forest_shuffle_interval}
    For any signature $\Signature$ and any leaning $\Signature$-forests $\ForestF_1$ and
    $\ForestF_2$, $\ForestF_1 \cshuffle \ForestF_2$ is an interval of the leaning forest
    $\Signature$-easterly wind poset of $w_1 \Conc w_2$, where the decoration word of
    $\ForestF_1$ (resp.\ $\ForestF_2$) is $n_1 \Conc w_1$ (resp.\ $n_2 \Conc w_2$) and $n_1$
    (resp.\ $n_2$) is the size of $\ForestF_1$ (resp.\ $\ForestF_2$). More precisely,
    \begin{equation} \label{equ:leaning_forest_shuffle_interval}
        \ForestF_1 \cshuffle \ForestF_2
        =
        \Han{\ForestF_1 \Over \ForestF_2, \ForestF_1 \Under \ForestF_2}.
    \end{equation}
\end{Statement}
\begin{Proof}
    First of all, by definition of $\cshuffle$, all $\Signature$-forests of $\ForestF_1
    \cshuffle \ForestF_2$ have $\Par{n_1 + n_2} \Conc w_1 \Conc w_2$ as decoration word.
    Hence, by Proposition~\ref{prop:all_leaning_forests}, $\ForestF_1 \cshuffle \ForestF_2$
    is a subset of the leaning forest $\Signature$-easterly wind poset of $w_1 \Conc w_2$.

    Let $\ForestF$ be an $\Signature$-forest. Assume that $\ForestF$ belongs to the set
    $\ForestF_1 \cshuffle \ForestF_2$. By definition of $\cshuffle$, we have $\Top \App n_1
    \App \ForestF = \ForestF_1$ and $\Bottom \App n_1 \App \ForestF = \ForestF_2$.
    Therefore, by Lemma~\ref{lem:over_under_interval},
    \begin{math}
        \ForestF_1 \Over \ForestF_2
        \EasterlyWindLeq \ForestF \EasterlyWindLeq
        \ForestF_1 \Under \ForestF_2.
    \end{math}
    Conversely, assume that
    \begin{math}
        \ForestF_1 \Over \ForestF_2
        \EasterlyWindLeq \ForestF \EasterlyWindLeq
        \ForestF_1 \Under \ForestF_2.
    \end{math}
    By Lemma~\ref{lem:restriction_order}, we have
    \begin{math}
        \Top \App n_1 \App \ULine{\ForestF_1 \Over \ForestF_2}
        \EasterlyWindLeq \Top \App n_1 \App \ForestF \EasterlyWindLeq
        \Top \App n_1 \ULine{\ForestF_1 \Under \ForestF_2}
    \end{math}
    and
    \begin{math}
        \Bottom \App n_1 \App \ULine{\ForestF_1 \Over \ForestF_2}
        \EasterlyWindLeq
        \Bottom \App n_1 \App \ForestF
        \EasterlyWindLeq
        \Bottom \App n_1 \ULine{\ForestF_1 \Under \ForestF_2}.
    \end{math}
    Directly from the definitions of the over and under operations, we have $\Top \App n_1
    \App \ULine{\ForestF_1 \Over \ForestF_2} = \ForestF_1 = \Top \App n_1 \App
    \ULine{\ForestF_1 \Under \ForestF_2}$ and $\Bottom \App n_1 \App \ULine{\ForestF_1 \Over
    \ForestF_2} = \ForestF_2 = \Bottom \App n_1 \App \ULine{\ForestF_1 \Under \ForestF_2} =
    \ForestF_2$. Hence, we have $\Top \App n_1 \App \ForestF = \ForestF_1$ and $\Bottom \App
    n_1 \App \ForestF = \ForestF_2$, showing that $\ForestF$ belongs to $\ForestF_1
    \cshuffle \ForestF_2$. Therefore, \eqref{equ:leaning_forest_shuffle_interval} holds.
\end{Proof}

\section{Natural Hopf algebras of nonsymmetric operads} \label{sec:natural_hopf_algebras}
In this final section we build on the preceding material to accomplish one of the main
objectives of this work: introducing new bases for the natural Hopf algebras of free
nonsymmetric operads. Existing descriptions of these Hopf algebras are given with respect to
an elementary basis, under which the product of two basis elements is expressed as a
concatenation of leaning forests. We then construct two additional bases, the fundamental
basis and the homogeneous basis, such that the product of two elements is, respectively, a
shuffle or a specialized concatenation of leaning forests.

\subsection{Nonsymmetric operads and natural Hopf algebras}
In this preliminary section, we recall the main elementary concepts of operad theory and of
the natural Hopf algebra construction.

\subsubsection{Nonsymmetric operads}
We use the standard definitions about nonsymmetric operads (called simply \Def{operads}
here), as found in~\cite{Gir18}. An operad $\Operad$ is above all considered to be a
signature. We denote by
\begin{equation}
    \FullComposition :
    \Operad \App n \to \Operad \App m_1 \to \dots \to \Operad \App m_n
    \to \Operad \App \ULine{m_1 + \dots + m_n}
\end{equation}
the composition map of $\Operad$ defined for any $n, m_1, \dots, m_n \in \N$, and by $\Unit$
the unit of $\Operad$.

Let $\Operad$ be an operad. When each $x \in \Operad$ admits finitely many factorizations of
the form
\begin{math}
    x = \FullComposition \App y \App y_1 \App \dots \App y_{\Arity \App y}
\end{math}
where $y, y_1, \dots, y_{\Arity \App y} \in \Operad$, $\Operad$ is \Def{finitely
factorizable}. When there exists a map $\Deg : \Operad \to \N$ such that $\Deg^{-1}
\App 0 = \Bra{\Unit}$ and, for any $y, y_1, \dots, y_{\Arity \App y} \in \Operad$,
\begin{equation}
    \Deg \App
    \ULine{\FullComposition \App y \App y_1 \App \dots \App y_{\Arity \App y}}
    = \Deg \App y + \Deg \App y_1 + \dots + \Deg \App y_{\Arity \App y},
\end{equation}
the map $\Deg$ is a \Def{grading} of~$\Operad$.

\subsubsection{Natural Hopf algebras of nonsymmetric operads}
The \Def{natural Hopf algebra} ~\cite{vdl04,ML14,BG16,Gir24} of a finitely factorizable
operad $\Operad$ admitting a grading $\Deg$ is the Hopf algebra $\NaturalHopfAlgebra \App
\Operad$ defined as follows.  Let $\Reduced : \Operad^* \to \Par{\Operad \setminus
\Bra{\Unit}}^*$ be the map such that $\Reduced \App w$ is the subword of $w \in \Operad^*$
consisting of its letters different from $\Unit$. Any fixed point of $\Reduced$ is
\Def{reduced}. Let $\NaturalHopfAlgebra \App \Operad$ be the $\K$-linear span of the set
$\Reduced \App \Operad^*$. The bases of $\NaturalHopfAlgebra \App \Operad$ are thus indexed
by $\Reduced \App \Operad^*$, and the \Def{elementary basis} (or \Def{$\BasisE$-basis} for
short) of $\NaturalHopfAlgebra \App \Operad$ is the set $\Bra{\BasisE_w : w \in \Reduced
\App \Operad^*}$.

This vector space is endowed with an associative algebra structure through
the product $\Product$ satisfying, for any $w_1, w_2 \in \Reduced \App \Operad^*$,
\begin{equation}
    \BasisE_{w_1} \Product \BasisE_{w_2} = \BasisE_{w_1 \Conc w_2}.
\end{equation}
The element $\BasisE_\epsilon$ is the identity w.r.t.\ the product $\Product$.

Moreover, $\NaturalHopfAlgebra \App \Operad$ is endowed with the coproduct $\Coproduct$
defined as the unique associative algebra morphism satisfying, for any $x \in \Operad$,
\begin{equation} \label{equ:natural_coproduct}
    \Coproduct \App \BasisE_x =
    \sum_{y \in \Operad}
    \;
    \sum_{w \in \Operad^{\Arity \App y}}
    \;
    \Iverson{
        x = \FullComposition \App y \App \ULine{w \App 1} \App \dots \App
            \ULine{w \App \ULine{\Arity \App y}}
    }
    \;
    \BasisE_{\Reduced \App y} \otimes \BasisE_{\Reduced \App w},
\end{equation}
where $\Iverson{-}$ is the Iverson bracket as defined at the end of
Section~\ref{sec:introduction}. Due to the fact that $\Operad$ is finitely factorizable,
\eqref{equ:natural_coproduct} is a finite sum. This coproduct endows $\NaturalHopfAlgebra
\App \Operad$ with the structure of a bialgebra. By extending additively $\Deg$ on
$\Operad^*$, the map $\Deg$ defines a grading of $\NaturalHopfAlgebra \App \Operad$. Thus,
$\NaturalHopfAlgebra \App \Operad$ admits an antipode and becomes a Hopf algebra.

\subsection{Natural Hopf algebras of free operads}
Here, we begin by describing the free operads on terms and then describe the natural Hopf
algebras of free operads in terms of leaning forests.

\subsubsection{Free operads on terms}
The \Def{free operad on $\Signature$} is the set $\SetTerms \App \Signature$ considered as a
signature through the arity map $\Arity$, with the composition map such that for any
$\TermT, \TermT_1, \dots, \TermT_{\Arity \App \TermT} \in \SetTerms \App \Signature$,
\begin{math}
    \FullComposition \App \TermT \App \TermT_1 \App \dots \App
    \TermT_{\Arity \App \TermT}
\end{math}
is the $\Signature$-term obtained by substituting each leaf of $\TermT$ from left to right
with $\TermT_1$, \dots, $\TermT_{\Arity \App \TermT}$, and with $\LeafSimple$ as unit. For
instance, in $\SetTerms \App \SignatureExample$, we have
\begin{equation}
    \FullComposition
    \, \App \,
    \ColA{
        \ULine{
            \GenA_2
            \ULine{\GenA_1 \LeafSimple} \ULine{\GenA_3 \LeafSimple \LeafSimple \LeafSimple}
        }
    }
    \; \App \;
    \ColB{\ULine{\GenA_2 \LeafSimple \LeafSimple}}
    \, \App \,
    \ColB{\LeafSimple}
    \, \App \,
    \ColB{\ULine{\GenA_1 \ULine{\GenA_1 \LeafSimple}}}
    \, \App \,
    \ColB{\ULine{\GenA_2 \LeafSimple \ULine{\GenA_3 \LeafSimple \LeafSimple \LeafSimple}}}
    =
    \ColA{
        \GenA_2 \ULine{\GenA_1 \ColB{\ULine{\GenA_2 \LeafSimple \LeafSimple}}}
        \ULine{\GenA_3
            \ColB{\LeafSimple}
            \ColB{\ULine{\GenA_1 \ULine{\GenA_1 \LeafSimple}}}
            \ColB{\ULine{\GenA_2 \LeafSimple \ULine{\GenA_3 \LeafSimple \LeafSimple
            \LeafSimple}}}
        }
    }.
\end{equation}

\subsubsection{Hopf algebras on leaning forests} \label{subsubsec:hopf_algebras_forests}
By construction, the Hopf algebra $\NaturalHopfAlgebra \App \ULine{\SetTerms \App
\Signature}$ is graded by $\Deg$ and its bases are indexed by the set of words
$\Par{\TermT_1, \dots, \TermT_k}$ such that $k \geq 0$ and for each $i \in [k]$, $\TermT_i
\in \SetTerms \App \Signature \setminus \Bra{\LeafSimple}$. The map sending such a word
$\Par{\TermT_1, \dots, \TermT_k}$ to the leaning $\Signature$-forest $n \, \TermT_1 \dots
\TermT_k \, \LeafSimple \dots \LeafSimple$ where $n := \Deg \App \TermT_1 + \dots + \Deg
\App \TermT_k$ is a one-to-one correspondence between the set of words on $\SetTerms \App
\Signature \setminus \Bra{\LeafSimple}$ and the set of leaning $\Signature$-forests. For
instance, the word
\begin{math}
    \Par{
        \GenA_2 \LeafSimple \ULine{\GenA_3 \LeafSimple \LeafSimple \LeafSimple},
        \GenA_1 \LeafSimple,
        \GenA_3 \ULine{\GenA_1 \LeafSimple} \LeafSimple \LeafSimple
    }
\end{math}
on $\SetTerms \App \SignatureExample \setminus \Bra{\LeafSimple}$ is sent to the leaning
$\SignatureExample$-forest
\begin{math}
    5 \,
    \ULine{\GenA_2 \LeafSimple \ULine{\GenA_3 \LeafSimple \LeafSimple \LeafSimple}}
    \ULine{\GenA_1 \LeafSimple}
    \ULine{\GenA_3 \ULine{\GenA_1 \LeafSimple} \LeafSimple \LeafSimple}
    \LeafSimple \LeafSimple.
\end{math}
For this reason, through this correspondence, we shall identify words on $\SetTerms \App
\Signature \setminus \Bra{\LeafSimple}$ with leaning $\Signature$-forests.

On the $\BasisE$-basis, the product of $\NaturalHopfAlgebra \App \ULine{\SetTerms \App
\Signature}$ expresses, for any leaning $\Signature$-forests $\ForestF_1$ and $\ForestF_2$,
as
\begin{equation} \label{equ:product_e_basis}
    \BasisE_{\ForestF_1} \Product \BasisE_{\ForestF_2}
    = \BasisE_{\ForestF_1 \Over \ForestF_2},
\end{equation}
where $\Over$ is the over operation defined in Section~\ref{subsubsec:over_under}. For
instance, in $\NaturalHopfAlgebra \App \ULine{\SetTerms \App \SignatureExample}$,
\begin{equation}
    \BasisE_{
        \ColA{
            4 \, \ULine{\GenA_1 \ULine{\GenA_2 \LeafSimple \LeafSimple}}
            \ULine{\GenA_3 \ULine{\GenA_1 \LeafSimple} \LeafSimple \LeafSimple}
            \LeafSimple
            \LeafSimple
        }
    }
    \Product
    \BasisE_{
        \ColB{2 \, \ULine{\GenA_2 \ULine{\GenA_1 \LeafSimple} \LeafSimple} \LeafSimple}
    }
    =
    \BasisE_{
        6 \,
        \ColA{
            \ULine{\GenA_1 \ULine{\GenA_2 \LeafSimple \LeafSimple}}
            \ULine{\GenA_3 \ULine{\GenA_1 \LeafSimple} \LeafSimple \LeafSimple}
        }
        \ColB{\ULine{\GenA_2 \ULine{\GenA_1 \LeafSimple} \LeafSimple}}
        \LeafSimple \LeafSimple \LeafSimple
    }.
\end{equation}

Let us now describe the coproduct of $\NaturalHopfAlgebra \App \ULine{\SetTerms \App
\Signature}$ on the $\BasisE$-basis. Given a leaning $\Signature$-forest $\ForestF$ of size
$n$, a pair $\Par{I_1, I_2}$ of sets is \Def{$\ForestF$-admissible} if $I_1 \sqcup I_2 =
[n]$, for any $i_1 \in I_1$, all ancestors of the internal node $i_1 + 1$ of $\ForestF$
except the root belong to $I_1$, and for any $i_2 \in I_2$, all descendants of the internal
node $i_2 + 1$ of $\ForestF$ belong to $I_2$. This property is denoted by $\Par{I_1, I_2}
\AdmissibleSet \ForestF$. For instance, by considering the leaning
$\SignatureExample$-forest
\begin{math}
    3 \, \ULine{\GenA_3 \LeafSimple \ULine{\GenA_1 \LeafSimple} \LeafSimple}
    \ULine{\GenA_2 \LeafSimple \LeafSimple} \LeafSimple,
\end{math}
the $\ForestF$-admissible pairs of sets are exactly $\Par{\Bra{1, 2, 3}, \emptyset}$,
$\Par{\Bra{1, 2}, \Bra{3}}$, $\Par{\Bra{1, 3}, \Bra{2}}$, $\Par{\Bra{1}, \Bra{2, 3}}$,
$\Par{\Bra{3}, \Bra{1, 2}}$, and $\Par{\emptyset, \Bra{1, 2, 3}}$.

It follows from a description of~\cite{Gir24} of the coproduct of $\NaturalHopfAlgebra \App
\ULine{\SetTerms \App \Signature}$ that for any leaning $\Signature$-forest $\ForestF$ of
size $n$,
\begin{equation}
    \Coproduct \App \BasisE_\ForestF
    =
    \sum_{I_1, I_2 \subseteq [n]}
    \Iverson{\Par{I_1, I_2} \AdmissibleSet \ForestF} \;
    \BasisE_{\ForestF \App I_1} \otimes \BasisE_{\ForestF \App I_2},
\end{equation}
where, for any leaning $\Signature$-forest $\ForestF$ of size $n$ and any subset $I$ of
$[n]$, the notation $\ForestF \App I$ refers to the restriction, as introduced in
Section~\ref{subsubsec:restrictions}. For instance, in $\NaturalHopfAlgebra \App
\ULine{\SetTerms \App \SignatureExample}$,
\begin{align}
    \Coproduct \App
    \BasisE_{
        3 \, \ULine{\GenA_3 \LeafSimple \ULine{\GenA_1 \LeafSimple} \LeafSimple}
        \ULine{\GenA_2 \LeafSimple \LeafSimple} \LeafSimple
    }
    & =
    \BasisE_0
    \otimes
    \BasisE_{
        3 \, \ULine{\GenA_3 \LeafSimple \ULine{\GenA_1 \LeafSimple} \LeafSimple}
        \ULine{\GenA_2 \LeafSimple \LeafSimple} \LeafSimple
    }
    +
    \BasisE_{1 \, \ULine{\GenA_3 \LeafSimple \LeafSimple \LeafSimple}}
    \otimes
    \BasisE_{2 \, \ULine{\GenA_1 \LeafSimple} \ULine{\GenA_2 \LeafSimple \LeafSimple}}
    \\
    & \quad +
    \BasisE_{1 \, \ULine{\GenA_2 \LeafSimple \LeafSimple}}
    \otimes
    \BasisE_{2 \, \ULine{\GenA_3 \LeafSimple \ULine{\GenA_1 \LeafSimple} \LeafSimple}
    \LeafSimple}
    +
    \BasisE_{2 \, \ULine{\GenA_3 \LeafSimple \ULine{\GenA_1 \LeafSimple}
    \LeafSimple}\LeafSimple}
    \otimes
    \BasisE_{1 \, \ULine{\GenA_2 \LeafSimple \LeafSimple}}
    \notag
    \\
    & \quad +
    \BasisE_{2 \, \ULine{\GenA_3 \LeafSimple \LeafSimple \LeafSimple} \ULine{\GenA_2
    \LeafSimple \LeafSimple}}
    \otimes
    \BasisE_{1 \, \ULine{\GenA_1 \LeafSimple}}
    +
    \BasisE_{
        3 \, \ULine{\GenA_3 \LeafSimple \ULine{\GenA_1 \LeafSimple} \LeafSimple}
        \ULine{\GenA_2 \LeafSimple \LeafSimple} \LeafSimple
    }
    \otimes
    \BasisE_0.
    \notag
\end{align}

\subsubsection{Noncommutative symmetric functions} \label{subsubsec:ncsf}
In particular, the Hopf algebra of noncommutative symmetric functions
$\NCSF$~\cite{GKLLRT95} can be understood as a natural Hopf algebra of a free operad.
Indeed, $\NCSF$ is isomorphic to the natural Hopf algebra $\NaturalHopfAlgebra \App
\ULine{\SetTerms \App \Signature}$ of the free operad on the signature $\Signature :=
\Bra{\GenA_1}$ where $\Arity \App \GenA_1 = 1$ (see~\cite{Gir24}). Indeed, by encoding a
leaning $\Signature$-forest $\ForestF$ by the integer composition $\Par{r_1, \dots, r_k}$,
$k \geq 0$, such that for any $i \in [k]$, $r_i$ is the degree of the $i$-th subterm of
$\ForestF$, $\NaturalHopfAlgebra \App \ULine{\SetTerms \App \Signature}$ and $\NCSF$ are
defined on the same vector space. Moreover, for any integer compositions $\Par{r_1, \dots,
r_k}$, $k \geq 0$, and $\Par{r'_1, \dots, r'_{k'}}$, $k' \geq 0$, we have
\begin{equation}
    \BasisE_{\Par{r_1, \dots, r_k}}
    \Product
    \BasisE_{\Par{r'_1, \dots, r'_{k'}}}
    =
    \BasisE_{\Par{r_1, \dots, r_k, r'_1, \dots, r'_{k'}}}
\end{equation}
and
\begin{equation}
    \Coproduct \App \BasisE_{\Par{r_1}}
    =
    \sum_{i \in \HanL{r_1}}
    \BasisE_{\Par{i}} \otimes \BasisE_{\Par{r_1 - i}},
\end{equation}
where $\Par{0}$ and the empty integer composition are identified. This product and coproduct
are the ones of $\NCSF$ through its elementary basis.

\subsection{Fundamental and homogeneous bases} \label{subsec:basis_f_h}
We introduce a two new bases of the natural Hopf algebra of a free operad. These bases mimic
well-known constructions of bases in combinatorial Hopf algebras defined by summing on
intervals of particular posets. Here, the posets intervening in this new basis are the
leaning forest $\Signature$-easterly wind posets.

\subsubsection{Fundamental basis}
Let us use the leaning forest $\Signature$-easterly wind posets to build a new basis of
$\NaturalHopfAlgebra \App \ULine{\SetTerms \App \Signature}$. For any leaning
$\Signature$-forest $\ForestF$, let
\begin{equation}
    \BasisF_\ForestF :=
    \sum_{\ForestF' \in \SetLeaningForests \App \Signature}
    \Iverson{\ForestF \EasterlyWindLeq \ForestF'} \;
    \ULine{\Mobius \App \Bra{1} \App \ForestF \App  \ForestF'} \;
    \BasisE_{\ForestF'},
\end{equation}
where $\SetLeaningForests \App \Signature$ is the set defined in
Section~\ref{subsubsec:leaning_forests}, $\EasterlyWindLeq$ is the partial order relation of
the leaning forest $\Signature$-easterly wind posets, and, in accordance with the notations
introduced at the end of Section~\ref{subsubsec:el_shellability}, $\Mobius \App \Bra{1}$ is
the Möbius function of these lattices. For instance, in $\NaturalHopfAlgebra \App
\ULine{\SetTerms \App \SignatureExample}$,
\begin{equation}
    \BasisF_{
        3 \,
        \ULine{\GenA_3 \LeafSimple \ULine{\GenA_1 \LeafSimple} \LeafSimple}
        \ULine{\GenA_2 \LeafSimple \LeafSimple} \LeafSimple
    }
    =
    \BasisE_{
        3 \,
        \ULine{\GenA_3 \LeafSimple \ULine{\GenA_1 \LeafSimple} \LeafSimple}
        \ULine{\GenA_2 \LeafSimple \LeafSimple} \LeafSimple
    }
    -
    \BasisE_{
        3 \, \ULine{\GenA_3 \LeafSimple \ULine{\GenA_1 \LeafSimple}
        \ULine{\GenA_2 \LeafSimple \LeafSimple}} \LeafSimple \LeafSimple
    }
    -
    \BasisE_{
        3 \, \ULine{\GenA_3 \ULine{\GenA_1 \LeafSimple} \LeafSimple \LeafSimple}
        \ULine{\GenA_2 \LeafSimple \LeafSimple} \LeafSimple
    }
    -
    \BasisE_{
        3 \, \ULine{\GenA_3 \ULine{\GenA_1 \LeafSimple} \LeafSimple
        \ULine{\GenA_2 \LeafSimple \LeafSimple}} \LeafSimple \LeafSimple
    }.
\end{equation}
By Möbius inversion and triangularity, for any leaning $\Signature$-forest~$\ForestF$,
\begin{equation}
    \BasisE_\ForestF =
    \sum_{\ForestF' \in \SetLeaningForests \App \Signature}
    \Iverson{\ForestF \EasterlyWindLeq \ForestF'} \; \BasisF_{\ForestF'},
\end{equation}
so that the set $\Bra{\BasisF_\ForestF : \ForestF \in \SetLeaningForests \App
\Signature}$ is a basis of $\NaturalHopfAlgebra \App \ULine{\SetTerms \App \Signature}$,
called the \Def{fundamental basis} (or \Def{$\BasisF$-basis} for short).

We can now state one of the main results of this work.

\begin{Statement}{Theorem}{thm:product_fundamental_basis}
    For any signature $\Signature$ and any leaning $\Signature$-forests $\ForestF_1$ and
    $\ForestF_2$,
    \begin{equation} \label{equ:product_fundamental_basis_1}
        \BasisF_{\ForestF_1} \Product \BasisF_{\ForestF_2}
        = \sum_{\ForestF \in \SetLeaningForests \App \Signature}
        \Iverson{
            \ForestF_1 \Over \ForestF_2 \EasterlyWindLeq \ForestF
            \EasterlyWindLeq \ForestF_1 \Under \ForestF_2
        }
        \; \BasisF_\ForestF.
    \end{equation}
\end{Statement}
\begin{Proof}
    Let $\Product'$ be the product on $\NaturalHopfAlgebra \App \ULine{\SetTerms \App
    \Signature}$ such that for any $\Signature$-forests $\ForestF_1$ and $\ForestF_2$,
    $\BasisF_{\ForestF_1} \Product' \BasisF_{\ForestF_2}$ is the right-hand side
    of~\eqref{equ:product_fundamental_basis_1}. For any leaning $\Signature$-forests
    $\ForestF_1$ and $\ForestF_2$, by denoting by $n_1$ the size of $\ForestF_1$ and by
    using Proposition~\ref{prop:leaning_forest_shuffle_interval} and
    Lemma~\ref{lem:over_under_restriction}, we have
    \begin{align}
        \BasisE_{\ForestF_1} \Product' \BasisE_{\ForestF_2}
        & =
        \sum_{\ForestF_1' \in \SetLeaningForests \App \Signature}
        \sum_{\ForestF_2' \in \SetLeaningForests \App \Signature}
        \Iverson{\ForestF_1 \EasterlyWindLeq \ForestF_1'}
        \Iverson{\ForestF_2 \EasterlyWindLeq \ForestF_2'}
        \;
        \BasisF_{\ForestF_1'} \Product' \BasisF_{\ForestF_2'}
        \\
        & =
        \sum_{\ForestF_1' \in \SetLeaningForests \App \Signature}
        \sum_{\ForestF_2' \in \SetLeaningForests \App \Signature}
        \sum_{\ForestF \in \SetLeaningForests \App \Signature}
        \Iverson{\ForestF_1 \EasterlyWindLeq \ForestF_1'}
        \Iverson{\ForestF_2 \EasterlyWindLeq \ForestF_2'}
        \Iverson{
            \ForestF_1' \Over \ForestF_2' \EasterlyWindLeq \ForestF
            \EasterlyWindLeq \ForestF_1' \Under \ForestF_2'
        }
        \;
        \BasisF_\ForestF
        \notag
        \\
        & =
        \sum_{\ForestF_1' \in \SetLeaningForests \App \Signature}
        \sum_{\ForestF_2' \in \SetLeaningForests \App \Signature}
        \Iverson{\ForestF_1 \EasterlyWindLeq \ForestF_1'}
        \Iverson{\ForestF_2 \EasterlyWindLeq \ForestF_2'}
        \sum_{\ForestF \in \ForestF_1' \cshuffle \ForestF_2'}
        \BasisF_\ForestF
        \notag
        \\
        & =
        \sum_{\ForestF_1' \in \SetLeaningForests \App \Signature}
        \sum_{\ForestF_2' \in \SetLeaningForests \App \Signature}
        \Iverson{\ForestF_1 \EasterlyWindLeq \ForestF_1'}
        \Iverson{\ForestF_2 \EasterlyWindLeq \ForestF_2'}
        \sum_{\ForestF \in \SetLeaningForests \App \Signature}
        \Iverson{\Top \App n_1 \App \ForestF = \ForestF_1'}
        \Iverson{\Bottom \App n_1 \App \ForestF = \ForestF_2'}
        \;
        \BasisF_\ForestF
        \notag
        \\
        & =
        \sum_{\ForestF \in \SetLeaningForests \App \Signature}
        \Iverson{\ForestF_1 \EasterlyWindLeq \Top \App n_1 \App \ForestF}
        \Iverson{\ForestF_2 \EasterlyWindLeq \Bottom \App n_1 \App \ForestF}
        \;
        \BasisF_\ForestF
        \notag
        \\
        & =
        \sum_{\ForestF \in \SetLeaningForests \App \Signature}
        \Iverson{\ForestF_1 \Over \ForestF_2 \EasterlyWindLeq \ForestF}
        \;
        \BasisF_\ForestF
        \notag
        \\
        & =
        \BasisE_{\ForestF_1 \Over \ForestF_2}.
        \notag
    \end{align}
    This shows that
    \begin{math}
        \BasisE_{\ForestF_1} \Product' \BasisE_{\ForestF_2}
        =
        \BasisE_{\ForestF_1} \Product \BasisE_{\ForestF_2},
    \end{math}
    so that $\Product'$ and $\Product$ are the same products. Therefore,
    \eqref{equ:product_fundamental_basis_1} holds.
\end{Proof}

The product of $\NaturalHopfAlgebra \App \ULine{\SetTerms \App \Signature}$ on the
$\BasisF$-basis is akin to a shuffle of reduced $\Signature$-forests. For instance, in
$\NaturalHopfAlgebra \App \ULine{\SetTerms \App \SignatureExample}$, we have
\begin{align}
    \BasisF_{
        \ColA{
            3 \,
            \ULine{\GenA_3 \LeafSimple \LeafSimple \LeafSimple}
            \ULine{\GenA_2 \LeafSimple \ULine{\GenA_2 \LeafSimple \LeafSimple}}
            \LeafSimple
        }
    }
    \Product
    \BasisF_{
        \ColB{
            4 \,
            \ULine{\GenA_3 \LeafSimple \ULine{\GenA_2 \LeafSimple \LeafSimple} \LeafSimple}
            \ULine{\GenA_1 \LeafSimple}
            \ULine{\GenA_2 \LeafSimple \LeafSimple}
            \LeafSimple
        }
    }
    & =
    \BasisF_{
        7 \,
        \ColA{
            \ULine{\GenA_3 \LeafSimple \LeafSimple \LeafSimple}
            \ULine{\GenA_2 \LeafSimple \ULine{\GenA_2 \LeafSimple \LeafSimple}}
        }
        \ColB{
            \ULine{\GenA_3 \LeafSimple \ULine{\GenA_2 \LeafSimple \LeafSimple} \LeafSimple}
            \ULine{\GenA_1 \LeafSimple}
            \ULine{\GenA_2 \LeafSimple \LeafSimple}
        }
        \LeafSimple
        \LeafSimple
    }
    \\
    & \quad +
    \BasisF_{
        7 \,
        \ColA{
            \ULine{\GenA_3 \LeafSimple \LeafSimple \LeafSimple}
            \ULine{
                \GenA_2 \LeafSimple \ULine{\GenA_2 \LeafSimple
                \ColB{\ULine{\GenA_3 \LeafSimple \ULine{\GenA_2 \LeafSimple \LeafSimple}
                \LeafSimple}}
            }}
        }
        \ColB{
            \ULine{\GenA_1 \LeafSimple}
            \ULine{\GenA_2 \LeafSimple \LeafSimple}
        }
        \LeafSimple
        \LeafSimple
        \LeafSimple
    }
    \notag
    \\
    & \quad +
    \BasisF_{
        7 \,
        \ColA{
            \ULine{\GenA_3 \LeafSimple \LeafSimple \LeafSimple}
            \ULine{
                \GenA_2 \LeafSimple \ULine{\GenA_2
                \ColB{\ULine{\GenA_3 \LeafSimple \ULine{\GenA_2 \LeafSimple \LeafSimple}
                \LeafSimple}}
                \LeafSimple
            }}
        }
        \ColB{
            \ULine{\GenA_1 \LeafSimple}
            \ULine{\GenA_2 \LeafSimple \LeafSimple}
        }
        \LeafSimple
        \LeafSimple
        \LeafSimple
    }
    \notag
    \\
    & \quad +
    \BasisF_{
        7 \,
        \ColA{
            \ULine{\GenA_3 \LeafSimple \LeafSimple \LeafSimple}
            \ULine{
                \GenA_2 \LeafSimple \ULine{\GenA_2
                \ColB{\ULine{\GenA_3 \LeafSimple \ULine{\GenA_2 \LeafSimple \LeafSimple}
                \LeafSimple}}
                \ColB{\ULine{\GenA_1 \LeafSimple}}
            }}
        }
        \ColB{
            \ULine{\GenA_2 \LeafSimple \LeafSimple}
        }
        \LeafSimple
        \LeafSimple
        \LeafSimple
        \LeafSimple
    }.
    \notag
\end{align}

Observe that the fundamental basis of $\NaturalHopfAlgebra \App \ULine{\SetTerms \App
\Signature}$ coincides with the ribbon basis of $\NCSF$~\cite{GKLLRT95}. Indeed, by
employing the notation introduced in Section~\ref{subsubsec:ncsf}, for any integer
compositions $\Par{r_1, \dots, r_k}$, $k \geq 0$, and $\Par{r'_1, \dots, r'_{k'}}$, $k' \geq
0$, we have
\begin{equation}
    \BasisF_{\Par{r_1, \dots, r_k}}
    \Product
    \BasisF_{\Par{r'_1, \dots, r'_{k'}}}
    =
    \BasisF_{\Par{r_1, \dots, r_k, r'_1, \dots, r'_{k'}}}
    +
    \BasisF_{\Par{r_1, \dots, r_k + r'_1, \dots, r'_{k'}}}.
\end{equation}
This product is the one of $\NCSF$ through its ribbon basis.

\subsubsection{Homogeneous basis}
Let us use again the leaning forest $\Signature$-easterly wind posets to build a new basis
of $\NaturalHopfAlgebra \App \ULine{\SetTerms \App \Signature}$. For any leaning
$\Signature$-forest $\ForestF$, let
\begin{equation}
    \BasisH_\ForestF :=
    \sum_{\ForestF' \in \SetLeaningForests \App \Signature}
    \Iverson{\ForestF' \EasterlyWindLeq \ForestF} \;
    \BasisF_{\ForestF'}.
\end{equation}
For instance, in $\NaturalHopfAlgebra \App \ULine{\SetTerms \App \SignatureExample}$,
\begin{equation}
    \BasisH_{
        3 \, \ULine{\GenA_3 \LeafSimple \ULine{\GenA_1 \LeafSimple} \LeafSimple}
        \ULine{\GenA_2 \LeafSimple \LeafSimple} \LeafSimple
    }
    =
    \BasisF_{
        3 \, \ULine{\GenA_3 \LeafSimple \LeafSimple \LeafSimple} \ULine{\GenA_1 \LeafSimple}
        \ULine{\GenA_2 \LeafSimple \LeafSimple}
    }
    +
    \BasisF_{
        3 \, \ULine{\GenA_3 \LeafSimple \LeafSimple \ULine{\GenA_1 \LeafSimple}}
        \ULine{\GenA_2 \LeafSimple \LeafSimple} \LeafSimple
    }
    +
    \BasisF_{
        3 \, \ULine{\GenA_3 \LeafSimple \ULine{\GenA_1 \LeafSimple} \LeafSimple}
        \ULine{\GenA_2 \LeafSimple \LeafSimple} \LeafSimple
    }.
\end{equation}
By Möbius inversion and triangularity, for any leaning $\Signature$-forest $\ForestF$,
\begin{equation}
    \BasisF_\ForestF =
    \sum_{\ForestF' \in \SetLeaningForests \App \Signature}
    \Iverson{\ForestF' \EasterlyWindLeq \ForestF} \;
    \ULine{\Mobius \App \Bra{1} \App \ForestF' \App \ForestF} \;
    \BasisH_{\ForestF'},
\end{equation}
so that the set $\Bra{\BasisH_\ForestF : \ForestF \in \SetLeaningForests \App \Signature}$
is a basis of $\NaturalHopfAlgebra \App \ULine{\SetTerms \App \Signature}$, called the
\Def{homogeneous basis} (or \Def{$\BasisH$-basis} for short).

We can now state one of the main results of this work.

\begin{Statement}{Theorem}{thm:product_homogeneous_basis}
    For any signature $\Signature$ and any leaning $\Signature$-forests $\ForestF_1$ and
    $\ForestF_2$,
    \begin{equation} \label{equ:product_homogeneous_basis_1}
        \BasisH_{\ForestF_1} \Product \BasisH_{\ForestF_2}
        = \BasisH_{\ForestF_1 \Under \ForestF_2}.
    \end{equation}
\end{Statement}
\begin{Proof}
    For any leaning $\Signature$-forests $\ForestF_1$ and $\ForestF_2$, by denoting by $n_1$
    the size of $\ForestF_1$ and by using Theorem~\ref{thm:product_fundamental_basis},
    Proposition~\ref{prop:leaning_forest_shuffle_interval}, and
    Lemma~\ref{lem:over_under_restriction}, we have
    \begin{align}
        \BasisH_{\ForestF_1} \Product \BasisH_{\ForestF_2}
        & =
        \sum_{\ForestF_1' \in \SetLeaningForests \App \Signature}
        \sum_{\ForestF_2' \in \SetLeaningForests \App \Signature}
        \Iverson{\ForestF_1' \EasterlyWindLeq \ForestF_1}
        \Iverson{\ForestF_2' \EasterlyWindLeq \ForestF_2}
        \;
        \BasisF_{\ForestF_1'} \Product \BasisF_{\ForestF_2'}
        \\
        & =
        \sum_{\ForestF_1' \in \SetLeaningForests \App \Signature}
        \sum_{\ForestF_2' \in \SetLeaningForests \App \Signature}
        \sum_{\ForestF \in \SetLeaningForests \App \Signature}
        \Iverson{\ForestF_1' \EasterlyWindLeq \ForestF_1}
        \Iverson{\ForestF_2' \EasterlyWindLeq \ForestF_2}
        \Iverson{
            \ForestF_1' \Over \ForestF_2' \EasterlyWindLeq \ForestF
            \EasterlyWindLeq \ForestF_1' \Under \ForestF_2'
        }
        \;
        \BasisF_\ForestF
        \notag
        \\
        & =
        \sum_{\ForestF_1' \in \SetLeaningForests \App \Signature}
        \sum_{\ForestF_2' \in \SetLeaningForests \App \Signature}
        \Iverson{\ForestF_1' \EasterlyWindLeq \ForestF_1}
        \Iverson{\ForestF_2' \EasterlyWindLeq \ForestF_2}
        \sum_{\ForestF \in \ForestF_1' \cshuffle \ForestF_2'}
        \BasisF_\ForestF
        \notag
        \\
        & =
        \sum_{\ForestF_1' \in \SetLeaningForests \App \Signature}
        \sum_{\ForestF_2' \in \SetLeaningForests \App \Signature}
        \Iverson{\ForestF_1' \EasterlyWindLeq \ForestF_1}
        \Iverson{\ForestF_2' \EasterlyWindLeq \ForestF_2}
        \sum_{\ForestF \in \SetLeaningForests \App \Signature}
        \Iverson{\Top \App n_1 \App \ForestF = \ForestF_1'}
        \Iverson{\Bottom \App n_1 \App \ForestF = \ForestF_2'}
        \;
        \BasisF_\ForestF
        \notag
        \\
        & =
        \sum_{\ForestF \in \SetLeaningForests \App \Signature}
        \Iverson{\Top \App n_1 \App \ForestF \EasterlyWindLeq \ForestF_1}
        \Iverson{\Bottom \App n_1 \App \ForestF \EasterlyWindLeq \ForestF_2}
        \;
        \BasisF_\ForestF
        \notag
        \\
        & =
        \sum_{\ForestF \in \SetLeaningForests \App \Signature}
        \Iverson{\ForestF \EasterlyWindLeq \ForestF_1 \Under \ForestF_2}
        \;
        \BasisF_\ForestF
        \notag
        \\
        & =
        \BasisH_{\ForestF_1 \Under \ForestF_2}.
        \notag
    \end{align}
    This shows that~\eqref{equ:product_homogeneous_basis_1} holds.
\end{Proof}

By employing the notation introduced at the end of
Section~\ref{subsubsec:hopf_algebras_forests}, the $\BasisH$-basis of $\NCSF$ admits the
following expression for its product. For any integer compositions $\Par{r_1,
\dots, r_k}$, $k \geq 0$, and $\Par{r'_1, \dots, r'_{k'}}$, $k' \geq 0$, we have
\begin{equation}
    \BasisH_{\Par{r_1, \dots, r_k}}
    \Product
    \BasisH_{\Par{r'_1, \dots, r'_{k'}}}
    =
    \BasisH_{\Par{r_1, \dots, r_k + r'_1, \dots, r'_{k'}}}.
\end{equation}

\section{Conclusion and open questions} \label{sec:conclusion}
We have introduced a new partial order relation $\EasterlyWindLeq$ on the underlying set
$\SetTerms \App \Signature$ of free operads, namely the easterly wind partial order. As
shown in this work, the resulting posets yield new bases of natural Hopf algebras
$\NaturalHopfAlgebra \App \ULine{\SetTerms \App \Signature}$ of free operads, sharing key
properties with a broad class of combinatorial Hopf algebras. We list here several open
questions and directions for future research in this context.

At the general level of the easterly wind posets, many properties remain unknown,
including an explicit expression for the Möbius function of $\SpecialNodes$-tilted
$\Signature$-easterly wind posets and the enumeration of the set $\TerminalInterval \App
\SpecialNodes \App \TermT$ of terms greater than or equal to the $\SpecialNodes$-tilted
$\Signature$-term~$\TermT$. This last question is linked with the enumeration of the
intervals of such posets.

In Section~\ref{subsubsec:rooted_tree_lattices}, we have shown that certain
$\SpecialNodes$-tilted $\Signature$-easterly wind posets contain, as maximal intervals, the
Tamari lattices. A natural question is whether one can similarly realize, as
$\SpecialNodes$-tilted $\Signature$-easterly wind posets, other classical lattices involving
treelike structures, such as the Kreweras lattices~\cite{Kre72}, the Stanley
lattices~\cite{Sta75,Knu04}, the $m$-Tamari lattices~\cite{BPR12}, the $m$-canyon
lattices~\cite{CG22}, and the pruning-grafting lattices~\cite{BP08}.

In Section~\ref{subsubsec:fuss_catalan_lattices}, we have defined lattices on the
combinatorial family of Fuss-Catalan objects. To the best of our knowledge, these lattices
are new and warrant a detailed combinatorial study, including the enumeration of their
intervals and their relationships with known structures on the same combinatorial family.

Besides, now linked with the natural Hopf algebra $\NaturalHopfAlgebra \App \ULine{\SetTerms
\App \Signature}$ of a free operad $\SetTerms \App \Signature$, one may ask how to express
the coproduct and the antipode of $\NaturalHopfAlgebra \App \ULine{\SetTerms \App
\Signature}$ on the new $\BasisF$-basis and $\BasisH$-basis. A further question concerns
using these two bases to investigate the cofreeness and self-duality of~$\NaturalHopfAlgebra
\App \ULine{\SetTerms \App \Signature}$, and obtain necessary and sufficient conditions for
these properties, depending on $\Signature$.

Finally, for now the easterly wind order is defined only at the level of terms. It would be
valuable to generalize this order on the underlying set of any operad $\Operad$, possibly
subject to certain restrictions, so that the analogous $\BasisF$-basis and $\BasisH$-basis
of $\NaturalHopfAlgebra \App \Operad$ satisfy generalizations of
Theorems~\ref{thm:product_fundamental_basis} and~\ref{thm:product_homogeneous_basis}. One
approach to build such a partial order relation on $\Operad$ is to choose a generating set
$\Signature_\Operad$ of $\Operad$ and consider the easterly wind order on treelike
factorizations of the elements of $\Operad$ as elements of the free operad $\SetTerms \App
\Signature_\Operad$. The main challenge is to choose a canonical factorization for each $x
\in \Operad$, since an element may admit, when $\Operad$ is not free, several
factorizations.


\bibliography{Bibliography}

\newcommand{\etalchar}[1]{$^{#1}$}
\begin{thebibliography}{GKL{\etalchar{+}}95}

\bibitem[BG16]{BG16}
J.-P. Bultel and S.~Giraudo.
\newblock {Combinatorial Hopf algebras from PROs}.
\newblock {\em Journal of Algebraic Combinatorics}, 44(2):455--493, 2016.

\bibitem[Bjo80]{Bjo80}
A.~Bjorner.
\newblock {Shellable and Cohen-Macaulay partially ordered sets}.
\newblock {\em Transactions of the American Mathematical Society},
  260(1):159--183, 1980.

\bibitem[BP08]{BP08}
J.~L. Baril and J.~M. Pallo.
\newblock {The pruning-grafting lattice of binary trees}.
\newblock {\em Theoretical Computer Science}, 409(3):382--393, 2008.

\bibitem[BPR12]{BPR12}
F.~Bergeron and L.-F. Preville-Ratelle.
\newblock {Higher trivariate diagonal harmonics via generalized Tamari posets}.
\newblock {\em Journal of Combinatorics}, 3(3):317--341, 2012.

\bibitem[BW96]{BW96}
A.~Bjorner and M.~L. Wachs.
\newblock {Shellable nonpure complexes and posets. I}.
\newblock {\em Transactions of the AMS}, 348(4):1299--1327, 1996.

\bibitem[CG22]{CG22}
C.~Combe and S.~Giraudo.
\newblock {Three Fuss-Catalan posets in interaction and their associative
  algebras}.
\newblock {\em Combinatorial Theory}, 2(1), 2022.

\bibitem[CGM15]{CGM15}
H.~Cheballah, S.~Giraudo, and R.~Maurice.
\newblock {Hopf algebra structure on packed square matrices}.
\newblock {\em Journal of Combinatorial Theory, Series A}, 133:139--182, 2015.

\bibitem[Com23]{Com23}
C.~Combe.
\newblock {Geometric realizations of {T}amari interval lattices via cubic
  coordinates}.
\newblock {\em Order}, 40(3):589--621, 2023.

\bibitem[Cra66]{Cra66}
H.~H. Crapo.
\newblock {The M\"obius function of a lattice}.
\newblock {\em Journal of Combinatorial Theory}, 1(1):126--131, 1966.

\bibitem[CS98]{CS98}
I.~Chajda and V.~Sn\'a{\v{s}}el.
\newblock {Congruences in Ordered Sets}.
\newblock {\em Mathematica Bohemica}, 123(1):95--100, 1998.

\bibitem[DHT02]{DHT02}
G.~Duchamp, F.~Hivert, and J.-Y. Thibon.
\newblock {Noncommutative Symmetric Functions VI: Free Quasi- Symmetric
  Functions and Related Algebras}.
\newblock {\em International Journal of Algebra and Computation}, 12:671--717,
  2002.

\bibitem[DP02]{DP02}
B.~A. Davey and H.~A. Priestley.
\newblock {\em {Introduction to Lattices and Order}}, pages 145--174.
\newblock Cambridge University Press, 2 edition, 2002.

\bibitem[FJN95]{FJN95}
R.~Freese, J.~Je{\v{z}}ek, and J.~B. Nation.
\newblock {\em {Free Lattices}}, volume~42 of {\em Mathematical Surveys and
  Monographs}.
\newblock American Mathematical Society, Providence, RI, 1995.

\bibitem[GC17]{CP17}
V.~Pilaud G.~Chatel.
\newblock Cambrian hopf algebras.
\newblock {\em Advances in Mathematics}, 311:598--633, 2017.

\bibitem[Gir11]{Gir11}
S.~Giraudo.
\newblock {\em {Combinatoire algébrique des arbres}}.
\newblock PhD thesis, Université de Marne-la-Vallée, 2011.

\bibitem[Gir12]{Gir12}
S.~Giraudo.
\newblock {Algebraic and combinatorial structures on pairs of twin binary
  trees}.
\newblock {\em Journal of Algebra}, 360:115--157, 2012.

\bibitem[Gir18]{Gir18}
S.~Giraudo.
\newblock {\em {Nonsymmetric Operads in Combinatorics}}.
\newblock Springer Nature Switzerland AG, 2018.
\newblock ix+172.

\bibitem[Gir24]{Gir24}
S.~Giraudo.
\newblock {Polynomial realizations of Hopf algebras built from nonsymmetric
  operads}.
\newblock {\em \Arxiv{2406.12559}}, 2024.

\bibitem[GKL{\etalchar{+}}95]{GKLLRT95}
I.~M. Gelfand, D.~Krob, A.~Lascoux, B.~Leclerc, V.~S. Retakh, and J.-Y. Thibon.
\newblock {Noncommutative Symmetric Functions}.
\newblock {\em Advances in Mathematics}, 112(2):218--348, 1995.

\bibitem[HNT05]{HNT05}
F.~Hivert, J.-C. Novelli, and J.-Y. Thibon.
\newblock {The Algebra of Binary Search Trees}.
\newblock {\em Theoretical Computer Science}, 339(1):129--165, 2005.

\bibitem[Knu04]{Knu04}
D.~Knuth.
\newblock {\em {The Art of Computer Programming. Volume 4, Fascicle 4.
  Generating all trees --- History of combinatorial generation}}.
\newblock Addison Wesley Longman, 2004.

\bibitem[Kre72]{Kre72}
G.~Kreweras.
\newblock {Sur les partitions non crois\'ees d'un cycle}.
\newblock {\em Discrete Mathematics}, 1(4):333--350, 1972.

\bibitem[LR98]{LR98}
J.-L. Loday and M.~Ronco.
\newblock {Hopf Algebra of the Planar Binary Trees}.
\newblock {\em Advances in Mathematics}, 139:293--309, 1998.

\bibitem[LR02]{LR02}
J.-L. Loday and M.~Ronco.
\newblock {Order Structure on the Algebra of Permutations and of Planar Binary
  Trees}.
\newblock {\em Journal of Algebraic Combinatorics}, 15(3):253--270, 2002.

\bibitem[ML14]{ML14}
M.~Mendez and J.~Liendo.
\newblock {An antipode formula for the natural Hopf algebra of a set operad}.
\newblock {\em Advances in Applied Mathematics}, 53:112--140, 2014.

\bibitem[MR95]{MR95}
C.~Malvenuto and C.~Reutenauer.
\newblock {Duality between Quasi-Symmetrical Functions and the Solomon Descent
  Algebra}.
\newblock {\em Journal of Algebra}, 177(3):967--982, 1995.

\bibitem[NT06]{NT06}
J.-C. Novelli and J.-Y. Thibon.
\newblock {Polynomial realizations of some trialgebras}.
\newblock {\em Formal Power Series and Algebraic Combinatorics}, 2006.

\bibitem[Rea04]{Rea04}
N.~Reading.
\newblock {Lattice Congruences of the Weak Order}.
\newblock {\em Order}, 21(4):315--344, 2004.

\bibitem[Sta75]{Sta75}
R.~P. Stanley.
\newblock {The Fibonacci lattice}.
\newblock {\em Fibonacci Quarterly}, 13(3):215--232, 1975.

\bibitem[Sta11]{Sta11}
R.~P. Stanley.
\newblock {\em {Enumerative Combinatorics}}, volume~1.
\newblock Cambridge University Press, second edition, 2011.

\bibitem[Tam62]{Tam62}
D.~Tamari.
\newblock {The algebra of bracketings and their enumeration}.
\newblock {\em Nieuw Archief voor Wiskunde}, 10(3):131--146, 1962.

\bibitem[vdL04]{vdl04}
P.~van~der Laan.
\newblock {\em {Operads. Hopf algebras and coloured Koszul duality}}.
\newblock PhD thesis, Universiteit Utrecht, 2004.

\end{thebibliography}

\end{document}